\documentclass[12pt]{amsart}
\usepackage{mathrsfs}
\usepackage{amsfonts}
\usepackage{latexsym,amsmath,amsthm,amssymb,amsfonts}

\numberwithin{equation}{section}
\allowdisplaybreaks

\def\endproof{$\hfill\Box$\\}
\def\s{\,\,\,\,}

\DeclareMathOperator{\tr}{tr}

\numberwithin{equation}{section}
\newtheorem{theorem}{Theorem}[section]
\newtheorem{lem}[theorem]{Lemma}
\newtheorem{thm}[theorem]{Theorem}

\newtheorem{rem}[theorem]{Remark}

\newtheorem{defnm}{Definition}[section]
\newcommand{\sch}{Schr\"odinger }
\newcommand{\n}{\nabla}
\newcommand{\p}{\partial}
\newcommand{\ii}{\mathbf{i}}

 \newcommand{\al}{\alpha}
 \newcommand{\be}{\beta}

 \newcommand{\Norm}[1]{\left\Vert#1\right\Vert}

\newcounter{Cnumber}

\title[ ]
{Local Nonautonomous Schr\"{o}dinger Flows on K\"{a}hler Manifolds}
\author[ ]
{Zonglin Jia,\quad\quad Youde Wang}
\date{}
\begin{document}

\maketitle

\begin{abstract}
$\,\,\,\,\,\,$In this paper, we prove that the nonautonomous Schr\"{o}dinger flow from a compact Riemannian manifold into a K\"ahler manifold admits a local solution. Under some certain conditions, the solution is unique and has higher regularity.
\end{abstract}

\section{Introduction}

Let $(M,g)$ be a closed Riemannian manifold with metric $g$ and $(N,J,h)$ be a complete K\"{a}hler manifold with complex structure $J$ and metric $h$.
In \cite{DW}, Ding and Wang introduced the concept of Schr\"{o}dinger flow for maps from $M$ into $N$. This flow is defined as follows:
\[\partial_tu(x,t)=J(u(x,t))\tau(u(x,t))\]
where $u$ is a smooth map from $M$ to $N$ and $\tau(u)$ is the tension field of $u$, which is defined by
\[\tau(u)=\tr(\nabla du).\]
Here $\n$ is the connection on the bundle $u^{*}TN\otimes T^*M$ which are the induced connection from Riemannian connections on $(M, g)$ and $(N,h)$. It is well known that
\[\tau(u)=0\]
if and only if $u$ is the critical point of
\[E(u):=\frac{1}{2}\int_M|du|^2\,dM.\]

Let $f:M\times[0,T_*]\rightarrow \mathbb{R}$ be $C^1$-smooth function with respect to the space variables. In the sequel we always suppose that there exist two constants $\eta\in[1,\infty)$ and $\delta\in(0,1]$(otherwise, one can replace $\eta$ by $\max\{\eta,1\}$ and $\delta$ by $\min\{\delta,1\}$) such that $\eta>f(x,t)>\delta$ on $M\times[0,T_*]$.
Similarly, we may define
\[E_f(u):=\frac{1}{2}\int_Mf(x, t)|du(x)|^2\,dM.\]
The corresponding tension field for every $t$ reads
\[\tau_f(u):=\tr(\nabla(fdu)).\]
In \cite{WW} (also see \cite{PWW3}), the following geometric flow was introduced
\begin{eqnarray}\label{Ai}
\left\{
\begin{array}{llll}
\partial_tu(x,t)=J(u(x,t))\tau_{f(x,t)}(u(x,t));\\
u(\cdot ,0)=u_0:M\rightarrow N.
\end{array}
\right.
\end{eqnarray}

This flow is of strong physical background and is called nonautonomous Schr\"{o}dinger flow with coupling function $f$ or NSF for short. In the case the coupling function $f(x,t) \equiv f(x)$ and $N$ is a two dimensional standard sphere $S^2$ in Euclidean space, the flow is just so-called inhomogeneous Heisenberg spin chain system, which also attracted one's attention for long time. Also, the study of inhomogeneous ferromagnets is of considerable intrinsic importance in higher dimensions in their own right(\cite{DPL, LD, LG, LG1, RS}). It is well-known that one dimensional Heisenberg ferromagnetic spin system in the classical continuum limit is found to be one of the interesting nonlinear systems which exhibit a very rich variety of dynamical properties. In particular, the simple chain model with isotropic bilinear exchange interaction is the pioneer among them whose elementary excitations are governed by solitons in addition to magnons (\cite{KIK}).

Now let us recall some previous results. On one hand, Heisenberg spin systems with different kinds of magnetic interactions such as anisotropy, interaction with external fields, Gilbert damping, etc., have been investigated by many authors. We refer to \cite{DPL, I, LD, SSB, WW, Y} and references therein.

On the other hand, as $f(x,t)\equiv 1$, (\ref{Ai}) is just the Schr\"odinger flow
\begin{eqnarray}\label{Ai1}
\left\{
\begin{array}{llll}
\partial_tu(x,t)=J(u(x,t))\tau(u(x,t));\\
u(\cdot ,0)=u_0:M\rightarrow N.
\end{array}
\right.
\end{eqnarray}
Ding and Wang discussed first in \cite{DW} the case that $M=S^1$ and $N$ is a complete K\"ahler manifold. They showed that, if the initial value $u_0\in C^{\infty}(S^1,N)$, the Schr\"odinger flow admits a unique smooth solution on $S^1\times[0,T)$ for some $T\in(0,\infty]$. Moreover, they proved that the solution can be extended globally under some assumptions on target manifolds. For the Schr\"odinger flow from $\mathbb{R}^m\times\mathbb{R}$ ($m=1,2$) into a Riemann surface, in \cite{CSU} Chang, Shatah and Uhlenbeck show that equation (\ref{Ai1}) has a unique global solution $u$ for every $u_0 \in H^1(\mathbb{R}^m)$. Moreover, if $u_0$ is smooth, then $u$ is smooth (also see \cite{PWW2}). In two-dimensional spaces ($m=2$), they assume $S^1$ symmetry of $N$ and consider symmetric maps and equivariant maps, and they obtain similar results for small energy data.

Later, Ding and Wang in \cite{DW1} proved that the Schr\"odinger flow from a closed Riemannian manifold into a complete K\"ahler manifold admits a local solution in a suitable Sobolev space. They also obtained local existence results for any Schr\"odinger flow from $\mathbb{R}^m$ into a complete K\"ahler manifold (also see \cite{McGahagan}). In particular, Rodnianski, Rubinstein and Staffilani considered the Schr\"odinger map flow from a one-dimensional domain into a complete K\"ahler manifold in \cite{RRS} and established the global well-posedness of the initial value problem for the Schr\"odinger map flow for maps from the real line into K\"ahler manifolds and for maps from the circle into Riemannian surfaces.

It is worthy to point out that Song and Wang improve the uniqueness results on the local solutions to the Schr\"odinger flows in \cite{DW, McGahagan} recently.\medskip

For the following Schr\"odinger flow (Landau-Lifshitz systems) from $\mathbb{R}^m\times\mathbb{R}$ into $S^2$
\begin{eqnarray}\label{Bi}
\left\{
\begin{array}{llll}
\partial_tu(x,t)=J(u(x,t))\tau(u(x,t)),\\
u(\cdot ,0)=u_0:\mathbb{R}^m\rightarrow S^2,
\end{array}
\right.
\end{eqnarray}
smooth solutions for (\ref{Bi}) were found by C. Sulem, J. P. Sulem, and C. Bardos \cite{SSB}. In particular, they established the local well-posedness of (\ref{Bi}) for initial data in $H^s(\mathbb{R}^m)$, $s > m/2+2$, $m\geqslant2$, by the difference method.

The first global well-posedness result for (\ref{Bi}) in critical spaces (precisely, global well-posedness for small data in the critical Besov spaces of dimensions $d\geqslant3$) was proved in \cite{BT}, and independently by Bejanaru in \cite{B1}. This was later improved to global regularity for small data in the critical Sobolev spaces of dimensions $d\geqslant4$ in \cite{BIK}. When the dimension is 2, in the case of equivariant data with energy close to the energy of the equivariant harmonic map, the existence of global equivariant solutions (and asymptotic stability) was proved in \cite{GKT}. In the case of radial or equivariant data of small energy, global well-posedness was proved in \cite{CSU}. Finally, in \cite{BIKT} the global well-posedness result for (\ref{Bi}) for small data in the critical Besov spaces of dimensions $d\geqslant2$ was addressed.

On the other hand, one also considered the dynamical behavior of the Schr\"odinger flows such as the stability in\cite{GKT, GS}. In particular, One discussed an energy critical Schr\"odinger map problem with a 2-sphere target for equivariant initial data of homotopy index $\kappa=1$. More precisely, in \cite{MRR} Merle, Rapha\"el and Rodnianski considered the Cauchy problem for the energy critical Schr\"odinger map from $\mathbb{R}^2\times\mathbb{R}$ into $S^2$ and gained some results on blowing-up (also see \cite{Pe}).\medskip

In \cite{PWW3}, the authors studied the so-called Inhomogeneous Schr\"{o}dinger Flow(ISF), which is a special case of NSF since in this case $f(x, t)$ does not depend on the variable $t$, i.e. $f(x,t)=\sigma(x)$.They showed that, if $M$ is a closed Riemannian manifold with $\dim(M)\leqslant3$ and $(N,J,h)$ is a K\"{a}hler manifold with nonpositive sectional curvature, and the coupling function $\sigma(x)>0$ (or $\sigma(x)<0$) is smooth enough, then ISF has a unique local solution in a suitable Sobolev space. In spite of these developments not much progress has been made in the case of higher-dimensional inhomogeneous Schr\"odinger flow.

In this paper, we will consider the local existence and uniqueness of NSF in a suitable Sobolev space. Before stating our main results, we need to introduce several definitions on the Sobolev spaces of sections with vector bundle values on $M$. Let $(E, M, \pi)$ be a vector bundle with base manifold $M$. If $(E, M, \pi)$ is equipped with a metric, then we may define so-called vector bundle value Sobolev spaces as follows:

\begin{defnm}
$H^{k,p}(M, E)$ is the completeness of the set of smooth sections with
compact supports denoted by $\{s| \,\,s\in C_0^\infty(M, E)\}$ with
respect to the norm
$$\|s\|^p_{H^{k,p}(M, E)}=\sum_{i=0}^k \int_M |\nabla^i s|^pdM.$$
Here $p\in[1,\infty)$ and $\nabla$ is the connection on $E$ which is compatible with the metric on $E$. For $p=\infty$, we also define
\[
||s||_{H^{k,\infty}(M,E)}:=\max\Big\{||\nabla^is||_{L^{\infty}(M,E)}\Big|0\leqslant i\leqslant k\Big\}
\]
For the sake of convenience, we usually denote $H^{k,2}$ by $H^k$.
\end{defnm}

When $(E, M, \pi)$ is a trivial bundle, i.e. $E=M\times \mathbb{R}^L$, We usually use $W^{k,p}(M, \mathbb{R}^L)$ to denote the space of Sobolev functions.

By Nash's embedding theorem, for large enough $L$, $(N, h)$ can always be embedded isometrically in $\mathbb{R}^L$. This means that $(N, h)$ may be regarded as a submanifold of $\mathbb{R}^L$ with original point $O$ in $N$ so that we are able to use theory of Sobolev spaces.

\begin{defnm} \label{def:2}Let $\mathbb{N}^+$ be the set of positive integers. For $k \in\mathbb{N}^+\cup\{0\}$ and $p\in[1,\infty]$, the Sobolev space of
maps from $(M, g)$ into a Riemannian manifold $(N, h)$ is defined by
$$W^{k,p}(M,N) = \{u\in W^{k,p}(M, \mathbb{R}^L)\Big|\,\, u(x)\in N\,\,\,\mbox{for a. e.}\,\,\, x\in M \}.$$
And we also define
$$\dot{W}^{k,p}(M,N) = \{u\Big|||\nabla^ku||_{L^p}<\infty,\,\, u(x)\in N\,\,\,\mbox{for a. e.}\,\,\, x\in M \}.$$
\end{defnm}

The main result of this paper is stated as follows:
\begin{thm}\label{NSFtheorem1.1}
Let $(M, g)$ be a $m$-dimensional closed Riemannian manifold and $(N, J, h)$ be a complete K\"ahler manifold. Assume that $u_0\in W^{r,2}(M,N)$ and $f\in C^1([0,T_*],C^{r+1}(M))$ and $0<\delta<f(x,t)<\eta<\infty$, where $r\geqslant m_0+2$ and $m_0:=[\frac{m}{2}]+1$. Here $[q]$ is the integral part of $q$. Then the Cauchy problem of NSF admits a local strong solution $u\in L^{\infty}([0,T],W^{r,2}(M,N))$ for some $T=T(N,||u_0||_{W^{m_0+2,2}})$. When $r\geqslant m_0+3$, the local solution is unique. Furthermore, if $u_0\in C^{\infty}(M,N)$ and $f\in C^{\infty}(M\times[0,T_*])$, the local solution $u\in C^{\infty}(M\times[0,T],N)$.
\end{thm}
\begin{thm}\label{NSFtheorem1.2}
Let $\mathbb{R}^m$ be a $m$-dimensional Euclidean space and $(N,J,h)$ be a compact K\"ahler manifold. Assume that $f\in C^1([0,T_*],C^{r+1}(\mathbb{R}^m))$ where $r\geqslant m_0+2$ is an integer number and $m_0:=[\frac{m}{2}]+1$. If there exist $R_1>0,\cdots, R_m>0$, such that for every $j=1, 2, \cdots, m$,
\[0<f(x^1,\cdots,x^j,\cdots,x^m,t)=f(x^1,\cdots,x^j+R_j,\cdots,x^m,t),\]
then (\ref{Ai}) with $u_0\in W^{r,2}(\mathbb{R}^m,N)$ admits a local strong solution
$$u\in L^{\infty}([0,T_{m_0}],W^{r,2}(\mathbb{R}^m,N))$$
for some $T_{m_0}>0$. If $r\geqslant m_0+3$, the local solution is unique. Furthermore, if $u_0\in\mathscr{H}:=\bigcap\limits_{k=1}^{\infty}W^{k,2}(\mathbb{R}^m,N)$ and $f\in C^1([0,T_*],C^{\infty}(\mathbb{R}^m))$, then $u\in C^{\infty}([0,T_{m_0}],\mathscr{H})$.
\end{thm}

\begin{thm}\label{thm0.1}
Suppose $M$ is an $m$-dimensional complete manifold with bounded Ricci curvature $Ric_M$, $N$ is a complete K\"{a}hler manifold with bounded geometry. Let $\mathscr{S}_{\infty}:=W^{2,2}(M,N)\cap\dot{W}^{1,\infty}(M,N)\cap\dot{W}^{2,\infty}(M,N)$. We assume that $f\in L^{\infty}([0,T],W^{1,\infty}(M))$, $f(t,x)>\delta>0$ and $\partial_tf\in L^{\infty}([0,T]\times M)$. If $u_1, u_2\in L^{\infty}([0,T],\mathscr{S}_{\infty})$ are two solutions to NSF with the same initial value $u_0\in\mathscr{S}_{\infty}$, then $u_1=u_2$ a.e. on $[0,T]\times M$.
\end{thm}

\begin{thm}\label{thm0.2}
Suppose $m\geqslant3$, $M$ is an $m$-dimensional complete manifold with bounded Riemannian curvature $R^M$ and positive injective radius $inj(M)>0$. $N$ is a  complete K\"{a}hler manifold with bounded geometry. Let $\mathscr{S}_m:=W^{[\frac{m}{2}]+1,2}(M,N)\cap\dot{W}^{1,\infty}(M,N)\cap\dot{W}^{2,m}(M,N)$. We assume that $f\in L^{\infty}([0,T],W^{1,2}(M)\cap \dot{W}^{1,m}(M))$, $0<\delta<f(t,x)<\eta<\infty$ and $\partial_tf\in L^{\infty}([0,T]\times M)$. If $u_1,u_2\in L^{\infty}([0,T],\mathscr{S}_m)$ are two solutions to NSF with the same initial value $u_0\in\mathscr{S}_m$, then $u_1=u_2$ a.e. on $[0,T]\times M$.
\end{thm}

We will follow the ideas and techniques from \cite{DW1, PWW3, WW} to approach the problem by employing the following approximate system:
\[
\left\{
\begin{array}{llll}
\partial_t u=\varepsilon\tau_f(u)+J(u)\tau_f(u),\\
u(\cdot ,0)=u_0\in C^{\infty}(M,N)
\end{array}
\right.
\]
with $f\in C^{\infty}(M\times[0,T_*])$. For each $\varepsilon\in(0,1)$, the above parabolic system admits a $C^{\infty}$-smooth local solution $u_{\varepsilon}$ and the existence interval of the solution is denoted by $[0,T_{\varepsilon})$. For the details we refer to the appendix of this paper.

The key ingredient of the proof of theorem \ref{NSFtheorem1.1} is how to introduce a norm on the Sobolev space $H^{k,2}(M, u_{\varepsilon}^{*}TN)$ such that we can obtain some uniform a priori estimates on $u_{\varepsilon}$ with respect to $\varepsilon$. In other words, we try to obtain some a priori estimates on $u_{\varepsilon}$ which do not depend on $\varepsilon$. By some complicate computations we find the following norm
\[
||\nabla u_{\varepsilon}||^2_{H^{k,2}(f)}:=\sum\limits_{l=1}^{k+1}\int_{M}f^l|\nabla^l u_{\varepsilon}|^2\,dM
\]
is a suitable choice. We can see easily
\begin{equation}\label{1.4}
\delta^{k+1}||\nabla u_{\varepsilon}||^2_{H^{k,2}}\leqslant||\nabla u_{\varepsilon}||^2_{H^{k,2}(f)}\leqslant\eta^{k+1}||\nabla u_{\varepsilon}||^2_{H^{k,2}}.
\end{equation}

As for Theorem \ref{NSFtheorem1.2}, we follow the idea of Theorem 1.2 of \cite{DW1} to approximate $\mathbb{R}^m$ by $\widetilde{\mathbb{T}}_k^m:=(\mathbb{R}^1/2kR_1\mathbb{Z})\times\cdots\times(\mathbb{R}^1/2kR_m\mathbb{Z})$. Since any map $u_0$ which is in $W^{r,2}(\mathbb{R}^m,N)$ can be approximated by a family of smooth maps $\{u_{i0}\}\subseteq C_0^{\infty}(\mathbb{R}^m,N)$ and $u_{i0}$ can be regarded as a map from $\widetilde{\mathbb{T}}_{k_i}^m$ into $N$ for some large $k_i$, we employ a system on $\widetilde{\mathbb{T}}_{k_i}^m\times[0,T_i)$ taking $u_{i0}$ as its initial value.

The key idea to prove Theorem \ref{thm0.1} and Theorem \ref{thm0.2} is to construct an energy functional as follows:
\[
E(t)=\int_Md^2(u_1(t),u_2(t))\,dM+\frac{1}{2}\int_Mf(t)\cdot|\mathcal{P}\nabla_2u_2(t)-\nabla_1u_1(t)|^2\,dM,
\]
where $d(y_1,y_2)$ means the distance of $y_1$ and $y_2$ on $N$. $\nabla_{\lambda}$ denotes the connection on $u_{\lambda}^*TN$ induced by the Levi-Civita connection $\nabla^N$ on $N$. And $\mathcal{P}$ is a global isomorphism between $u_1^*TN\otimes T^*M$ and $u_2^*TN\otimes T^*M$. It is defined by parallel transportation in $N$. Our goal is to show that this functional satisfies a Gronwall type inequality and hence vanishes on $[0,T]$.
\section{Preliminary}
\textbf{In the next we appoint that the same index appearing twice means to sum it.}
\subsection{\sch Flow in Moving Frame}

Let $T>0$ and $I=[0,T]$ be an interval. Suppose $u:M\times I\to N$ is a solution to the \sch flow
\begin{equation}\label{e:sch}
  \p_tu=J(u)\tau(u).
\end{equation}
We are going to rewrite the above equation in a moving frame, namely, a chosen gauge of the pull-back bundle $u^*TN$.

To fix our notations, we let roman numbers $i,j,k$ be indices ranging from $1$ to $m$ and Greek letters $\al,\be$ ranging from $1$ to $n$, where $n$ is the dimension of $N$. Let $M\times I$ be endowed with the natural product metric. We will use $\n$ to denote connections on different vector bundles which are naturally induced by the Levi-Civita connections on $M$ and $N$. In particular, this includes the pull-back bundle $u^*TN$ on $M\times I$, the pull-back bundle $u(t)^*TN$ on some time slice $M\times \{t\}$ for $t\in I$ and their tensor product bundles with the cotangent bundle $T^*M$. Sometimes in the context, we also use more specific notations such as $\n^N$ and $\n^M$ to emphasize which connection we are using.

Locally on an open geodesic ball $U\subset M$, we may choose an orthonormal frame $\{e_i\}_{i=1}^m$ of the tangent bundle $TM$. Set $e_0:=\p_t$ such that $\{e_i\}_{i=0}^m$ forms a local orthonormal basis of $T(U\times I)$. For convenience, we denote $\n_i:=\n_{e_i}$ and $\n_t=\n_0$. Then $\n_t e_i=0, 0\leqslant i\leqslant m$ with $\n$ the Levi-Civita connection on $M\times I$.

Recall that the tension field is
$$\tau(u)=\tr_g \n^2u = \n_k\n_k u,$$
where $\n_k u$ denotes the covariant derivative of $u$ and is a section of the bundle $u^*TN\otimes T^*M$. Then the \sch flow (\ref{e:sch}) has the form
\begin{equation*}
  \n_tu=J(u)\n_k\n_ku.
\end{equation*}
Differentiating the equation, we get
\begin{equation*}
  \begin{aligned}
    \n_t\n_i u&=\n_i\n_t u\\
    &=\n_i(J(u)\n_k\n_k u)=J(u)\n_i\n_k\n_ku\\
    &=J(u)(\n_k\n_i\n_k u + R^N(\n_i u, \n_k u)\n_k u + R^M(e_i, e_k, e_k, e_l)\n_l u)\\
    &=J(u)(\n_k\n_k\n_i u + R^N(\n_i u, \n_k u)\n_k u + Ric^M(e_i, e_l)\n_l u),
  \end{aligned}
\end{equation*}
where $R^M, R^N$ are the curvature of $M$ and $N$, respectively, and $Ric^M$ is the Ricci curvature of $M$. Here we have used the fact that $\n$ is torsion free and $\n^N J=0$ since the target manifold $N$ is K\"ahler.

Next we choose a local frame $\{f_\al\}_{\al=1}^{n}$ of the pull-back bundle $u^*TN$, such that the complex structure $J$ in this frame is reduced to a constant skew-symmetric matrix which we denote by $J_0$. Letting $\n_\ii u=:\phi_\ii^\al f_\al$, we may further rewrite the above equation for $\n_i u$ as
\begin{equation*}\label{e:sch1}
  \n_t\phi_i=J_0(\Delta_x\phi_i+R^N(\phi_i, \phi_k)\phi_k + Ric^M_{ij}\phi_j),
\end{equation*}
where $\Delta_x=\n_k\n_k$ is the Laplacian operator on $u(t)^*TN\otimes T^*M$. Obviously, this is a nonlinear Schr\"odinger system.

\subsection{Propositions on Sobolev Norms}
We also need to recall an important theorem proved in \cite{DW1}. This is a generalized Gagliardo-Nirenberg inequality.

\begin{thm}\label{thm:int} (\cite{DW1})
Suppose $s \in C^\infty(E)$ is a section where $E$ is a vector
bundle on $M$. Then we have
\begin{equation}\label{NSF2.2}
\Norm{\n^j s}_{L^p} \leq C\Norm{s}^a_{H^{k,q}}\Norm{s}^{1-a}_{L^r},
\end{equation}
where $1\leq p,q,r\leq \infty$, and $j/k \leq a \leq 1 (j/k \leq a <
1$ if $q=m/(k-j) \neq 1)$ are numbers such that
\[ \frac 1 p = \frac j m + \frac 1 r + a\left(\frac 1 q - \frac 1 r - \frac k m\right). \]
The constant $C$ only depends on $M$ and the numbers $j,k,q,r,a$.
\end{thm}

Hence, for $\n u \in \Gamma(u^*(TN))$, by the above theorem we have
\begin{equation}\label{equ:int}
\Norm{\n^{j+1} u}_{L^p} \leq C\Norm{\n u}^a_{H^{k,q}}\Norm{\n u}^{1-a}_{L^r}.
\end{equation}
Ding and Wang also showed that the $H^{k,p}$ norm of section $\n u$ is equivalent to the normal Sobolev $W^{k+1, p}$ norm of the map $u$. Precisely, we have

\begin{lem}\label{lem:equ} (\cite{DW1})
Assume that $k > m/2$. Then there exists a constant $C = C(N,k)$
such that for all $u \in C^\infty(M,N)$,
\[ \Norm{D u}_{W^{k-1,2}} \leq C\sum_{i=1}^k \Norm{\n u}^i_{H^{k-1,2}}, \]
and
\[ \Norm{\n u}_{H^{k-1,2}} \leq C\sum_{i=1}^k \Norm{D u}^i_{W^{k-1,2}}. \]
Here $D$ denotes the usually derivative.
\end{lem}

In order to prove Theorem \ref{NSFtheorem1.2}, we need the following lemma. It is almost the same as Proposition 2.1 of \cite{DW1}, so we omit its proof.

\begin{lem}\label{NSFlemma2.3}(\cite{DW1})
For given $R_1>0,\cdots, R_m>0$, let
\[M=\widetilde{\mathbb{T}}_k^m:=(\mathbb{R}^1/2kR_1\mathbb{Z})\times\cdots\times(\mathbb{R}^1/2kR_m\mathbb{Z}),\]
where $k\geqslant1$ and the Riemannian metric of $\widetilde{\mathbb{T}}_k^m$ is just the Euclidean metric. Then the constant in (\ref{equ:int}) does not depend on the diameter $k$(That is to say, it depends only upon the geometry of $\widetilde{\mathbb{T}}_1^m$ and $j,p,a,k,q,r$).
\end{lem}
\section{The Proofs of Theorems}
\textbf{Proof of theorem \ref{NSFtheorem1.1}.} It is known that, for each $\varepsilon\in(0,1)$, the approximate system
\begin{equation}\label{T:1}
\left\{
\begin{array}{llll}
\partial_t u=\varepsilon\tau_f(u)+J(u)\tau_f(u),\\
u(\cdot ,0)=u_0\in C^{\infty}(M,N)
\end{array}
\right.
\end{equation}
with $f\in C^{\infty}(M\times[0,T_*])$ is a uniformly parabolic system. By the classical theory(the details will be written in the appendix), there exists a $T_{\varepsilon}>0$ such that (\ref{T:1}) admits a smooth local solution $u_{\varepsilon} \in C(M\times[0,T_{\varepsilon}),N)\cap  C^{\infty}(M\times(0,T_{\varepsilon}),N)$.

To prove the main theorem we need to estimate the uniform upper bound of $||\nabla u_{\varepsilon}||^2_{H^{k,2}(f)}$ and the uniform lower bound of $T_{\varepsilon}$ with respect to $\varepsilon$. For this purpose, we try to use the interpolation inequality for Sobolev sections on vector bundle to obtain a nonlinear differential inequality with respect to $||\nabla u_{\varepsilon}||^2_{H^{m_0,2}(f)}$, where
$$m_0:=\left[\frac{m}{2}\right]+1$$
and $[q]$ is the integral part of $q$. In the sequel, we also denote $\{q\}:=q-[q]$.

Once we obtain the uniform estimates on the quantity $||\nabla u_{\varepsilon}||^2_{H^{m_0,2}(f)}$ with respect to $T$, we can infer a linear differential inequality about $||\nabla u_{\varepsilon}||^2_{H^{k,2}(f)}$ for $k\geqslant m_0+1$(its coefficients can be expressed nonlinearly by $||\nabla u_{\varepsilon}||_{H^{k-1,2}(f)}$). Using the comparison theorem of ODE iteratively, we can get the uniform lower bound $T$ of $T_{\varepsilon}$, which is independent of $\varepsilon$ and $k$, and the upper bound of $||\nabla u_{\varepsilon}||^2_{H^{k,2}(f)}$ not depending on $\varepsilon$. Therefore, to prove the main theorem we need to establish several lemmas as follows.
\medskip

For the sake of simplicity, in the sequel we always assume that $M=\mathbb{T}^m$ be an $m$-dimensional flat torus without loss of generality. Before stating the proof, let us define following notations. We denote $u_{\varepsilon}$ by $u$ and let
\[\Omega:=\{y\in N|dist_{N}(y,u_0(M))<1\}\]
so that $\overline{\Omega}$ is compact. Let
\[T'_{\varepsilon}:=\sup\{t|u(M\times[0,t])\subseteq\Omega\}.\]
For convenience, we denote $\nabla_{\frac{\partial}{\partial x_{a_i}}}$ by $\nabla_{a_i}$. Setting $\overrightarrow{a}=(a_1,\cdots,a_l)$, we denote $\nabla_{a_1}\cdots\nabla_{a_l}$ by $\nabla_{\overrightarrow{a}}$. Denote $||\cdot||_{L^{p}}$ by $||\cdot||_p$. Let $|\overrightarrow{a}|=l$ be the length of $\overrightarrow{a}$ where $l=\sum\limits_ia_i$. For a set $A$, let $|A|$ denote the number of elements in $A$.

\begin{lem}\label{NSFlemma3.1}
Let $u$ be the solution of $(\ref{T:1})$ and $E_f(u)$ is defined in the introduction. Then for any $T>0$, $E_f(u)$ is bounded in $[0,T]$. More precisely,
\[
E_f(u)\leqslant E_{f_0}(u_0)\exp\Big\{\int_0^t\bar{C}_1(s)ds\Big\},
\]
for some continuous function $\bar{C}_1(t)$ which does not depend on $\varepsilon$ and for all $t\in[0,T]$, where $f_0:=f(\cdot,0)$.
\end{lem}

\textbf{Proof:} We compute
\[
\aligned
\frac{d}{dt}E_f(u)&=\frac{1}{2}\int_Mf_t|\nabla u|^2\,dM+\int_Mf\langle\nabla_t\nabla_au,\nabla_au\rangle\,dM\\
&=\frac{1}{2}\int_M\frac{f_t}{f}\cdot f|\nabla u|^2\,dM+\int_M\langle\nabla_a\nabla_tu,f\nabla_au\rangle\,dM\\
&=\frac{1}{2}\int_M\frac{f_t}{f}\cdot f|\nabla u|^2\,dM-\int_M\langle\nabla_tu,\tau_f(u)\rangle\,dM\\
&=\frac{1}{2}\int_M\frac{f_t}{f}\cdot f|\nabla u|^2\,dM-\int_M\langle(\varepsilon+J(u))\tau_f(u),\tau_f(u)\rangle\,dM\\
&=\frac{1}{2}\int_M\frac{f_t}{f}\cdot f|\nabla u|^2\,dM-\varepsilon\int_M|\tau_f(u)|^2\\
&\leqslant \bar{C}_1(t)\cdot E_f(u).
\endaligned
\]
Here $$\bar{C}_1(t):=\max\limits_{x\in M}\left\{\left|\frac{f_t(x,t)}{f(x,t)}\right|\right\}$$ is a continuous function and we have used the antisymmetry of $J$. By Gronwall inequality, we get the required inequality.\endproof

Obviously, we have that
\[E_f(u)=||\nabla u||^2_{H^{0,2}(f)}.\]
For $||\nabla u||^2_{H^{k-1,2}(f)}$ with $k\geqslant2$, we want to get a similar estimate. So, for $2\leqslant l\leqslant k$, we consider:
\begin{equation}\label{T:2}
\aligned
&\frac{1}{2}\frac{d}{dt}\int_Mf^l\langle\nabla_{a_1}\cdots\nabla_{a_l}u,\nabla_{a_1}\cdots\nabla_{a_l}u\rangle\,dM\\
=&\frac{1}{2}\int_Mlf^{l-1}f_t|\nabla^lu|^2\,dM+\int_Mf^l\langle\nabla_t\nabla_{a_1}\cdots\nabla_{a_l}u,\nabla_{a_1}\cdots\nabla_{a_l}u\rangle\,dM.
\endaligned
\end{equation}
Exchanging the order of covariant differentiation we have
\begin{equation}\label{T:3}
\nabla_t\nabla_{a_1}\cdots\nabla_{a_l}u=\nabla_{a_1}\cdots\nabla_{a_l}\nabla_tu+Q_1(u),
\end{equation}
where
\[Q_1(u):=\sum(\nabla^{r_1}R)(\nabla_{\overrightarrow{b^1_1}}u,\cdots,\nabla_{\overrightarrow{b^1_{r_1}}}u)(\nabla_{\overrightarrow{c^1}}\nabla_tu,\nabla_{\overrightarrow{d^1}}u)\nabla_{\overrightarrow{e^1}}\nabla_{a_l}u\]
and $\overrightarrow{b^1_1},\, \cdots,\, \overrightarrow{b^1_{r_1}},\,\overrightarrow{c^1},\,\overrightarrow{d^1},\,\overrightarrow{e^1}$ are multi-indexes satisfying that $(\overrightarrow{b^1_1},\cdots,\overrightarrow{b^1_{r_1}},\overrightarrow{c^1},\overrightarrow{d^1},\overrightarrow{e^1})$ is a permutation of $(a_1,\cdots,a_{l-1})$. We should remind that, in the above, $|\overrightarrow{d^1}|>0$, $|\overrightarrow{b^1_1}|,\cdots,|\overrightarrow{b^1_{r_1}}|$, $|\overrightarrow{c^1}|$, $|\overrightarrow{e^1}|$ are nonnegative integers.
\medskip

For any $t\in[0,T'_{\varepsilon})$, there holds
\[
|Q_1(u)|\leqslant c^*_1(l,\Omega)\sum|\nabla^{j_1}u|\cdots|\nabla^{j_s}u|\cdot|\nabla^{j_{s+1}}\nabla_tu|
\]
where $s\geqslant2$ and $j_1+j_2+\cdots+j_s+j_{s+1}=l$ with $1\leqslant j_i\leqslant l-1$ for $1\leqslant i\leqslant s$ and $0\leqslant j_{s+1}\leqslant l-2$.
\medskip

Taking $j_{s+1}$ order covariant derivative with respect to both sides of $(\ref{T:1})$ and using the integrability of the complex structure $J$ of K\"ahler manifold $N$, we have
\[\nabla^{j_{s+1}}\nabla_tu=(\varepsilon+J(u))\nabla^{j_{s+1}}\tau_f(u).\]
Taking inner product on both sides of the above equation and by the antisymmetry of $J$, one can obtain
\[|\nabla^{j_{s+1}}\nabla_tu|^2=(1+\varepsilon^2)|\nabla^{j_{s+1}}\tau_f(u)|^2.\]
So
\[|\nabla^{j_{s+1}}\nabla_tu|=\sqrt{1+\varepsilon^2}|\nabla^{j_{s+1}}\nabla_{\beta}(f\nabla_{\beta}u)|\leqslant c(j_{s+1})\sum\limits_{p+q=j_{s+1}+1}|\nabla^pf|\cdot|\nabla^{q+1}u|.\]
Letting
\[c_1(l,\Omega):=c_1^*(l,\Omega)\cdot c(j_{s+1}),\]
we get
\begin{equation}\label{T:4}
|Q_1(u)|\leqslant c_1(l,\Omega)\sum|\nabla^{j_1}u|\cdots|\nabla^{j_s}u|\cdot|\nabla^{q+1}u|\cdot|\nabla^pf|
\end{equation}
with
\[s\geqslant2,\,\,\,\,j_1+\cdots+j_s+p+q=l+1,\]
and
\[1\leqslant j_i\leqslant l-1\,\,\,\,\mbox{for}\,\,\,\,1\leqslant i\leqslant s, \,\,\,\,1\leqslant p+q\leqslant l-1,\,\,0\leqslant p\leqslant l-1.\]
So
\begin{equation}\label{T*:5}
\aligned
A:=&\int_Mf^l\langle\nabla_t\nabla_{a_1}\cdots\nabla_{a_l}u,\nabla_{a_1}\cdots\nabla_{a_l}u\rangle\,dM\\
=&\int_M\langle f^l\nabla_{a_1}\cdots\nabla_{a_l}u,(\varepsilon+J(u))\nabla_{a_1}\cdots\nabla_{a_l}\nabla_{\beta}(f\nabla_{\beta}u)\rangle\,dM\\
&+\int_M\langle Q_1(u),f^l\nabla_{a_1}\cdots\nabla_{a_l}u\rangle\,dM.
\endaligned
\end{equation}
Exchanging the order of covariant differentiation we have
\begin{equation}\label{T:5}
\nabla_{a_1}\cdots\nabla_{a_l}\nabla_{\beta}(f\nabla_{\beta}u)=\nabla_{\beta}\nabla_{a_1}\cdots\nabla_{a_l}(f\nabla_{\beta}u)+Q_2(u),
\end{equation}
where
\[
Q_2(u):=\sum(\nabla^{r_2}R)(\nabla_{\overrightarrow{b^2_1}}u,\cdots,\nabla_{\overrightarrow{b^2_{r_2}}}u)(\nabla_{\overrightarrow{c^2}}u,\nabla_{\overrightarrow{d^2}}\nabla_{\beta}u)\nabla_{\overrightarrow{e^2}}(f\nabla_{\beta}u)
\]
and $\overrightarrow{b^2_1},\cdots,\overrightarrow{b^2_{r_2}},\overrightarrow{c^2},\overrightarrow{d^2},\overrightarrow{e^2}$ are multi-indexes satisfying that $(\overrightarrow{b^2_1},\cdots,\overrightarrow{b^2_{r_2}},\overrightarrow{c^2},\overrightarrow{d^2},\overrightarrow{e^2})$ is a permutation of $(a_1,\cdots,a_{l})$. Here, $|\overrightarrow{c^2}|>0$, $|\overrightarrow{b^2_1}|,\cdots,|\overrightarrow{b^2_{r_2}}|,|\overrightarrow{d^2}|$ and $|\overrightarrow{e^2}|$ are nonnegative integers.
\begin{equation}\label{T:6}
|Q_2(u)|\leqslant c_2(l,\Omega)\sum|\nabla^{j_1}u|\cdots|\nabla^{j_s}u|\cdot|\nabla^{q+1}u|\cdot|\nabla^pf|
\end{equation}
with
\[s\geqslant2,\s j_1+\cdots+j_s+p+q=l+1, \s 1\leqslant j_i\leqslant l\s \mbox{for}\s 1\leqslant i\leqslant s,\]
\[ \s 0\leqslant p+q\leqslant l-1\s\mbox{and}\s 0\leqslant p\leqslant l-1.\]
In view of (\ref{T:3}) and (\ref{T:5}), we infer from (\ref{T*:5}) by integrating by parts
\begin{equation}\label{T:7}
\aligned
A=&\int_Mf^l\langle\nabla_t\nabla_{a_1}\cdots\nabla_{a_l}u,\nabla_{a_1}\cdots\nabla_{a_l}u\rangle\,dM\\
=&\int_M\langle f^l\nabla_{a_1}\cdots\nabla_{a_l}u,\,(\varepsilon+J(u))\nabla_{\beta}\nabla_{a_1}\cdots\nabla_{a_l}(f\nabla_{\beta}u)\rangle\,dM\\
&+\int_M\langle Q_1(u)+(\varepsilon+J(u))Q_2(u),f^l\nabla_{a_1}\cdots\nabla_{a_l}u\rangle\,dM\\
=&-\int_M\langle(\varepsilon+J(u))\nabla_{a_1}\cdots\nabla_{a_l}(f\nabla_{\beta}u),\,f^l\nabla_{\beta}\nabla_{a_1}\cdots\nabla_{a_l}u\rangle\,dM\\
&-\int_M\langle(\varepsilon+J(u))\nabla_{a_1}\cdots\nabla_{a_l}(f\nabla_{\beta}u),\,lf^{l-1}\nabla_{\beta} f\cdot\nabla_{a_1}\cdots\nabla_{a_l}u\rangle\,dM\\
&+\int_M\langle Q_1(u)+(\varepsilon+J(u))Q_2(u),\,f^l\nabla_{a_1}\cdots\nabla_{a_l}u\rangle\,dM.
\endaligned
\end{equation}

By exchanging the order of covariant differentiation again, we have
\begin{equation}\label{T:8}
\nabla_{\beta}\nabla_{a_1}\cdots\nabla_{a_l}u=\nabla_{a_1}\cdots\nabla_{a_l}\nabla_{\beta}u+Q_3(u),
\end{equation}
where
\[Q_3(u):=\sum(\nabla^{r_3}R)(\nabla_{\overrightarrow{b^3_1}}u,\cdots,\nabla_{\overrightarrow{b^3_{r_3}}}u)(\nabla_{\overrightarrow{c^3}}\nabla_{\beta}u,\nabla_{\overrightarrow{d^3}}u)\nabla_{\overrightarrow{e^3}}\nabla_{a_l}u\]
and $\overrightarrow{b^3_1},\cdots,\overrightarrow{b^3_{r_3}},\overrightarrow{c^3},\overrightarrow{d^3},\overrightarrow{e^3}$ are multi-indexes satisfying that $(\overrightarrow{b^3_1},\cdots,\overrightarrow{b^3_{r_3}},\overrightarrow{c^3},\overrightarrow{d^3},\overrightarrow{e^3})$ is just a permutation of $(a_1,\cdots,a_{l-1})$. Here, $|\overrightarrow{b^3_1}|,\cdots,|\overrightarrow{b^3_{r_3}}|,|\overrightarrow{c^3}|$ and $|\overrightarrow{e^3}|$ are nonnegative integers and $|\overrightarrow{d^3}|>0$. Moreover, we have that
\begin{equation}\label{Q:3}
|Q_3(u)|\leqslant c_3(l,\Omega)\sum|\nabla^{j_1}u|\cdots|\nabla^{j_s}u|,
\end{equation}
where
\[s\geqslant3,\,\,\,\,j_1+\cdots+j_s=l+1,\,\,\,\,1\leqslant j_i\leqslant l-1.\]
It is also easy to see that
\begin{equation}\label{T:11}
|\nabla Q_3(u)|\leqslant c_4(l,\Omega)\sum|\nabla^{j_1}u|\cdots|\nabla^{j_s}u|
\end{equation}
with \[s\geqslant3,\s j_1+\cdots+j_s=l+2,\s 1\leqslant j_i\leqslant l.\]

Instituting (\ref{T:8}) into (\ref{T:7}), we obtain
\begin{equation}\label{T:9}
\aligned
A=&\int_Mf^l\langle\nabla_t\nabla_{a_1}\cdots\nabla_{a_l}u,\,\nabla_{a_1}\cdots\nabla_{a_l}u\rangle\,dM\\
=&-\int_M\langle(\varepsilon+J(u))\nabla_{a_1}\cdots\nabla_{a_l}(f\nabla_{\beta}u),\, f^l\nabla_{a_1}\cdots\nabla_{a_l}\nabla_{\beta}u\rangle\,dM\\
&-\int_M\langle(\varepsilon+J(u))\nabla_{a_1}\cdots\nabla_{a_l}(f\nabla_{\beta}u),\,lf^{l-1}\nabla_{\beta}f\cdot\nabla_{a_1}\cdots\nabla_{a_l}u\rangle\,dM\\
&+\int_M\langle(\varepsilon+J(u))\nabla_{a_2}\cdots\nabla_{a_l}(f\nabla_{\beta}u),\, f^l\nabla_{a_1}Q_3(u)+lf^{l-1}\nabla_{a_1}f\cdot Q_3(u)\rangle\,dM\\
&+\int_M\langle Q_1(u)+(\varepsilon+J(u))Q_2(u),\,f^l\nabla_{a_1}\cdots\nabla_{a_l}u\rangle\,dM.
\endaligned
\end{equation}

Now we consider the sum of the first term and the second term on the right hand side of (\ref{T:9}), denoted by $S$. Then, we have
\begin{equation}\label{T:12}
\aligned
S=&-\int_M\langle(\varepsilon+J(u))\nabla_{a_1}\cdots\nabla_{a_l}(f\nabla_{\beta}u),\\
&\,\,\,\,\,\,\,\,\,\,\,\,\,\,\,\,\,f^l\nabla_{a_1}\cdots\nabla_{a_l}\nabla_{\beta}u+lf^{l-1}\nabla_{\beta}f\cdot\nabla_{a_1}\cdots\nabla_{a_l}u\rangle\,dM\\
=&-\int_M\langle(\varepsilon+J(u))(f\nabla_{a_1}\cdots\nabla_{a_l}\nabla_{\beta}u+\sum\limits_{j=1}^l\nabla_{a_j}f\cdot\nabla_{a_1}\cdots\widehat{\nabla}_{a_j}\cdots\nabla_{a_l}\nabla_{\beta}u\\
&+\sum\limits_{i<j}\nabla_{a_i}\nabla_{a_j}f\cdot\nabla_{a_1}\cdots\widehat{\nabla}_{a_i}\cdots\widehat{\nabla}_{a_j}\cdots\nabla_{a_l}\nabla_{\beta}u+
\cdots+\nabla_{a_1}\cdots\nabla_{a_l}f\cdot\nabla_{\beta}u),\\
&\s f^l\nabla_{a_1}\cdots\nabla_{a_l}\nabla_{\beta}u+lf^{l-1}\nabla_{\beta}f\cdot\nabla_{a_1}\cdots\nabla_{a_l}u\rangle\,dM,
\endaligned
\end{equation}
where $\widehat{\nabla}_{a_j}$ means deleting $\nabla_{a_j}$.

Set
\[
Q_4(u):=\sum\limits_{i<j}\nabla_{a_i}\nabla_{a_j}f\cdot\nabla_{a_1}\cdots\widehat{\nabla}_{a_i}\cdots\widehat{\nabla}_{a_j}\cdots\nabla_{a_l}\nabla_{\beta}u+\cdots+\nabla_{a_1}\cdots\nabla_{a_l}f\cdot\nabla_{\beta}u.
\]
It follows
\begin{equation}\label{S:1}
\aligned
S=&-\int_M\langle(\varepsilon+J(u))(f\nabla_{a_1}\cdots\nabla_{a_l}\nabla_{\beta}u+\sum\limits_{j=1}^l\nabla_{a_j}f\cdot\nabla_{a_1}\cdots\widehat{\nabla}_{a_j}\cdots\nabla_{a_l}
\nabla_{\beta}u\\
&+Q_4(u)),\, f^l\nabla_{a_1}\cdots\nabla_{a_l}\nabla_{\beta}u+lf^{l-1}\nabla_{\beta}f\cdot\nabla_{a_1}\cdots\nabla_{a_l}u\rangle\,dM\\
=&-\int_M\langle(\varepsilon+J(u))(f\nabla_{a_1}\cdots\nabla_{a_l}\nabla_{\beta}u+\sum\limits_{j=1}^l\nabla_{a_j}f\cdot\nabla_{a_1}\cdots\widehat{\nabla}_{a_j}\cdots\nabla_{a_l}
\nabla_{\beta}u,\\
& f^l\nabla_{a_1}\cdots\nabla_{a_l}\nabla_{\beta}u+lf^{l-1}\nabla_{\beta}f\cdot\nabla_{a_1}\cdots\nabla_{a_l}u\rangle\,dM\\
&-\int_M\langle(\varepsilon+J(u))Q_4(u), \, f^l\nabla_{a_1}\cdots\nabla_{a_l}\nabla_{\beta}u+lf^{l-1}\nabla_{\beta}f\cdot\nabla_{a_1}\cdots\nabla_{a_l}u\rangle\,dM\\
=&\, J_1 + J_2.
\endaligned
\end{equation}
Then, by taking integration by parts, we have
\begin{equation}\label{T:13}
\aligned
J_2=&-\int_M\langle(\varepsilon+J(u))Q_4(u),\, f^l\nabla_{a_1}\cdots\nabla_{a_l}\nabla_{\beta}u+lf^{l-1}\nabla_{\beta}f\cdot\nabla_{a_1}\cdots\nabla_{a_l}u\rangle\\
=&\int_M\langle(\varepsilon+J(u))(f^l\nabla_{a_1}Q_4(u)+lf^{l-1}\nabla_{a_1}f\cdot Q_4(u)),\, \nabla_{a_2}\cdots\nabla_{a_l}\nabla_{\beta}u\rangle\,dM\\
 & - \int_M\langle(\varepsilon+J(u))Q_4(u),\, lf^{l-1}\nabla_{\beta}f\cdot\nabla_{a_1}\cdots\nabla_{a_l}u\rangle\,dM.
\endaligned
\end{equation}
Note that there hold
\begin{equation}\label{Q:4}
|Q_4(u)|\leqslant c_5(l,\Omega)\sum\limits_{p+q=l+1,1\leqslant q\leqslant l-1}|\nabla^pf|\cdot|\nabla^qu|
\end{equation}
and
\begin{equation}\label{T:15}
|\nabla Q_4(u)|\leqslant c_6(l,\Omega)\sum\limits_{p+q=l+2,1\leqslant q\leqslant l}|\nabla^pf|\cdot|\nabla^qu|.
\end{equation}

By antisymmetry of $J$ we get the following
\begin{equation*}
\aligned
J_1=&-\int_M\langle(\varepsilon+J(u))(f\nabla_{a_1}\cdots\nabla_{a_l}\nabla_{\beta}u+\sum\limits_{j=1}^l\nabla_{a_j}f\cdot\nabla_{a_1}\cdots\widehat{\nabla}_{a_j}\cdots\nabla_{a_l}\nabla_{\beta}u),\\
&\s f^l\nabla_{a_1}\cdots\nabla_{a_l}\nabla_{\beta}u+lf^{l-1}\nabla_{\beta}f\cdot\nabla_{a_1}\cdots\nabla_{a_l}u\rangle\,dM\\
=&-\varepsilon\int_Mf^{l+1}|\nabla^{l+1}u|^2\,dM\\
&-\sum\limits_{j=1}^l\int_M\langle(\varepsilon+J(u))\nabla_{a_j}f\cdot\nabla_{a_1}\cdots\widehat{\nabla}_{a_j}\cdots\nabla_{a_l}\nabla_{\beta}u,\, f^l\nabla_{a_1}\cdots\nabla_{a_l}\nabla_{\beta}u\rangle\,dM\\
&-\int_M\langle(\varepsilon+J(u))(f\nabla_{a_1}\cdots\nabla_{a_l}\nabla_{\beta}u),\, lf^{l-1}\nabla_{\beta}f\cdot\nabla_{a_1}\cdots\nabla_{a_l}u\rangle\,dM\\
&-\int_M\langle(\varepsilon+J(u))\sum\limits_{j=1}^l\nabla_{a_j}f\cdot\nabla_{a_1}\cdots\widehat{\nabla}_{a_j}\cdots\nabla_{a_l}\nabla_{\beta}u,\, lf^{l-1}\nabla_{\beta}f\cdot\nabla_{a_1}\cdots\nabla_{a_l}u\rangle.
\endaligned
\end{equation*}
This equality can also be written as
\begin{equation}\label{T:16}
\aligned
J_1=&-\varepsilon\int_Mf^{l+1}|\nabla^{l+1}u|^2\,dM -J_1^1-J_1^2\\
&-\int_M\langle(\varepsilon+J(u))(f\nabla_{a_1}\cdots\nabla_{a_l}\nabla_{\beta}u),\, lf^{l-1}\nabla_{\beta}f\cdot\nabla_{a_1}\cdots\nabla_{a_l}u\rangle\,dM.
\endaligned
\end{equation}
Here,
$$J_1^1=\sum\limits_{j=1}^l\int_M\langle(\varepsilon+J(u))\nabla_{\beta}f\cdot\nabla_{a_1}\cdots\widehat{\nabla}_{a_j}\cdots\nabla_{a_l}\nabla_{a_j}u, \, f^l\nabla_{a_1}\cdots\nabla_{a_{j-1}}\nabla_{\beta}\nabla_{a_{j+1}}\cdots\nabla_{a_l}\nabla_{a_j}u\rangle$$
and
$$J_1^2=\int_M\langle(\varepsilon+J(u))\sum\limits_{j=1}^l\nabla_{a_j}f\cdot\nabla_{a_1}\cdots\widehat{\nabla}_{a_j}\cdots\nabla_{a_l}\nabla_{\beta}u, \, lf^{l-1}\nabla_{\beta}f\cdot\nabla_{a_1}\cdots\nabla_{a_l}u\rangle\,dM.$$

Now we compute $J_1^1$. By exchanging the order of covariant differentiation, we have
\begin{equation}\label{T:17}
J_1^1=\int_M\langle(\varepsilon+J(u))\nabla_{\beta}f\cdot(\nabla_{a_1}\cdots\nabla_{a_l}u+Q_5^j(u)),\, f^l(\nabla_{a_1}\cdots\nabla_{a_l}\nabla_{\beta}u+Q_6^j(u))\rangle,
\end{equation}
where
\begin{equation}\label{Q:5}
|Q_5^j(u)|\leqslant c_7(l,\Omega)\sum|\nabla^{j_1}u|\cdots|\nabla^{j_s}u|
\end{equation}
with
\[s\geqslant3,\,\,\,\,j_1+\cdots+j_s=l,\,\,\,\,1\leqslant j_i\leqslant l-2;\]
and
\begin{equation}\label{Q:6}
|Q_6^j(u)|\leqslant c_9(l,\Omega)\sum|\nabla^{j_1}u|\cdots|\nabla^{j_s}u|
\end{equation}
with
\[s\geqslant3,\,\,\,\,j_1+\cdots+j_s=l+1,\,\,\,\,1\leqslant j_i\leqslant l-1.\]
Later, we will also use the following
\begin{equation}\label{3.22}
|\nabla Q_5^j(u)|\leqslant c_8(l,\Omega)\sum|\nabla^{j_1}u|\cdots|\nabla^{j_s}u|
\end{equation}
with
\[s\geqslant3,\,\,\,\,j_1+\cdots+j_s=l+1,\,\,\,\,1\leqslant j_i\leqslant l-1.\]

Substituting (\ref{T:17}) into (\ref{T:16}) and using the antisymmetry of $J$ we obtain the following identity
\[\aligned
J_1=&-\varepsilon\int_Mf^{l+1}|\nabla^{l+1}u|^2\,dM\\
&-\sum\limits_{j=1}^l\int_M\langle(\varepsilon+J(u))\nabla_{\beta}f(\nabla_{a_1}\cdots\nabla_{a_l}u+Q_5^j(u)),\,f^l(\nabla_{a_1}\cdots\nabla_{a_l}\nabla_{\beta}u+Q^j_6(u))\rangle\,dM\\
&-\int_M\langle(\varepsilon+J(u))(f\nabla_{a_1}\cdots\nabla_{a_l}\nabla_{\beta}u),\, lf^{l-1}\nabla_{\beta}f\cdot\nabla_{a_1}\cdots\nabla_{a_l}u\rangle\,dM\\
&-\int_M\langle(\varepsilon+J(u))\sum\limits_{j=1}^l\nabla_{a_j}f\cdot\nabla_{a_1}\cdots\widehat{\nabla}_{a_j}\cdots\nabla_{a_l}\nabla_{\beta}u,\,lf^{l-1}\nabla_{\beta}f\cdot
\nabla_{a_1}\cdots\nabla_{a_l}u\rangle\,dM.
\endaligned\]
Furthermore, we write $J_1$ as
\begin{eqnarray*}
J_1&=&-\varepsilon\int_Mf^{l+1}|\nabla^{l+1}u|^2\,dM\\
&&-l\int_M\langle(\varepsilon+J(u))\nabla_{\beta}f\cdot\nabla_{a_1}\cdots\nabla_{a_l}u,\, f^l\nabla_{a_1}\cdots\nabla_{a_l}\nabla_{\beta}u\rangle\,dM\\
&&-\sum\limits_{j=1}^l\int_M\langle(\varepsilon+J(u))\nabla_{\beta}f\cdot Q_5^j(u),\, f^l\nabla_{a_1}\cdots\nabla_{a_l}\nabla_{\beta}u \rangle\,dM\\
&&-\sum\limits_{j=1}^l\int_M\langle(\varepsilon+J(u))\nabla_{\beta}f\cdot\nabla_{a_1}\cdots\nabla_{a_l}u,\,f^l\cdot Q_6^j(u)\rangle\,dM\\
&&-\sum\limits_{j=1}^l\int_M\langle(\varepsilon+J(u))\nabla_{\beta}f\cdot Q_5^j(u),\, f^l\cdot Q_6^j(u)\rangle\,dM\\
&&-l\int_M\langle(\varepsilon-J(u))\nabla_{\beta} f\cdot\nabla_{a_1}\cdots\nabla_{a_l}u,f^l\nabla_{a_1}\cdots\nabla_{a_l}\nabla_{\beta}u\rangle\,dM\\
&&-\int_M\langle(\varepsilon+J(u))\sum\limits_{j=1}^l\nabla_{a_j}f\cdot\nabla_{a_1}\cdots\widehat{\nabla}_{a_j}\cdots\nabla_{a_l}\nabla_{\beta}u,\, lf^{l-1}\nabla_{\beta}f\cdot\nabla_{a_1}\cdots\nabla_{a_l}u\rangle\,dM.
\end{eqnarray*}
By rearranging the terms in the above equality we have
\begin{eqnarray*}
J_1&=&-\varepsilon\int_Mf^{l+1}|\nabla^{l+1}u|^2\,dM-2\varepsilon l\int_M\langle\nabla_{\beta}f\cdot\nabla_{a_1}\cdots\nabla_{a_l}u,\,f^l\nabla_{a_1}\cdots\nabla_{a_l}\nabla_{\beta}u\rangle\,dM\\
&&+\sum\limits_{j=1}^l\int_M\langle(\varepsilon+J(u))(\nabla_{a_1}f\cdot lf^{l-1}\nabla_{\beta}f\cdot Q_5^j(u)+f^l\nabla_{a_1}\nabla_{\beta}f\cdot Q_5^j(u)\\
&&+f^l\nabla_{\beta}f\cdot\nabla_{a_1}Q_5^j(u)),\, \nabla_{a_2}\cdots\nabla_{a_l}\nabla_{\beta}u\rangle\,dM\\
&&-\sum_{j=1}^l\int_M\langle(\varepsilon+J(u))\nabla_{\beta}f\cdot\nabla_{a_1}\cdots\nabla_{a_l}u,\,f^lQ_6^j(u)\rangle\,dM\\
&&-\sum_{j=1}^l\int_M\langle(\varepsilon+J(u))\nabla_{\beta}f\cdot Q_5^j(u),\,f^lQ_6^j(u)\rangle\,dM\\
&&-\int_M\langle(\varepsilon+J(u))\sum\limits_{j=1}^l\nabla_{a_j}f\cdot\nabla_{a_1}\cdots\widehat{\nabla}_{a_j}\cdots\nabla_{a_l}\nabla_{\beta}u,\,
lf^{l-1}\nabla_{\beta}f\cdot\nabla_{a_1}\cdots\nabla_{a_l}u\rangle\,dM.
\end{eqnarray*}
Substituting the above expression of $J_1$ and the expression (\ref{T:13}) of $J_2$ into $(\ref{T:12})$ one can obtain an expression of $S$ as follows
\begin{eqnarray*}
S&=&-\varepsilon\int_Mf^{l+1}|\nabla^{l+1}u|^2\,dM-2\varepsilon l\int_M\langle\nabla_{\beta}f\cdot\nabla_{a_1}\cdots\nabla_{a_l}u,\,f^l\nabla_{a_1}\cdots\nabla_{a_l}\nabla_{\beta}u\rangle\,dM\\
& &+\sum\limits_{j=1}^l\int_M\langle(\varepsilon+J(u))(\nabla_{a_1}f\cdot lf^{l-1}\nabla_{\beta}f\cdot Q_5^j(u)+f^l\nabla_{a_1}\nabla_{\beta}f\cdot Q_5^j(u)\\
& &+f^l\nabla_{\beta}f\cdot\nabla_{a_1}Q_5^j(u)),\, \nabla_{a_2}\cdots\nabla_{a_l}\nabla_{\beta}u\rangle\,dM\\
& &-\sum_{j=1}^l\int_M\langle(\varepsilon+J(u))\nabla_{\beta}f\cdot\nabla_{a_1}\cdots\nabla_{a_l}u,\,f^lQ_6^j(u)\rangle\,dM\\
& &-\sum_{j=1}^l\int_M\langle(\varepsilon+J(u))\nabla_{\beta}f\cdot Q_5^j(u),\,f^lQ_6^j(u)\rangle\,dM\\
& &-\int_M\langle(\varepsilon+J(u))\sum\limits_{j=1}^l\nabla_{a_j}f\cdot\nabla_{a_1}\cdots\widehat{\nabla}_{a_j}\cdots\nabla_{a_l}\nabla_{\beta}u,\,
lf^{l-1}\nabla_{\beta}f\cdot\nabla_{a_1}\cdots\nabla_{a_l}u\rangle\,dM\\
& &+\int_M\langle(\varepsilon+J(u))(f^l\nabla_{a_1}Q_4(u)+lf^{l-1}\nabla_{a_1}f\cdot Q_4(u)),\, \nabla_{a_2}\cdots\nabla_{a_l}\nabla_{\beta}u\rangle\,dM\\
& &- \int_M\langle(\varepsilon+J(u))Q_4(u),\, lf^{l-1}\nabla_{\beta}f\cdot\nabla_{a_1}\cdots\nabla_{a_l}u\rangle\,dM.
\end{eqnarray*}

On the other hand, (\ref{T:9}) can be written as
\[\aligned
A=&\int_Mf^l\langle\nabla_t\nabla_{a_1}\cdots\nabla_{a_l}u,\,\nabla_{a_1}\cdots\nabla_{a_l}u\rangle\,dM\\
=&S+\int_M\langle(\varepsilon+J(u))\nabla_{a_2}\cdots\nabla_{a_l}(f\nabla_{\beta}u),\, f^l\nabla_{a_1}Q_3(u)+lf^{l-1}\nabla_{a_1}f\cdot Q_3(u)\rangle\,dM\\
&+\int_M\langle Q_1(u)+(\varepsilon+J(u))Q_2(u),\, f^l\nabla_{a_1}\cdots\nabla_{a_l}u\rangle\,dM.
\endaligned\]
Then, substituting the expression of $S$ into the above identity we derive the following
\begin{eqnarray*}
A&=&\int_Mf^l\langle\nabla_t\nabla_{a_1}\cdots\nabla_{a_l}u,\nabla_{a_1}\cdots\nabla_{a_l}u\rangle\,dM\\
&=&-\varepsilon\int_Mf^{l+1}|\nabla^{l+1}u|^2\,dM-2\varepsilon l\int_M\langle\nabla_{\beta}f\cdot\nabla_{a_1}\cdots\nabla_{a_l}u,\, f^l\nabla_{a_1}\cdots\nabla_{a_l}\nabla_{\beta}u\rangle\,dM\\
& &+\int_M\langle Q_1(u)+(\varepsilon+J(u))Q_2(u),\, f^l\nabla_{a_1}\cdots\nabla_{a_l}u\rangle\,dM\\
& &+\int_M\langle(\varepsilon+J(u))\nabla_{a_2}\cdots\nabla_{a_l}(f\nabla_{\beta}u),\,f^l\nabla_{a_1}Q_3(u)+lf^{l-1}\nabla_{a_1}f\cdot Q_3(u)\rangle\,dM\\
& &-\int_M\langle(\varepsilon+J(u))Q_4(u),\, lf^{l-1}\nabla_{\beta}f\cdot\nabla_{a_1}\cdots\nabla_{a_l}u\rangle\,dM\\
& &+\int_M\langle(\varepsilon+J(u))(f^l\nabla_{a_1}Q_4(u)+lf^{l-1}\nabla_{a_1}f\cdot Q_4(u)),\, \nabla_{a_2}\cdots\nabla_{a_l}\nabla_{\beta}u\rangle\,dM\\
& &+\sum\limits_{j=1}^l\int_M\langle(\varepsilon+J(u))(\nabla_{a_1}f\cdot lf^{l-1}\nabla_{\beta}f\cdot Q_5^j(u)+f^l\nabla_{a_1}\nabla_{\beta}f\cdot Q_5^j(u)\\
& &+f^l\nabla_{\beta}f\cdot\nabla_{a_1}Q_5^j(u)),\, \nabla_{a_2}\cdots\nabla_{a_l}\nabla_{\beta}u\rangle\,dM\\
& &-\sum\limits_{j=1}^l\int_M\langle(\varepsilon+J(u))\nabla_{\beta}f\cdot Q_5^j(u),\,f^lQ_6^j(u)\rangle\,dM\\
& &-\sum\limits_{j=1}^l\int_M\langle(\varepsilon+J(u))\nabla_{\beta}f\cdot \nabla_{a_1}\cdots\nabla_{a_l}u,\,f^lQ_6^j(u)\rangle\,dM\\
& &-\int_M\langle(\varepsilon+J(u))\sum\limits_{j=1}^l\nabla_{a_j}f\cdot\nabla_{a_1}\cdots\widehat{\nabla}_{a_j}\cdots\nabla_{a_l}\nabla_{\beta}u,\,
lf^{l-1}\nabla_{\beta}f\cdot\nabla_{a_1}\cdots\nabla_{a_l}u\rangle\,dM.
\end{eqnarray*}
Consequently, substituting the above identity into (\ref{T:2}), i.e. the following
\[
\aligned
&\frac{1}{2}\frac{d}{dt}\int_Mf^l\langle\nabla_{a_1}\cdots\nabla_{a_l}u,\,\nabla_{a_1}\cdots\nabla_{a_l}u\rangle\,dM\\
=&\frac{1}{2}\int_Mlf^{l-1}f_t|\nabla^lu|^2\,dM+\int_Mf^l\langle\nabla_t\nabla_{a_1}\cdots\nabla_{a_l}u,\,\nabla_{a_1}\cdots\nabla_{a_l}u\rangle\,dM,
\endaligned
\]
and then using the following inequality
\[
\aligned
&-\varepsilon\int_Mf^{l+1}|\nabla^{l+1}u|^2\,dM-2\varepsilon l\int_M\langle f^{\frac{l-1}{2}}\nabla_{\beta}f\cdot\nabla_{a_1}\cdots\nabla_{a_l}u,f^{\frac{l+1}{2}}\nabla_{a_1}\cdots\nabla_{a_l}\nabla_{\beta}u\rangle\,dM\\
\leqslant&-\varepsilon\int_Mf^{l+1}|\nabla^{l+1}u|^2\,dM+2\varepsilon l\int_Mf^{\frac{l-1}{2}}|\nabla f|\cdot|\nabla^lu|\cdot f^{\frac{l+1}{2}}|\nabla^{l+1}u|\,dM\\
\leqslant&-\varepsilon\int_Mf^{l+1}|\nabla^{l+1}u|^2\,dM+\varepsilon\int_Mf^{l+1}|\nabla^{l+1}u|^2\,dM+\varepsilon l^2\int_Mf^{l-1}|\nabla f|^2\cdot|\nabla^l u|^2\,dM\\
=&\,\varepsilon l^2\int_Mf^{l-1}|\nabla f|^2\cdot|\nabla^l u|^2\,dM,
\endaligned
\]
We obtain
\begin{eqnarray*}
& &\frac{1}{2}\frac{d}{dt}\int_Mf^l|\nabla^lu|^2\,dM\\
&\leqslant&\frac{l}{2}\int_M\frac{|f_t|}{f}f^l|\nabla^lu|^2\,dM+\varepsilon l^2\int_M\frac{|\nabla f|^2}{f}f^l|\nabla^lu|^2\,dM\\
&&+c(l,\Omega)\int_Mf^l|\nabla^lu|(|Q_1(u)|+|Q_2(u)|)\,dM\\
&&+c(l,\Omega)\int_M\sum\limits_{p+q=l-1}|\nabla^pf|\cdot|\nabla^{q+1}u|(f^l|\nabla Q_3(u)|+f^{l-1}|\nabla f|\cdot|Q_3(u)|)\,dM\\
&&+c(l,\Omega)\int_M|Q_4(u)|f^{l-1}|\nabla f|\cdot|\nabla^lu|\,dM\\
&&+c(l,\Omega)\int_M(f^l|\nabla Q_4(u)|+f^{l-1}|\nabla f|\cdot|Q_4(u)|)|\nabla^lu|\,dM\\
&&+c(l,\Omega)\sum\limits_{j=1}^l\int_M(|\nabla f|^2f^{l-1}|Q_5^j(u)|+f^l|\nabla^2f|\cdot|Q_5^j(u)|+f^l|\nabla f|\cdot|\nabla Q_5^j(u)|)|\nabla^lu|\,dM\\
&&+c(l,\Omega)\sum\limits_{j=1}^l\int_M|\nabla f|\cdot f^l\cdot|Q_5^j(u)|\cdot|Q_6^j(u)|\,dM\\
&&+c(l,\Omega)\sum\limits_{j=1}^l\int_M|\nabla f|\cdot f^l\cdot|\nabla^lu|\cdot|Q_6^j(u)|\,dM\\
&&+c(l,\Omega)\int_M|\nabla f|^2|\nabla^lu|^2f^{l-1}\,dM.
\end{eqnarray*}
Here $$c(l,\Omega):=l\cdot\max\{1,\max\limits_{1\leqslant i\leqslant9}\{c_i(l,\Omega)\},\,\sup\limits_{y\in\Omega}(1+|J(y)|)\}.$$
Therefore, by taking a complicate computation and using H\"older inequality we infer
\begin{eqnarray}\label{ZL}
&&\frac{1}{2}\frac{d}{dt}\int_Mf^l|\nabla^lu|^2\,dM\nonumber\\
&\leqslant&\frac{c_l(t,f)}{2}\int_Mf^l|\nabla^lu|^2\,dM+c_l(t,f)\sqrt{\int_Mf^l|\nabla^lu|^2\,dM}(||Q_1(u)||_2+||Q_2(u)||_2)\nonumber\\
 &&+c_l(t,f)||\nabla u||_{H^{l-1,2}}(||\nabla Q_3(u)||_2+||Q_3(u)||_2)\\
 &&+c_l(t,f)||\nabla^l u||_2(||\nabla Q_4(u)||_2+||Q_4(u)||_2)\nonumber\\
 &&+c_l(t,f)\sum\limits_j(||Q_5^j(u)||_2+||\nabla Q_5^j(u)||_2)||\nabla^lu||_2+c_l(t,f)\sum\limits_j||Q_5^j(u)||_2||Q_6^j(u)||_2\nonumber\\
 &&+c_l(t,f)\sum\limits_j||Q_6^j(u)||_2||\nabla^lu||_2,\nonumber
\end{eqnarray}
where
\[\aligned
&c_l(t,f):=&(c(l,\Omega)+1)\cdot\max\Big\{\kappa_1,\,\kappa_2,\, \kappa_3,\, \kappa_4, \, \kappa_5\Big\},
\endaligned\]

$$\kappa_1:=\frac{l}{\delta}\cdot\max\limits_{x\in M}(|f_t(x,t)|)+\frac{2l^2+2}{\delta}\cdot\max\limits_{x\in M}(|\nabla f(x,t)|^2),$$

$$\kappa_2:=\max\{\max_{x\in M}(f^{\frac{l}{2}}(x,t)),\,\max\limits_{x\in M}(f^l(x,t)),\, \max\limits_{x\in M}(f^l(x,t)|\nabla f(x,t)|)\},$$

$$\kappa_3:=\max\limits_{0\leqslant i\leqslant l-1}[\max\limits_{x\in M}(|\nabla^if(x,t)|)]\cdot\max\{\max\limits_{x\in M}(f^l(x,t)),\,\max\limits_{x\in M}(f^{l-1}(x,t)|\nabla f(x,t)|)\},$$

$$\kappa_4:=\max\{\max\limits_{x\in M}(f^l(x,t)), \s 2\max\limits_{x\in M}(|\nabla f(x,t)|\cdot f^{l-1}(x,t))\},$$
and

$$\kappa_5:=\max\limits_{x\in M}(|\nabla f(x,t)|^2f^{l-1}(x,t)+f^l(x,t)|\nabla^2f(x,t)|).$$
\medskip

In order to estimate $$\frac{1}{2}\frac{d}{dt}\int_Mf^l|\nabla^lu|^2\,dM$$ by the Sobolev norms of $u$, we need to estimate the quantities $||Q_1(u)||_2$, $||Q_2(u)||_2$, $||Q_3(u)||_2$, $||\nabla Q_3(u)||_2$, $||Q_4(u)||_2$, $||\nabla Q_4(u)||_2$, $||Q_5^j(u)||_2$, $||\nabla Q_5^j(u)||_2$, $||Q^j_6(u)||_2$ in (\ref{ZL}). Therefore, we need to
establish the following lemma.

\begin{lem}\label{NSFlemma3.2}
If $m\geqslant2$ and $2\leqslant l\leqslant m_0$, there exist constants $c(l,\Omega)$, $\tilde{c}_l(t,f)$, $C(l,\Omega)$, $c(M,l,\Omega)$ which are independent of $\varepsilon$ such that
\medskip

(1). for $Q_1(u)$ which satisfies (\ref{T:4}) with $(j_1,\cdots,j_s,p,q)$ as stated there, there holds true
\[||Q_1(u)||_2\leqslant\tilde{c}_l(t,f)C(l,\Omega)\sum\limits_{p=0}^{l-1}\sum\limits_{s=2}^{l+1-p}||\nabla u||_{H^{m_0,2}}^{\alpha(l,p,s)}||\nabla u||_2^{\beta(l,p,s)},
\]
where
\[\alpha(l,p,s):=\frac{l+1-p}{m_0}+\frac{m-2}{2m_0}s\]
and
\[\beta(l,p,s):=s+1-\alpha(l,p,s).\]
\medskip

(2). for $Q_2(u)$ which satisfies (\ref{T:6}) with $(j_1,\cdots,j_s, p, q)$ satisfy the corresponding conditions stated there, there holds true
\[
||Q_2(u)||_2\leqslant\tilde{c}_l(t,f)C(l,\Omega)\sum\limits_{p=0}^{l-1}\sum\limits_{s=2}^{l+1-p}||\nabla u||_{H^{m_0,2}}^{\alpha(l,p,s)}||\nabla u||_2^{\beta(l,p,s)}.
\]
\medskip

(3). for $Q_3(u)$ which satisfies (\ref{Q:3}) with $(j_1,\cdots,j_s)$ satisfy the corresponding conditions stated there, there hold true
\[||Q_3(u)||_2\leqslant c(M,l,\Omega)\sum\limits_{s=3}^{l+1}||\nabla u||_{H^{m_0,2}}^{\gamma(l,s)}||\nabla u||_2^{s-\gamma(l,s)},\]
where
\[\gamma(l,s):=[l+1+(\frac{m}{2}-1)s-\frac{m}{2}]/m_0;\]
and
\[||\nabla Q_3(u)||_2\leqslant c(M,l,\Omega)\sum\limits_{s=3}^{l+2}||\nabla u||_{H^{m_0,2}}^{\theta(l,s)}||\nabla u||_2^{s-\theta(l,s)}, \]
where
\[\theta(l,s):=[l+2+(\frac{m}{2}-1)s-\frac{m}{2}]/m_0.\]
\medskip

(4). for $Q_4(u)$ which satisfies (\ref{Q:4}) there hold true
\[||Q_4(u)||_2\leqslant c(l,\Omega)\tilde{c}_l(t,f)||\nabla u||_{H^{l-2,2}},\]
and
\[||\nabla Q_4(u)||_2\leqslant c(l,\Omega)\tilde{c}_l(t,f)||\nabla u||_{H^{l-1,2}}.\]
\medskip

(5). for $Q_5^j(u)$ which satisfies (\ref{Q:5}) with $(j_1,\cdots,j_s)$ satisfy the corresponding conditions stated there, there hold true
\[||Q_5^j(u)||_2\leqslant c(M,l,\Omega)\sum\limits_{s=3}^l||\nabla u||_{H^{m_0,2}}^{\psi(l,s)}||\nabla u||_2^{s-\psi(l,s)},\]
where
\[\psi(l,s):=\left[l+(\frac{m}{2}-1)s-\frac{m}{2}\right]/m_0;\]
and
\[||\nabla Q_5^j(u)||_2\leqslant c(M,l,\Omega)\sum\limits_{s=3}^{l+1}||\nabla u||_{H^{m_0,2}}^{\gamma(l,s)}||\nabla u||_2^{s-\gamma(l,s)}.\]
\medskip

(6). for $Q^j_6(u)$ which satisfies (\ref{Q:6}) with $(j_1,\cdots,j_s)$, there holds true
\[||Q_6^j(u)||_2\leqslant c(M,l,\Omega)\sum\limits_{s=3}^{l+1}||\nabla u||_{H^{m_0,2}}^{\gamma(l,s)}||\nabla u||_2^{s-\gamma(l,s)}.\]
\end{lem}

\medskip
\noindent\textbf{Proof:}
We set
\[\tilde{c}_l(t,f):=\max\limits_{x\in M}(\max\{1,f(x,t),|\nabla f(x,t)|,\cdots,|\nabla^{l+1}f(x,t)|\}).\]
From (\ref{T:4}) we get
\[|Q_1(u)|\leqslant c(l,\Omega)\cdot\tilde{c}_l(t,f)\sum\limits_{p=0}^{l-1}I_p,\]
where
\[I_p:=\sum|\nabla^{j_1}u|\cdots|\nabla^{j_s}u|\cdot|\nabla^{q+1}u|\]
with $s\geqslant 2$ and $$j_1+\cdots +j_s+q=l+1-p,$$ where
 $$1\leqslant j_i\leqslant\min\{l-1,l-p\}\s\s\mbox{and}\s\s \max\{1-p,0\}\leqslant q\leqslant l-1-p.$$
Obviously, we have
\[||Q_1(u)||_2\leqslant c(l,\Omega)\cdot\tilde{c}_l(t,f)\sum\limits_{p=0}^{l-1}||I_p||_2.\]
We want to derive the following
\begin{equation}\label{C:0}
||I_p||_2\leqslant\sum||\nabla^{j_1}u||_{p_1}\cdots||\nabla^{j_s}u||_{p_s}\cdot||\nabla^{q+1}u||_{t^*},
\end{equation}
where $p_i$ ($i=1,\cdots, s$) and $t^*$, belonging to $[1,\infty]$, will be determined later and satisfy
\begin{equation}\label{C:1}
\frac{1}{p_1}+\cdots+\frac{1}{p_s}+\frac{1}{t^*}=\frac{1}{2}.
\end{equation}
And then we employ Theorem 2.1 due to \cite{DW1} to obtain
\[
||\nabla^{j_i}u||_{p_i}\leqslant C(j_i,a_i)||\nabla u||^{a_i}_{H^{m_0,2}}||\nabla u||_2^{1-a_i},
\]
and
\[||\nabla^{q+1}u||_{t^*}\leqslant C(q,b)||\nabla u||^b_{H^{m_0,2}}||\nabla u||^{1-b}_2.\]
We hope that $p_i(i=1,\cdots,s)$ and $t^*$ satisfy the following conditions:
\[
\frac{1}{p_i}=\frac{j_i-1}{m}+\frac{1}{2}-a_i\frac{m_0}{m}\s \s\mbox{with}\s\s a_i\in\Big[\frac{j_i-1}{m_0},\, 1\Big)\]
which is equivalent to
\begin{equation}\label{C:2}
\frac{1}{p_i}\in\Big(\frac{1}{2}+\frac{j_i-1-m_0}{m},\, \frac{1}{2}\Big];
\end{equation}
and
\[
\frac{1}{t^*}=\frac{q}{m}+\frac{1}{2}-b\frac{m_0}{m}\s\s\mbox{with}\s\s b\in\Big[\frac{q}{m_0},\, 1\Big)\]
which is equivalent to
\begin{equation}\label{C:3}
\frac{1}{t^*}\in\Big(\frac{1}{2}+\frac{q-m_0}{m},\,\frac{1}{2}\Big].
\end{equation}

\medskip
We claim that there always exist $p_1, \cdots, p_s, t^*$ which are in $[1,\infty]$ and satisfy (\ref{C:1}), (\ref{C:2}) and (\ref{C:3}).

Next we prove the claim. Let $A:=\{i\,|\,j_i\geqslant m_0+1-\frac{m}{2}\}$. We need to consider the following two cases:
\medskip

\noindent (a). The case  $q=0$.  It is easy to see that $q=0$ implies that $$\frac{1}{2}+\frac{q-m_0}{m}<0.$$
From (\ref{C:2})and (\ref{C:3}) we see that there holds in this case
\[\frac{1}{p_1}+\cdots+\frac{1}{p_s}+\frac{1}{t^*}\in\Big(\sum\limits_{i\in A}\left(\frac{1}{2}+\frac{j_i-1-m_0}{m}\right),\, \frac{s+1}{2}\Big].\]
That is to say that, when $q=0$, there exist $p_1, \cdots, p_s, t^*\in [1,\infty]$ satisfying
$$\frac{1}{p_i}\in[0,1]\cap\Big(\frac{1}{2}+\frac{j_i-1-m_0}{m},\, \frac{1}{2}\Big]\s\s
\mbox{and}\s\s \frac{1}{t^*}\in[0,1]\cap\Big(\frac{1}{2}+\frac{q-m_0}{m},\,\frac{1}{2}\Big]$$ such that $$\frac{1}{p_1}+\cdots+\frac{1}{p_s}+\frac{1}{t^*}=\frac{1}{2}$$
if and only if
\[
\sum\limits_{i\in A}\left(\frac{1}{2}+\frac{j_i-1-m_0}{m}\right)<\frac{1}{2}.
\]
The above inequality can also be written as
\[\left(\frac{1}{2}-\frac{1+m_0}{m}\right)|A|+\frac{\sum\limits_{i\in A}j_i}{m}<\frac{1}{2}.\]
\medskip

\noindent (b). The case $q\geqslant1$. Obviously, $q\geqslant1$ implies that $$\frac{1}{2}+\frac{q-m_0}{m}\geqslant0.$$
Therefore, we can also see easily from (\ref{C:2}) and (\ref{C:3})
\[\frac{1}{p_1}+\cdots+\frac{1}{p_s}+\frac{1}{t^*}\in\Big(\sum\limits_{i\in A}\left(\frac{1}{2}+\frac{j_i-1-m_0}{m}\right)+\frac{1}{2}+\frac{q-m_0}{m},\,\frac{s+1}{2}\Big].\]
This means that there exist $p_1, \cdots, p_s, t^*\in [1,\infty]$ satisfying
$$\frac{1}{p_i}\in[0,1]\cap\Big(\frac{1}{2}+\frac{j_i-1-m_0}{m},\,\frac{1}{2}\Big]\s\s \mbox{and}\s\s \frac{1}{t^*}\in[0,1]\cap\Big(\frac{1}{2}+\frac{q-m_0}{m},\,\frac{1}{2}\Big]$$
such that $$\frac{1}{p_1}+\cdots+\frac{1}{p_s}+\frac{1}{t^*}=\frac{1}{2}$$ if and only if
\[\sum\limits_{i\in A}\Big(\frac{1}{2}+\frac{j_i-1-m_0}{m}\Big)+\frac{1}{2}+\frac{q-m_0}{m}<\frac{1}{2}.\]
\medskip

In the case $q=0$, we have
\[
\aligned
(\frac{1}{2}-\frac{1+m_0}{m})|A|+\frac{\sum\limits_{i\in A}j_i}{m}
=&(\frac{1}{2}-\frac{1+m_0}{m})|A|+\frac{l+1-p-\sum\limits_{i\in A^c}j_i}{m}\\
\leqslant&(\frac{1}{2}-\frac{1+m_0}{m})|A|+\frac{l+1-p-|A^c|}{m}\\
=&\frac{1}{m}[(\left\{\frac{m}{2}\right\}-1)|A|+l+1-p-s],
\endaligned
\]
where $A^c:=\{1,2,\cdots,s\}\setminus A$ and $$\left\{\frac{m}{2}\right\}:=\frac{m}{2}-\left[\frac{m}{2}\right].$$

Noting that $\{\frac{m}{2}\}-1=-\frac{1}{2}$ or $-1$ as $m$ is integer, we obtain
\[
\frac{1}{m}\left[(\left\{\frac{m}{2}\right\}-1)|A|+l+1-p-s\right]\leqslant\frac{1}{m}\left[-\frac{|A|}{2}+l+1-p-s\right].
\]

If $|A|\geqslant1$, then we have
\[|A|>2\Big[\frac{m}{2}\Big]-m,\]
which is equivalent to that
\[2(m_0-1)-m<|A|.\]
Since $l\leqslant m_0$, the following holds
\[2(l+1-2)-m<|A|,\]
so, we have
\[2(l+1-s)-m<|A|.\]
This implies that
\[l+1-p-s-\frac{m}{2}<\frac{|A|}{2},\]
which is equivalent to that
\[\frac{1}{m}[-\frac{|A|}{2}+l+1-p-s]<\frac{1}{2}.\]
So, we get
\[\frac{1}{m}\left[\left(\left\{\frac{m}{2}\right\}-1\right)|A|+l+1-p-s\right]<\frac{1}{2}.\]
This means that we can find the required $p_i$ and $t^*$.
\medskip

If $|A|=0$, the conclusion is  obviously true.
\medskip

In the case $q\geqslant1$, we have
\[
\aligned
&\sum\limits_{i\in A}\left(\frac{1}{2}+\frac{j_i-1-m_0}{m}\right)+\frac{1}{2}+\frac{q-m_0}{m}\\
=&\left(\frac{1}{2}-\frac{m_0+1}{m}\right)|A|+\frac{1}{m}[l+1-p-m_0-\sum\limits_{i\in A^c}j_i]+\frac{1}{2}\\
\leqslant&\left(\frac{1}{2}-\frac{m_0+1}{m}\right)|A|+\frac{1}{m}[l+1-p-m_0-|A^c|]+\frac{1}{2}\\
=&\left(\frac{1}{2}-\frac{m_0+1}{m}\right)|A|+\frac{1}{m}[l+1-p-m_0-s+|A|]+\frac{1}{2}.
\endaligned
\]
Since $$\frac{|A|}{2}+1>0,$$
we have
\[-\frac{|A|}{2}+l-m_0-1<0.\]
So, in view of the fact $p\geqslant0$ and $s\geqslant2$ we obtain that
\[-\frac{|A|}{2}+l+1-p-m_0-s<0.\]
It follows that
\[\left(\left\{\frac{m}{2}\right\}-1\right)|A|+l+1-p-m_0-s<0,\]
which is equivalent to that
\[\left(\frac{m}{2}-m_0\right)|A|+l+1-p-m_0-s<0.\]
Moreover, the above inequality is also equivalent to that
\[\left(\frac{1}{2}-\frac{m_0+1}{m}\right)|A|+\frac{1}{m}[l+1-p-m_0-s+|A|]<0.\]
So we can pick out the required $p_i$ and $t^*$. This completes the proof of the claim.
\medskip

Now let us return to the proof of Lemma \ref{NSFlemma3.2}. Note that $a_i$ is a function of $p_i$ and $b$ is a function of $t^*$. Since we have shown that, for given $(j_1,\cdots,j_s,q)$ which satisfies the conditions adhering to (\ref{T:4}), there always exist $p_i=p_i(j_1,\cdots,j_s,q)$ ($i=1, \cdots, s$) and $t^*=t^*(j_1,\cdots,j_s,q)$ such that (\ref{C:0}) holds true . So $a_i$ and $b$ can be written as
\[a_i=a_i(j_1,\cdots,j_s,q)\s\s \mbox{and}\s\s b=b(j_1,\cdots,j_s,q).\]
Setting
\[
\aligned
m^*(l,p):=\max\{&C(j_1,a_1),\cdots,C(j_s,a_s),C(q,b)|s\geqslant 2,\, j_1+\cdots +j_s+q=l+1-p,\\
&1\leqslant j_i\leqslant\min\{l-1,l-p\},\, \max\{1-p,0\}\leqslant q\leqslant l-1-p\}.
\endaligned
\]
Then we get
\[
||I_p||_2\leqslant \sum\limits_{s=2}^{l+1-p}(m^*(l,p))^{s+1}||\nabla u||_{H^{m_0,2}}^{a_1+\cdots+a_s+b}||\nabla u||_2^{s+1-a_1-\cdots-a_s-b}.
\]
Since
\[a_1+\cdots+a_s+b=\frac{l+1-p}{m_0}+\frac{m-2}{2m_0}s:=\alpha(l,p,s),\]
we get
\[
||Q_1||_2\leqslant\tilde{c}_l(t,f)C(l,\Omega)\sum\limits_{p=0}^{l-1}\sum\limits_{s=2}^{l+1-p}||\nabla u||_{H^{m_0,2}}^{\alpha(l,p,s)}||\nabla u||_2^{\beta(l,p,s)},
\]
where
\[\beta(l,p,s):=s+1-\alpha(l,p,s)\]
and
\[C(l,\Omega):=c(l,\Omega)\cdot\max\limits_{0\leqslant p\leqslant l-1}\Big\{\max\limits_{2\leqslant s\leqslant l+1-p}\{(m^*(l,p))^{s+1}\}\Big\}.\]
Thus, we finish the proof of (1) of Lemma \ref{NSFlemma3.2}.
\medskip

By the same argument as above, similarly we may also prove
\[
||Q_2||_2\leqslant\tilde{c}_l(t,f)C(l,\Omega)\sum\limits_{p=0}^{l-1}\sum\limits_{s=2}^{l+1-p}||\nabla u||_{H^{m_0,2}}^{\alpha(l,p,s)}||\nabla u||_2^{\beta(l,p,s)}.
\]

Using the same method as in Lemma 3.2 in \cite{DW1}, we get the remaining estimates in this lemma. Here, we omit the details.
\begin{rem}
If $m=1$, then $m_0=1$. At this time, Lemma \ref{NSFlemma3.2} does not work and readers can refer directly to Lemma \ref{NSFlemma3.3}.
\end{rem}

\begin{lem}\label{NSFlemma3.3}
Let $||Q_i||_2\,(i=1, 2, \cdots, 6)$, $||\nabla Q_3||_2$, $||\nabla Q_4||_2$ and $||\nabla Q_5^j||$ be the quantities in (\ref{ZL}). If $m\geqslant1$ and $m_0+1\leqslant l\leqslant k$, there exist constants $\tilde{c}_1(l,\Omega,t,f)$, $\cdots$, $\tilde{c}_7(l,\Omega,t,f)$, $c(l,\Omega)$, $\tilde{c}_l(t,f),\, \omega_1(l),\cdots,\omega_6(l)$ which are independent of $\varepsilon$ such that
\[||Q_1||_2\leqslant\tilde{c}_1(l,\Omega,t,f)(1+||\nabla u||_{H^{l-1,2}})(1+||\nabla u||_{H^{\max\{l-2,m_0\},2}})^{\omega_1(l)},\]

\[||Q_2||_2\leqslant\tilde{c}_2(l,\Omega,t,f)(1+||\nabla u||_{H^{l-1,2}})(1+||\nabla u||_{H^{\max\{l-2,m_0\},2}})^{\omega_2(l)},\]

\[||Q_3||_2\leqslant\tilde{c}_3(l,\Omega,t,f)(1+||\nabla u||_{H^{\max\{l-2,m_0\},2}})(1+||\nabla u||_{H^{\max\{l-3,m_0\},2}})^{\omega_3(l)},\]

\[||\nabla Q_3||_2\leqslant\tilde{c}_4(l,\Omega,t,f)(1+||\nabla u||_{H^{l-1,2}})(1+||\nabla u||_{H^{\max\{l-2,m_0\},2}})^{\omega_4(l)},\]

\[||Q_4||_2\leqslant c(l,\Omega)\tilde{c}_l(t,f)||\nabla u||_{H^{l-2,2}},\]

\[||\nabla Q_4||_2\leqslant c(l,\Omega)\tilde{c}_l(t,f)||\nabla u||_{H^{l-1,2}},\]

\[||Q_5^j||_2\leqslant\tilde{c}_5(l,\Omega,t,f)(1+||\nabla u||_{H^{\max\{l-3,m_0\},2}})(1+||\nabla u||_{H^{\max\{l-4,m_0\},2}})^{\omega_5(l)},\]

\[||\nabla Q_5^j||_2\leqslant\tilde{c}_6(l,\Omega,t,f)(1+||\nabla u||_{H^{\max\{l-2,m_0\},2}})(1+||\nabla u||_{H^{\max\{l-3,m_0\},2}})^{\omega_6(l)},\]
and
\[||Q_6^j||_2\leqslant\tilde{c}_7(l,\Omega,t,f)(1+||\nabla u||_{H^{\max\{l-2,m_0\},2}})(1+||\nabla u||_{H^{\max\{l-3,m_0\},2}})^{\omega_6(l)}.\]
\end{lem}

\medskip

\noindent\textbf{Proof:} First of all, let us consider $Q_1$.

When $l=m_0+1$, from $(\ref{T:4})$ one can obtain
\[|Q_1|\leqslant c(m_0+1,\Omega)\tilde{c}_{m_0+1}(t,f)\sum\limits_{p=0}^{m_0}I_p,\]
where
\[I_p=\sum|\nabla^{j_1}u|\cdots|\nabla^{j_s}u|\cdot|\nabla^{q+1}u|\]
with
\[s\geqslant2, \s\s j_1+\cdots+j_s+q=m_0+2-p,\s\s 1\leqslant j_i\leqslant\min\{m_0,m_0+1-p\},\]
and
\[\max\{1-p,0\}\leqslant q\leqslant m_0-p.\]
It is convenient to write $I_p$ as the following expression if one denotes $\nabla^{q+1}$ by $\nabla^{j_i}$
\[I_p=\sum|\nabla^{j_1}u|\cdots|\nabla^{j_s}u|\]
with
\[s\geqslant3,\s\s j_1+\cdots+j_s=m_0-p+3,\s\s 1\leqslant j_i\leqslant m_0-p+1.\]
\medskip

We need to estimate $I_p$ ($p=0, \cdots, m_0$). By Lemma 3.2 of \cite{DW1}, we have
\[||I_p||_2\leqslant c(M,m_0-p)||\nabla u||_{H^{m_0,2}}^{\varphi(s)}||\nabla u||_2^{s-\varphi(s)},\]
where
\[\varphi(s)=\left[m_0-p+3+\left(\frac{m}{2}-1\right)s-\frac{m}{2}\right]/m_0.\]
So,
\[||Q_1||_2\leqslant c(m_0+1,\Omega)\tilde{c}_{m_0+1}(t,f)\sum\limits_{p=0}^{m_0}\sum\limits_{s=3}^{m_0-p+3}c(M,m_0-p)||\nabla u||_{H^{m_0,2}}^{\varphi(s)}||\nabla u||_2^{s-\varphi(s)}.\]

When $l\geqslant m_0+2$, we have
\[|Q_1|\leqslant c(l,\Omega)\tilde{c}_l(t,f)\left(\sum\limits_{p=0}^{l-2-m_0}I_p+\sum\limits_{p=l-1-m_0}^{l-1}I_p\right).\]
We need to consider following two cases.
\medskip

\noindent\textsf{Case 1: $l-1-m_0\leqslant p\leqslant l-1$.}
\medskip

In this case, we have
\[I_p=\sum|\nabla^{j_1}u|\cdots|\nabla^{j_s}u|,\]
where
\[s\geqslant3,\,\,\,\,\,\,j_1+\cdots+j_s=l+2-p,\,\,\,\,\,\,1\leqslant j_i\leqslant l-p.\]
Since $l-1-p\leqslant m_0$, by Lemma 3.2 in \cite{DW1}, there holds true
\[||I_p||_2\leqslant c(M,l-1-p)\sum\limits_{s=3}^{l+2-p}||\nabla u||_{H^{m_0,2}}^{\chi(l,p,s)}||\nabla u||_2^{s-\chi(l,p,s)},\]
where
\[\chi(l,p,s):=\left[l+2-p+\left(\frac{m}{2}-1\right)s-\frac{m}{2}\right]/m_0.\]
\medskip

\noindent\textsf{Case 2: $0\leqslant p\leqslant l-2-m_0$.}
\medskip

For this case, we have
\[I_p=\sum|\nabla^{j_1}u|\cdots|\nabla^{j_s}u|,\]
where
\begin{equation}\label{H:0}
s\geqslant3,\s j_1+\cdots+j_s=l-p+2,\s 1\leqslant j_i\leqslant l-p,\s l-p\geqslant j_1\geqslant\cdots\geqslant j_s\geqslant1.
\end{equation}
In the present situation, we need to discuss the following two subcases.
\medskip

\noindent\textsf{Subcase 1. $1\leqslant p\leqslant l-2-m_0$.}

In order to obtain the required estimate, we need to apply Theorem 2.1 of \cite{DW1} and use H\"older inequality delicately. This ask us to pick out $p_1,\cdots,p_s\in[1,\infty]$ such that
\begin{equation}\label{H:1}
\frac{1}{p_1}+\cdots+\frac{1}{p_s}=\frac{1}{2},
\end{equation}

\begin{equation}\label{H:2}
\frac{1}{p_i}=\frac{j_i-1}{m}+\frac{1}{2}-\frac{l-2}{m}a_i,
\end{equation}
and
\begin{equation}\label{H:3}
a_i\in\Big[\frac{j_i-1}{l-2},\,\,1\Big).
\end{equation}
Once that $(\ref{H:1}),(\ref{H:2})$ and $(\ref{H:3})$ are met, we have immediately
\[||I_p||_2\leqslant\sum\left\||\nabla^{j_1}u|\cdots|\nabla^{j_s}u|\right\|_2\leqslant\sum||\nabla^{j_1}u||_{p_1}\cdots||\nabla^{j_s}u||_{p_s}\]
and
\[||\nabla^{j_i}u||_{p_i}\leqslant C(j_i,l,a_i)||\nabla u||^{a_i}_{H^{l-2,2}}||\nabla u||_2^{1-a_i}\leqslant C(j_i,l,a_i)||\nabla u||_{H^{l-2,2}}.\]
In other words, we need to pick $p_1, \cdots, p_s\in[1,\,\infty]$ such that
\begin{equation}\label{H:4}
\frac{1}{p_i}\in\left(\frac{1}{2}+\frac{j_i-l+1}{m},\,\frac{1}{2}\right]\cap[0,1]\s\s
\mbox{and}\s\s \frac{1}{p_1}+\cdots+\frac{1}{p_s}=\frac{1}{2}.\end{equation}
\medskip

For this purpose we take the same argument as in the previous. Let $$A:=\left\{i|\,j_i\geqslant l-1-\frac{m}{2}\right\}.$$
The fact there exists $(p_1,\cdots, p_s)$ satisfying (\ref{H:4}) is equivalent to that
\[\sum\limits_{i\in A}\left(\frac{1}{2}+\frac{j_i-l+1}{m}\right)<\frac{1}{2}.\]
Noting (\ref{H:0}), we can write the above inequality as
\[\left(\frac{1}{2}+\frac{1-l}{m}\right)|A|+\frac{\sum\limits_{i\in A}j_i}{m}=\left(\frac{1}{2}+\frac{1-l}{m}\right)|A|+\frac{l-p+2-\sum\limits_{i\in A^c}j_i}{m}<\frac{1}{2}.\]

Indeed, we have
\begin{equation}\label{H:5}
\aligned
&\left(\frac{1}{2}+\frac{1-l}{m}\right)|A|+\frac{l-p+2-\sum\limits_{i\in A^c}j_i}{m}\\
\leqslant&\left(\frac{1}{2}+\frac{1-l}{m}\right)|A|+\frac{l-p+2-|A^c|}{m}\\
=&\left(\frac{1}{2}+\frac{1-l}{m}\right)|A|+\frac{l-p+2-(s-|A|)}{m}\\
=&\left(\frac{1}{2}+\frac{2-l}{m}\right)|A|+\frac{l-p+2-s}{m}.
\endaligned
\end{equation}
Noting
\[\frac{1}{2}+\frac{2-l}{m}<0\]
and keeping that $p\geqslant1$ and $s\geqslant3$ in mind we have that, when $|A|\geqslant2$,
\begin{equation}\label{H:6}
\aligned
&\left(\frac{1}{2}+\frac{2-l}{m}\right)|A|+\frac{l-p+2-s}{m}\\
<&\frac{1}{2}+\frac{2-l}{m}+\frac{l-p+2-s}{m}\\
=&\frac{1}{2}+\frac{4-p-s}{m}\leqslant\frac{1}{2}.
\endaligned
\end{equation}
Hence, from the above (\ref{H:5}) and (\ref{H:6}) we know that the following holds true
\[\left(\frac{1}{2}+\frac{1-l}{m}\right)|A|+\frac{\sum\limits_{i\in A}j_i}{m}<\frac{1}{2}.\]

We need to discuss the remaining two cases $|A|=0$ and $|A|=1$.

For $|A|=0$, obviously, we have
\[\left(\frac{1}{2}+\frac{1-l}{m}\right)|A|+\frac{\sum\limits_{i\in A}j_i}{m}=0<\frac{1}{2}.\]

For $|A|=1$, we have
\[l-p\geqslant j_1\geqslant l-1-\frac{m}{2}>j_2\geqslant\cdots\geqslant j_s\geqslant1.\]
In the present situation, we have
\[\sum\limits_{i\in A}\left(\frac{1}{2}+\frac{j_i-l+1}{m}\right)=\frac{1}{2}+\frac{j_1-l+1}{m}\leqslant\frac{1}{2}+\frac{1-p}{m}.\]
So, in the case $|A|=1$, obviously there holds true for $p\geqslant2$
\[\sum\limits_{i\in A}\left(\frac{1}{2}+\frac{j_i-l+1}{m}\right)\leqslant\frac{1}{2}+\frac{1-p}{m}<\frac{1}{2}.\]
Yet, for $|A|=1$ and $p=1$ we have
\[\sum\limits_{i\in A}\left(\frac{1}{2}+\frac{j_i-l+1}{m}\right)\leqslant\frac{1}{2}.\]
However, for this special case we have
\[
\aligned
\left\||\nabla^{j_1}u|\cdots|\nabla^{j_s}u|\right\|_2
\leqslant&||\nabla^{j_1}u||_2||\nabla^{j_2}u||_{\infty}\cdots||\nabla^{j_s}u||_{\infty}\\
\leqslant&||\nabla u||_{H^{l-2,2}}||\nabla^{j_2}u||_{\infty}\cdots||\nabla^{j_s}u||_{\infty}.
\endaligned
\]
On the other side, since we have for $2\leqslant i\leqslant s$
\[\frac{1}{2}+\frac{j_i-l+1}{m}<0,\]
the Sobolev imbedding theorem implies that
\[||\nabla^{j_i}u||_{\infty}\leqslant C(j_i,l,a_i)||\nabla u||^{a_i}_{H^{l-2,2}}||\nabla u||_2^{1-a_i}\leqslant C(j_i,l,a_i)||\nabla u||_{H^{l-2,2}}.\]

From the above arguments we see that, for the case $1\leqslant p\leqslant l-2-m_0$, there always holds true
\[\|I_p\|_2\leqslant\sum\limits_{s=3}^{l-p+2}(M^*(l,p))^s||\nabla u||^s_{H^{l-2,2}},\]
where
\[M^*(l,p)=\max\left\{C(j_i,l,a_i)|\, i=1,\cdots,s; \s \sum_{i=1}^sj_i=l-p+2,1\leqslant j_i\leqslant l-p;\s s\geqslant3\right\}.\]

\noindent\textsf{Subcase 2: $p=0$}.
\medskip

Now, we have
\[I_0=\sum|\nabla^{j_1}u|\cdots|\nabla^{j_s}u|,\]
where
\[s\geqslant3,\,\,\,\,\,\,j_1+\cdots+j_s=l+2,\,\,\,\,\,\,1\leqslant j_i\leqslant l.\]
By the proof of Lemma 3.3 in \cite{DW1}, we infer
\[||I_0||_2\leqslant c(M,l)(1+||\nabla u||_{H^{l-1,2}})(1+||\nabla u||^{\sigma(l)}_{H^{l-2,2}}).\]
Therefore, when $l\geqslant m_0+2$, we also prove that there exists a $N^*(l,\Omega)$ such that
\[||Q_1||_2\leqslant N^*(l,\Omega)\cdot\tilde{c}_l(t,f)(1+||\nabla u||_{H^{l-1,2}})(1+||\nabla u||_{H^{l-2,2}})^{\omega(l)}.\]

By summarizing the above arguments we conclude that, for $k\geqslant l\geqslant m_0+1$, there always exists $\omega_1(l)$ such that
\[||Q_1||_2\leqslant\tilde{c}_1(l,\Omega,t,f)(1+||\nabla u||_{H^{l-1,2}})(1+||\nabla u||_{H^{\max\{l-2,m_0\},2}})^{\omega_1(l)}.\]
Here, \[\tilde{c}_1(l,\Omega,t,f)=N^*(l,\Omega)\cdot\tilde{c}_l(t,f).\]

By the same method, we can also get similar estimate for $Q_2$. Using the same trick of lemma 3.3 in \cite{DW1}, one may easily have those estimates for $Q_3$, $\nabla Q_3$, $Q_5^j$, $\nabla Q_5^j$, $Q_6^j$. As for $Q_4$ and $\nabla Q_4$, since of (\ref{Q:4}) and (\ref{T:15}), it is obvious. Thus, we complete the proof of this lemma.\\

Next we will use Lemma \ref{NSFlemma3.2} and Lemma \ref{NSFlemma3.3} to derive the positive lower bound of $T_{\varepsilon}$ and the upper bounds of the Sobolev norms of $u$.
\begin{lem}\label{NSFlemma3.4}
Let $m_0=[\frac{m}{2}]+1$ where $[q]$ is the integral part of a real number $q$. $u$ is a solution of $(\ref{T:1})$ and $[0,T_{\varepsilon})$ is its existing interval. Then
\[T_{\varepsilon}\geqslant T(f,m_0,||\nabla u_0||^2_{H^{m_0,2}(f_0)},T_*,M,N):=T\]
and
\[||u(t)||_{W^{k+1,2}(M,\mathbb{R}^L)}\leqslant C(M,N,f,k,T_*,||\nabla u_0||_{H^{k,2}(f_0)}):=F(k)\]
for all $k\geqslant m_0$ and $t\in[0,T]$.
\end{lem}

\textbf{Proof:} For $k=m_0$, by substituting the upper bounds of $||Q_i||_2$ ($i=1,\cdots, 6$) and $||\nabla Q_j||_2$($j=3, 4, 5$) in Lemma \ref{NSFlemma3.2} and Lemma \ref{NSFlemma3.3} into
$$\sum\limits_{l=1}^{m_0+1}\frac{d}{dt}\int_Mf^l|\nabla^lu|^2\,dM,$$
from (\ref{ZL}) we can see that there exist $Q(m_0)\geqslant1$, which depends only on $m_0$, and a function $D_{m_0}(t,f)$ such that
there holds
\begin{equation}\label{B:1}
\left\{
\begin{array}{llll}
\displaystyle\frac{d}{dt}||\nabla u||^2_{H^{m_0,2}(f)}\leqslant D_{m_0}(t,f)(1+||\nabla u||^2_{H^{m_0,2}(f)})^{Q(m_0)},\\
||\nabla u||^2_{H^{m_0,2}(f)}(0)=||\nabla u_0||^2_{H^{m_0,2}(f_0)}.
\end{array}
\right.
\end{equation}
Here,
\[
\aligned
D_{m_0}(t,f):=&a(m_0,\Omega)\cdot\max\{c_l(t,f)|1\leqslant l\leqslant m_0+1\}\\
&\cdot\max\{\tilde{c}_i(l,\Omega,t,f)|1\leqslant i\leqslant7,1\leqslant l\leqslant m_0+1\}
\endaligned
\]
and $a(m_0,\Omega)$ is a positive number depends only on $m_0$ and $\Omega$.

Noting that $D_{m_0}(\cdot,f)\in C^0([0,T_*])$, we set
\[D_{m_0}(f):=\max\limits_{s\in[0,T_*]}\{D_{m_0}(s,f)\}\]
and consider the following
\begin{equation}\label{B:2}
\left\{
\begin{array}{llll}
\displaystyle\frac{dU}{dt}=D_{m_0}(f)(1+U)^{Q(m_0)},\\
U(0)=||\nabla u_0||^2_{H^{m_0,2}(f_0)}.
\end{array}
\right.
\end{equation}
Let
\[T_{m_0}:=\min\Big\{\frac{1}{2(Q(m_0)-1)\cdot D_{m_0}(f)\cdot(1+||\nabla u_0||^2_{H^{m_0,2}(f_0)})^{Q(m_0)-1}},T_*\Big\}.\]
It is easy to see that, for $t\in[0,T_{m_0}]$,
\[
U(t):=\frac{1+||\nabla u_0||^2_{H^{m_0,2}(f_0)}}{\Big(1-(Q(m_0)-1)\cdot D_{m_0}(f)\cdot(1+||\nabla u_0||^2_{H^{m_0,2}(f_0)})^{Q(m_0)-1}\cdot t\Big)^{\frac{1}{Q(m_0)-1}}}-1
\]
is the solution to $(\ref{B:2})$. Therefore, the comparison principle tells us that, for all $t\in[0,\,\min\{T'_{\varepsilon},T_{m_0}\}]$, there holds
\[||\nabla u(t)||^2_{H^{m_0,2}(f)}\leqslant U(t)\leqslant U(T_{m_0}).\]

Similarly, for $k\geqslant m_0+1$, there exists $G(k)$ depending only on $k$ and $b(k,\Omega)$ such that, if we let
\[
\aligned
\tilde{h}(k,\Omega,t,f):=&b(k,\Omega)\cdot\max\{1,c_l(t,f)|\,1\leqslant l\leqslant k+1\}\\
&\cdot\max\{1,\,\tilde{c}_i(l,\Omega,t,f)|\,1\leqslant i\leqslant7,\,1\leqslant l\leqslant k+1\}
\endaligned
\]
and
\[\tilde{h}(k,\Omega,f):=\max\limits_{s\in[0,T_*]}\{\tilde{h}(k,\Omega,s,f)\},\]
then there holds
\begin{equation}\label{B:3}
\aligned
&\frac{d}{dt}||\nabla u||^2_{H^{k,2}(f)}
=\frac{d}{dt}||\nabla u||^2_{H^{m_0,2}(f)}+\sum\limits_{l=m_0+2}^{k+1}\frac{d}{dt}\int_Mf^l|\nabla^lu|^2\,dM\\
\leqslant&D_{m_0}(f)(1+U(t))^{Q(m_0)}+\tilde{h}(k,\Omega,f)(1+||\nabla u||^2_{H^{k,2}(f)})(1+||\nabla u||^2_{H^{k-1,2}(f)})^{G(k)}.
\endaligned
\end{equation}
\medskip

As $k=m_0+1$, we have
\[
\aligned
\frac{d}{dt}||\nabla u||^2_{H^{m_0+1,2}(f)}
\leqslant&D_{m_0}(f)(1+U(T_{m_0}))^{Q(m_0)}\\
&+\tilde{h}(m_0+1,\Omega,f)(1+||\nabla u||^2_{H^{m_0+1,2}(f)})(1+U(T_{m_0}))^{G(m_0+1)}.
\endaligned
\]
Letting
\[A_{m_0}:=D_{m_0}(f)(1+U(T_{m_0}))^{Q(m_0)}+\tilde{h}(m_0+1,\Omega,f)(1+U(T_{m_0}))^{G(m_0+1)}\]
and
\[B_{m_0}:=\tilde{h}(m_0+1,\Omega,f)(1+U(T_{m_0}))^{G(m_0+1)}.\]
Now, we consider the initial value problem of the following ordinary differential equation
\[
\left\{
\begin{array}{llll}
\displaystyle\frac{d}{dt}U_{m_0+1}=A_{m_0}+B_{m_0}\cdot U_{m_0+1},\\
U_{m_0+1}(0)=||\nabla u_0||^2_{H^{m_0+1,2}(f_0)}.
\end{array}
\right.
\]
The solution of the initial value problem is
\[
U_{m_0+1}(t)=\exp\{B_{m_0}t\}\cdot||\nabla u_0||^2_{H^{m_0+1,2}(f_0)}+A_{m_0}\frac{\exp\{B_{m_0}t\}-1}{B_{m_0}}.
\]
Hence, we have that, for all $t\in[0,\,\min\{T'_{\varepsilon},T_{m_0}\}]$, there holds true
\[||\nabla u(t)||^2_{H^{m_0+1,2}(f)}\leqslant U_{m_0+1}(t).\]

Next we take an induction argument. Supposing that there exists a monotonously increasing function $U_{k-1}\in C^{\infty}(\mathbb{R}^1)$ such that, for all $t\in[0,\,\min\{T'_{\varepsilon},T_{m_0}\}]$,
\[||\nabla u(t)||^2_{H^{k-1,2}(f)}\leqslant U_{k-1}(t),\]
then, by using $(\ref{B:3})$ we obtain that
\[
\aligned
&\frac{d}{dt}||\nabla u||^2_{H^{k,2}(f)}\\
\leqslant&D_{m_0}(f)(1+U(T_{m_0}))^{Q(m_0)}+\tilde{h}(k,\Omega,f)(1+||\nabla u||^2_{H^{k,2}(f)})(1+U_{k-1}(T_*))^{G(k)}.
\endaligned
\]
By Letting
\[
A_k:=D_{m_0}(f)(1+U(T_{m_0}))^{Q(m_0)}+\tilde{h}(k,\Omega,f)(1+U_{k-1}(T_*))^{G(k)}
\]
and
\[B_k:=\tilde{h}(k,\Omega,f)(1+U_{k-1}(T_*))^{G(k)},\]
we consider
\[
\left\{
\begin{array}{llll}
\frac{dU_k}{dt}=A_k+B_k\cdot U_k,\\
U_k(0)=||\nabla u_0||^2_{H^{k,2}(f_0)}.
\end{array}
\right.
\]
Then, the solution to this initial problem can be express precisely as
\[
U_k(t)=||\nabla u_0||^2_{H^{k,2}(f_0)}\cdot\exp\{B_kt\}+\frac{A_k}{B_k}(\exp\{B_kt\}-1)
\]
which is in $C^{\infty}(\mathbb{R}^1)$ and monotonously increasing. By the comparison theorem of ODE, we know that, for all $t\in[0,\,\min\{T'_{\varepsilon},T_{m_0}\}]$,
\[||\nabla u(t)||^2_{H^{k,2}(f)}\leqslant U_k(t).\]
So, for any $n\geqslant m_0+1$, there holds true
\begin{equation}\label{B:4}
||\nabla u(t)||^2_{H^{n,2}(f)}\leqslant U_n(t), \s\s\s t\in[0,\,\min\{T'_{\varepsilon},T_{m_0}\}].
\end{equation}

\medskip

Now we turn to deriving the lower bound of $T_{\varepsilon}$. By Theorem \ref{thm:int}, there exist $a\in(0,1)$ and a constant $c(M)$ such that
\[||\nabla_tu||_{\infty}\leqslant c(M)||\nabla_tu||_{H^{m_0,2}}^a||\nabla_tu||_2^{1-a}\leqslant c(M)||\nabla_tu||_{H^{m_0,2}}.\]
Using $(\ref{T:1})$, one can easily get
\[
\aligned
|\nabla^l\nabla_tu|&=\sqrt{1+\varepsilon^2}|\nabla^l\nabla_{\beta}(f\nabla_{\beta}u)|\\
&\leqslant c(l)\sum\limits_{p+q=l+1,0\leqslant p\leqslant l+1}|\nabla^pf|\cdot|\nabla^{q+1}u|\\
&\leqslant c(l)\cdot\tilde{c}_l(t,f)\sum\limits_{q=0}^{l+1}|\nabla^{q+1}u|.
\endaligned
\]
So, we can easily see from the above inequality that, for all $t\in[0,\min\{T'_{\varepsilon},T_{m_0}\}]$, there holds
\begin{equation}\label{B:5}
\aligned
||\nabla_tu(t)||_{H^{m_0,2}}
\leqslant&\hat{c}(m_0)\cdot\tilde{c}_{m_0}(t,f)||\nabla u(t)||_{H^{m_0+1,2}}\\
\leqslant&\frac{\hat{c}(m_0)}{\delta^{(m_0+2)/2}}\cdot\tilde{c}_{m_0}(t,f)U_{m_0+1}(t)\\
\leqslant&\frac{\hat{c}(m_0)}{\delta^{(m_0+2)/2}}\cdot\max\limits_{s\in[0,T_*]}\{\tilde{c}_{m_0}(s,f)\}\cdot U_{m_0+1}(T_*):=\mu,
\endaligned
\end{equation}
where
\[\hat{c}(m_0):=(m_0+1)\cdot\max\{c(l)|0\leqslant l\leqslant m_0\}.\]
Then, we have
\[\sup\limits_{x\in M}d_N(u(x,t),u_0(x))\leqslant c(M)\mu t.\]

If $T'_{\varepsilon}\geqslant T_{m_0}$, then, obviously we have $T_{\varepsilon}\geqslant T_{m_0}$. This means that $T_{\varepsilon}$ can be bounded uniformly from the below with respect to $\varepsilon\in(0,\,1)$.
\medskip

If $T'_{\varepsilon}<T_{m_0}$, then we claim that $$c(M)\mu\cdot T'_{\varepsilon}\geqslant1.$$
In fact, if $u(M,T'_{\varepsilon})\subseteq\Omega$, then, recalling the definition of $T'_{\varepsilon}$, we get that $u(M\times[0,T'_{\varepsilon}])\subseteq\Omega$. Since $M$ is compact, we know from the continuity of $u$ that $u(M\times[0,T'_{\varepsilon}])$ is also compact and $$\mbox{dist}_N(u(M\times[0,T'_{\varepsilon}]), \partial\Omega)>0.$$
Hence, there exists $\delta_0>0$ such that $u(M\times[0,T'_{\varepsilon}+\delta_0])\subseteq\Omega$. So we get a contradiction. It means that there exists a $x_0\in M$ such that $u(x_0,T'_{\varepsilon})\notin\Omega$. This is equivalent to \[\mbox{dist}_N(u(x_0,T'_{\varepsilon}),u_0(M))\geqslant1.\]
It's easy to see that
\[d_N(u(x_0,T'_{\varepsilon}),u_0(x_0))\geqslant1.\]
It follows that $c(M)\mu\cdot T'_{\varepsilon}\geqslant 1$. So the claim is true.

To summarize the above arguments, we obtain
\[T_{\varepsilon}\geqslant T:=\min\Big\{T_{m_0},\frac{1}{c(M)\mu}\Big\}.\]
We conclude that, for all $t\in[0,T]$, there holds true for $k\geq m_0+1$
\[||\nabla u(t)||^2_{H^{k,2}(f)}\leqslant U_k(t).\]
This implies
\[||\nabla u(t)||^2_{H^{k,2}}\leqslant\frac{U_k(t)}{\delta^{k+1}}.\]
It follows immediately from Lemma \ref{lem:equ}(Proposition 2.2 of \cite{DW1})
\[||Du(t)||_{W^{k,2}}\leqslant C(N,k)\sum\limits_{r=1}^{k+1}\Big(\frac{U_k(T_*)}{\delta^{k+1}}\Big)^{\frac{r}{2}}:=\bar{C}(k).\]

Since $\overline{\Omega}$ is compact, there exists a $R_u$ such that for all $t\in[0,T]$,
\[\max_{x\in M}|u(t)|\leqslant R_u,\]
then
\[
||u(t)||_{W^{k+1,2}(M,\,\mathbb{R}^L)}\leqslant\sqrt{|R_u|^2\cdot \mbox{vol}(M)+\bar{C}(k)^2}:=F(k).
\]
This completes proof of Lemma \ref{NSFlemma3.4}.\\

Now we return to prove Theorem \ref{NSFtheorem1.1} and replace $u$ by $u_{\varepsilon}$. Using (\ref{T:1}), we transform the derivatives with respect to time variables $t$ in $||u_{\varepsilon}||^2_{W^{k,2}(M\times[0,T],\, \mathbb{R}^L)}$ into derivatives with respect to space variables $x$. So there exist a positive number $q(k)$, which is large enough, and a function $V_k(\lambda)$, which is monotonously increasing, such that
\[
\aligned
||u_{\varepsilon}||^2_{W^{k,2}(M\times[0,T],\,\mathbb{R}^L)}
\leqslant&\int_0^TV_k(||u_{\varepsilon}(t)||^2_{W^{q(k),2}(M,\,\mathbb{R}^L)})dt\\
\leqslant&T\cdot V_k(F^2(q(k)-1)).
\endaligned
\]
We point out that $V_k(\lambda)$ depends on $f$ but not on $\varepsilon$. Set
\[k(l):=\min\{k\,|\,W^{k,2}(M\times[0,T],\mathbb{R}^L)\hookrightarrow C^l(M\times[0,T],\mathbb{R}^L)\,\,\mbox{compactly}\}.\]
Since $\{u_{\varepsilon}\}$ is a bounded sequence in $W^{k(0),2}(M\times[0,T],\,\mathbb{R}^L)$, then, there exists a subsequence, denoted by $\{u_{\varepsilon_{i,0}}| \varepsilon_{i,0} > 0 \,\,\mbox{and}\,\,\lim_{i\rightarrow +\infty}\varepsilon_{i,0}=0\}\subseteq\{u_{\varepsilon}\}$, and
$$u^0\in W^{k(0),2}(M\times[0,T],\,\mathbb{R}^L)$$
such that $u_{\varepsilon_{i,0}}\rightarrow u^0$ in $C^0(M\times[0,T],\,\mathbb{R}^L)$ as $i\rightarrow+\infty$. We can choose inductively
\[\{u_{\varepsilon_{i,r}}\}\subseteq\{u_{\varepsilon_{i,r-1}}\}\subseteq\cdots\subseteq\{u_{\varepsilon_{i,0}}\}\subseteq\{u_{\varepsilon}\}\]
such that $u_{\varepsilon_{i,r}}\rightarrow u^r$ strongly in $C^r(M\times[0,T],\mathbb{R}^L)$. By the Cantor's diagonal method, we can pick out a subsequence $\{u_{\varepsilon_{r,r}}\}$. For any $\tilde{r}\geqslant2$, $\{u_{\varepsilon_{r,r}}\}$ is a Cauchy sequence in $C^{\tilde{r}}(M\times[0,T],\mathbb{R}^L)$ and tends to $u^{\tilde{r}}:=u$. It is easy to see that $u^{\tilde{r}}=u^{\tilde{r}+1}=u$.

Since $u_{\varepsilon}(M\times[0,T])\subseteq N$ for $\varepsilon\in (0,\, 1)$, we know that $u$ satisfies $u(M\times[0,T])\subseteq N$ and
\begin{equation}\label{E:1}
\left\{
\begin{array}{llll}
\partial_t u=J(u)\tau_f(u),\\
u(\cdot ,0)=u_0.
\end{array}
\right.
\end{equation}

If $u_0$ belongs to $W^{k,2}(M,N)$ instead of $C^{\infty}(M,N)$ and $f$ is in $C^1([0,T_*],C^{k+1}(M))$ instead of $C^{\infty}(M\times[0,T_*])$ where $k\geqslant m_0+1$, we can choose $u_{i0}\in C^{\infty}(M,N)$ and $f_i\in C^{\infty}(M\times[0,T_*])$ such that $u_{i0}\longrightarrow u_0$ in $W^{k,2}(M,\mathbb{R}^L)$ and $f_i\longrightarrow f$ in $C^1([0,T_*],C^{k+1}(M))$. By the fact which we have shown, there exist $T_i>0$ and $u_i\in C^{\infty}(M\times[0,T_i],N)$ satisfying
\begin{equation}\label{E:2}
\left\{
\begin{array}{llll}
\partial_t u_i=J(u_i)\tau_{f_i}(u_i),\\
u_i(\cdot ,0)=u_{i0}.
\end{array}
\right.
\end{equation}
Using Lemma \ref{lem:equ} (Proposition 2.2 in \cite{DW1}), we infer that $\nabla u_{i0}\rightarrow\nabla u_0$ in $H^{k-1,2}(M,N)$. So, we have that $||\nabla u_{i0}||_{H^{k-1,2}}\rightarrow||\nabla u_0||_{H^{k-1,2}}$ and, for sufficiently large $i$, there hold true
$$\delta<f_i<\eta,$$
$$||f_i||_{C^1([0,T_*],C^k(M))}\leqslant||f||_{C^1([0,T_*],C^k(M))}+1,$$
$$||\nabla u_{i0}||_{H^{k-1,2}}^2\leqslant||\nabla u_0||_{H^{k-1,2}}^2+1,$$
$$||\nabla u_{i0}||_{H^{k-1,2}(f_0)}^2\leqslant\eta^k(||\nabla u_0||_{H^{k-1,2}}^2+1).$$
Since the approach by which we deal with (\ref{T:1}) also works for $\varepsilon=0$, for $(\ref{E:2})$ we can also derive the following estimate which is similar with $(\ref{B:1})$
\[
\left\{
\begin{array}{llll}
\aligned
\frac{d}{dt}||\nabla u_i||^2_{H^{m_0,2}(f)}\leqslant&D_{m_0}(t,f_i)(1+||\nabla u_i||^2_{H^{m_0,2}(f)})^{Q(m_0)}\\
\leqslant&D_{m_0}(f_i)(1+||\nabla u_i||^2_{H^{m_0,2}(f)})^{Q(m_0)},\\
||\nabla u_i||^2_{H^{m_0,2}(f)}(0)\leqslant&\eta^{m_0+1}(||\nabla u_0||^2_{H^{m_0,2}}+1).
\endaligned
\end{array}
\right.
\]
Recalling the expression of $D_{m_0}(f_i)$, we find that there exists $D_{m_0}$ such that
\[D_{m_0}\geqslant D_{m_0}(f_i).\]
So,
\[\frac{d}{dt}||\nabla u_i||^2_{H^{m_0,2}(f)}\leqslant D_{m_0}(1+||\nabla u_i||^2_{H^{m_0,2}(f)})^{Q(m_0)}.\]
Then there is $\hat{T}_{m_0}>0$ and a monotonously increasing function $U_{*}\in C^{\infty}([0,\hat{T}_{m_0}])$ satisfying
\[
\left\{
\begin{array}{llll}
\displaystyle\frac{dU_{*}}{dt}=D_{m_0}(1+U_{*})^{Q(m_0)},\\
U_{*}(0)=\eta^{m_0+1}(||\nabla u_0||^2_{H^{m_0,2}}+1).
\end{array}
\right.
\]
Indeed, we can take
\[\hat{T}_{m_0}:=\min\Big\{\frac{1}{2(Q(m_0)-1)\cdot D_{m_0}\cdot[1+\eta^{m_0+1}(1+||\nabla u_0||^2_{H^{m_0,2}})]^{Q(m_0)-1}},\,\,T_*\Big\}\]
and
\[
U_*(t)=\frac{1+\eta^{m_0+1}(1+||\nabla u_0||^2_{H^{m_0,2}})}{\Big(1-(Q(m_0)-1)\cdot D_{m_0}\cdot[1+\eta^{m_0+1}(1+||\nabla u_0||^2_{H^{m_0,2}})]^{Q(m_0)-1}\cdot t\Big)^{\frac{1}{Q(m_0)-1}}}-1.
\]
Therefore, the comparison principle tells us that, for all $t\in[0,\,\min\{T'_i,\,\hat{T}_{m_0}\}]$, there holds true
\[||\nabla u_i(t)||_{H^{m_0,2}(f)}\leqslant\sqrt{U_{*}(t)}\leqslant\sqrt{U_{*}(\hat{T}_{m_0})},\]
where
\[T'_i:=\sup\{t\,|\,u_i(M,[0,t])\subseteq\Omega\}\]
and
\[\Omega=\{y\in N\,|\,\mbox{dist}_{N}(y,u_0(M))<1\}.\]

When $i$ is large enough, we have
\[\sup\limits_{x\in M}d_N(u_{i0}(x),u_0(x))<\frac{1}{2}.\]
So, letting
\[\Omega_i:=\{y\in N\,|\,\mbox{dist}_N(y,u_{i0}(M))<\frac{1}{2}\},\]
we have
\[T''_i:=\sup\{t\,|\,u_i(M,[0,T])\subseteq\Omega_i\}\leqslant T'_i.\]

By the same way as we deal with $$\frac{d}{dt}||\nabla u_{\varepsilon}||^2_{H^{m_0+1,2}(f)},$$
we obtain the following
\[
\left\{
\begin{array}{llll}
\aligned
\frac{d}{dt}||\nabla u_i||^2_{H^{m_0+1,2}(f)}\leqslant& D_{m_0}(1+U_{*}(\hat{T}_{m_0}))^{Q(m_0)}\\
&+\tilde{h}(m_0+1,\Omega,f_i)(1+||\nabla u_i||^2_{H^{m_0+1,2}(f)})(1+U_{*}(\hat{T}_{m_0}))^{G(m_0+1)},\\
||\nabla u_i||^2_{H^{m_0+1,2}(f)}(0)\leqslant&\eta^{m_0+2}(||\nabla u_0||^2_{H^{m_0+1,2}}+1).
\endaligned
\end{array}
\right.
\]
From the expression of $\tilde{h}(m_0+1,\Omega,f_i)$, it is not difficult to see that there exists a positive number $\tilde{h}(m_0+1,\Omega)$ such that $\tilde{h}(m_0+1,\Omega,f_i)
\leqslant\tilde{h}(m_0+1,\Omega)$ for $i$ large enough. So, it follows
\[
\aligned
&\frac{d}{dt}||\nabla u_i||^2_{H^{m_0+1,2}(f)}-D_{m_0}(1+U_{*}(\hat{T}_{m_0}))^{Q(m_0)}\\
\leqslant& \tilde{h}(m_0+1,\Omega)(1+||\nabla u_i||^2_{H^{m_0+1,2}(f)})(1+U_{*}(\hat{T}_{m_0}))^{G(m_0+1)}.
\endaligned
\]
Then, the comparison principle tells us that there exists a monotonously increasing function $\bar{U}_{m_0+1}\in C^{\infty}(\mathbb{R}^1)$ such that, for all $t\in[0,\, \min\{T'_i,\hat{T}_{m_0}\}]$,
\[||\nabla u_i(t)||^2_{H^{m_0+1,2}(f)}\leqslant\bar{U}_{m_0+1}(t).\]
By a similar discussion with the previous, we can infer from the above inequality and the equation (\ref{T:1})
\[
\aligned
||\nabla_tu_i||_{H^{m_0,2}}
\leqslant&\hat{c}(m_0)\cdot\tilde{c}_{m_0}(t,f)\cdot\bar{U}_{m_0+1}(t)\\
\leqslant&\hat{c}(m_0)\cdot\max\limits_{s\in[0,T_*]}\{\tilde{c}_{m_0}(s,f)\}\cdot\bar{U}_{m_0+1}(T_*)
:=&\mu_{m_0+1}.
\endaligned
\]
Hence, by the same argument as in the previous we obtain that there holds true for any $t\in[0,\min\{T''_i,\hat{T}_{m_0}\}]$,
\[\sup\limits_{x\in M}d_N(u_i(x,t),u_{i0}(x))\leqslant c(M)\cdot\mu_{m_0+1}\cdot t.\]
Therefore, we have
\[c(M)\cdot\mu_{m_0+1}\cdot T''_i\geqslant\frac{1}{2}.\]
It follows
\[T_i\geqslant T_0:=\min\left\{\hat{T}_{m_0},\frac{1}{2c(M)\mu_{m_0+1}}\right\}.\]
\medskip

By an inductive method we know that there exists a $\bar{U}_{k-1}\in C^{\infty}(\mathbb{R}^1)$, which is independent of $i$ and is monotonously increasing with respect to the argument $t$, such that for all $t\in[0,T_0]$
\[||\nabla u_i(t)||^2_{H^{k-1,2}(f)}\leqslant\bar{U}_{k-1}(t)\]
and
\begin{equation}\label{E:6}
||\nabla u_i(t)||^2_{H^{k-1,2}}\leqslant\frac{\bar{U}_{k-1}(t)}{\delta^k}\leqslant\frac{\bar{U}_{k-1}(T_*)}{\delta^k}.
\end{equation}
It means that $\{u_i\}$ is a bounded sequence in $L^{\infty}([0,T_0],W^{k,2}(M,N))$. So, there is an $u\in L^{\infty}([0,T_0],W^{k,2}(M,\mathbb{R}^L))$ such that
\[u_i\rightharpoonup u\s\s\mbox{weakly* in}\s\s L^{\infty}([0,T_0],W^{k,2}(M,\mathbb{R}^L)).\]
Since
\[\partial_tu_i=J(u_i)\tau_{f_i}(u_i)\]
and for each $l\geqslant1$
\[D_{a_1}\cdots D_{a_l}D_tu_i=\sum\overrightarrow{B}_{\sigma(\vec{a})}(u_i)(\nabla_{\vec{b}_1}u_i,\cdots,\nabla_{\vec{b}_s}\nabla_tu_i),\]
where $\vec{a}=(a_1,\cdots,a_l)$, $(\vec{b}_1,\cdots,\vec{b}_s)=\sigma(\vec{a})$ is a permutation of $\vec{a}$ and $\overrightarrow{B}_{\sigma(\vec{a})}$ is a multi-linear form on $TN$, we know that there exists an $M(k,\Omega)$ and an $\epsilon(k)$ such that
\[||D_tu_i(t)||_{W^{k-2,2}(M,\mathbb{R}^L)}\leqslant\tilde{c}_{k-2}(t,f_i)M(k,\Omega)\sum\limits_{s=0}^{\epsilon(k)}||\nabla u_i(t)||^s_{H^{k-1,2}(M,N)}.\]

By virtue of $(\ref{E:6})$ and the expression of $\tilde{c}_{k-2}(t,f_i)$, we know that there exists a $\nu_{k-2}(t)$, which is not smaller than $\tilde{c}_{k-2}(t,f_i)$, such that $\{D_tu_i\}$ is a bounded sequence in $L^{\infty}([0,T_0],W^{k-2,2}(M,\mathbb{R}^L))$. By Aubin-Lions Lemma, there exists a subsequence which is still denoted by $u_i$ such that
 $$u_i\rightarrow u \s\s\mbox{strongly in}\s\s L^{\infty}([0,T_0],W^{k-1,2}(M,\mathbb{R}^L)).$$
So, from the Sobolev embedding theorem we have
 $$u_i\rightarrow u \s\s\mbox{strongly in}\s\s L^{\infty}([0,T_0],C^0(M,\mathbb{R}^L)),$$
$u(\cdot,0)=u_0(\cdot)$ and $u([0,T_0]\times M)\subseteq N$.
\medskip

In the next, we will prove that $u$ is a strong solution to NSF. Define
\[
\tilde{J}(y)z:=\left\{
\begin{array}{cccc}
J(y)z, &z\in T_yN,\\
z, &z\in (T_yN)^{\bot}.
\end{array}
\right.
\]
Now $\{u_i\}$ satisfy the following equation:
\[
\left\{
\begin{array}{llll}
D_tv=f_i\cdot J(v)P(v)\Delta v+D_{\beta}f_i\cdot \tilde{J}(v)D_{\beta}v,\\
v(\cdot,0)=u_0.
\end{array}
\right.
\]
Here, $P(y)$ is the orthogonal projection operator from $\mathbb{R}^L$ to $T_yN$.

For any test function $\varphi\in C^{\infty}(M\times[0,T_0],\mathbb{R}^L)$, we have
\[
\aligned
&\Big|\int_0^{T_0}\int_M(f_i\cdot\varphi)J(u_i)P(u_i)\Delta u_i\,dM-\int_0^{T_0}\int_M(f\cdot\varphi)J(u)P(u)\Delta u\,dM\Big|\\
\leqslant&\int_0^{T_0}\int_Mf_i\cdot|\varphi|\cdot|J(u_i)P(u_i)-J(u)P(u)|\cdot|\Delta u_i|\,dM\\
&+\int_0^{T_0}\int_M|\varphi|\cdot|f_i-f|\cdot|J(u)P(u)|\cdot|\Delta u_i|\,dM\\
&+\Big|\int_0^{T_0}\int_M(f\varphi)J(u)P(u)(\Delta u_i-\Delta u)\,dM\Big|\\
\leqslant&\eta||\varphi||_{\infty}\sqrt{\int_0^{T_0}\int_M|\Delta u_i|^2\,dM}\cdot\sqrt{T_0\cdot vol(M)}\cdot||J(u_i)P(u_i)-J(u)P(u)||_{C^0}\\
&+||\varphi||_{\infty}\cdot||f_i-f||_{C^0}\cdot||J(u)P(u)||_{C^0}\cdot\sqrt{T_0\cdot vol(M)}\cdot\sqrt{\int_0^{T_0}\int_M|\Delta u_i|^2\,dM}\\
&+\Big|\int_0^{T_0}\int_M(f\varphi)J(u)P(u)(\Delta u_i-\Delta u)\,dM\Big|.
\endaligned
\]
It follows that, as $i\rightarrow\infty$,
\begin{equation}\label{S:1}
\int_0^{T_0}\int_M(f_i\cdot\varphi)J(u_i)P(u_i)\Delta u_i\,dM\longrightarrow\int_0^{T_0}\int_M(f\cdot\varphi)J(u)P(u)\Delta u\,dM.
\end{equation}
Besides, we also have
\[
\aligned
&\Big|\int_0^{T_0}\int_M(D_{\beta}f_i\cdot\varphi)(J(u_i)D_{\beta}u_i)\,dM-\int_0^{T_0}\int_M(D_{\beta}f\cdot\varphi)(J(u)D_{\beta}u)\,dM\Big|\\
\leqslant&\int_0^{T_0}\int_M|D_{\beta}f_i|\cdot|\varphi|\cdot|\tilde{J}(u_i)-\tilde{J}(u)|\cdot|D_{\beta}u_i|\,dM\\
&+\Big|\int_0^{T_0}\int_M(D_{\beta}f\cdot\varphi)\tilde{J}(u)(D_{\beta}u_i-D_{\beta}u)\,dM\Big|\\
&+\int_0^{T_0}\int_M|D_{\beta}f_i-D_{\beta}f|\cdot|\varphi|\cdot|\tilde{J}(u)|\cdot|D_{\beta}u_i|\,dM\\
\leqslant&||D_{\beta}f_i||_{\infty}||\varphi||_{\infty}||\tilde{J}(u_i)-\tilde{J}(u)||_{\infty}\cdot\sqrt{T_0\cdot vol(M)}\cdot\sqrt{\int_0^{T_0}\int_M|D_{\beta}u_i|^2\,dM}\\
&+\Big|\int_0^{T_0}\int_M(D_{\beta}f\cdot\varphi)\tilde{J}(u)(D_{\beta}u_i-D_{\beta}u)\,dM\Big|\\
&+||D_{\beta}f_i-D_{\beta}f||_{\infty}||\varphi||_{\infty}||\tilde{J}(u)||_{\infty}\cdot\sqrt{T_0\cdot vol(M)}\cdot\sqrt{\int_0^{T_0}\int_M|D_{\beta}u_i|^2\,dM}.
\endaligned
\]
The above inequality implies that
\begin{equation}\label{S:2}
\int_0^{T_0}\int_M(D_{\beta}f_i\cdot\varphi)(J(u_i)D_{\beta}u_i)\,dM\longrightarrow\int_0^{T_0}\int_M(D_{\beta}f\cdot\varphi)(J(u)D_{\beta}u)\,dM.
\end{equation}
Moreover,
\[
\aligned
\int_0^{T_0}\int_M\partial_tu_i\cdot\varphi\,dM&=-\int_0^{T_0}\int_Mu_i\cdot\partial_t\varphi\,dM+\int_Mu_i(T_0)\varphi(T_0)\,dM-\int_Mu_{i0}\cdot\varphi(0)\,dM\\
&\longrightarrow-\int_0^{T_0}\int_Mu\cdot\partial_t\varphi\,dM+\int_Mu(T_0)\varphi(T_0)\,dM-\int_Mu_0\cdot\varphi(0)\,dM\\
&=\int_0^{T_0}\int_M\partial_tu\cdot\varphi\,dM.
\endaligned
\]
Since
\[\int_0^{T_0}\int_M\partial_tu_i\cdot\varphi\,dM=\int_0^{T_0}\int_MJ(u_i)\tau_{f_i}(u_i)\cdot\varphi\,dM,\]
we get from (\ref{S:1}) and (\ref{S:2}) that
\[\int_0^{T_0}\int_M\partial_tu\cdot\varphi\,dM=\int_0^{T_0}\int_MJ(u)\tau_f(u)\cdot\varphi\,dM.\]
This means that $u$ is a strong solution.
\\

Now we are in the position to discuss the uniqueness. Assume that $w$, $v$ both satisfy
\[
\left\{
\begin{array}{llll}
D_tu=f\cdot J(u)P(u)\Delta u+D_{\beta}f\cdot J(u)D_{\beta}u,\\
u(\cdot,0)=u_0.
\end{array}
\right.
\]
Then, we get
\[
\aligned
D_t(w-v)=&f[J(w)P(w)-J(v)P(v)]\Delta w\\
&+f\cdot J(v)P(v)\Delta(w-v)\\
&+D_{\beta}f\cdot(\tilde{J}(w)-\tilde{J}(v))D_{\beta}w\\
&+D_{\beta}f\cdot\tilde{J}(v)D_{\beta}(w-v).
\endaligned
\]
So, applying the above equality we obtain
\begin{equation}\label{U:1}
\aligned
\frac{1}{2}\frac{d}{dt}\int_M|w-v|^2\,dM
=&\int_M\langle w-v,f[J(w)P(w)-J(v)P(v)]\Delta w\rangle\,dM\\
&+\int_M\langle w-v,f\cdot J(v)P(v)\Delta (w-v)\rangle\,dM\\
&+\int_M\langle w-v,D_{\beta}f\cdot(\tilde{J}(w)-\tilde{J}(v))D_{\beta}w\rangle\,dM\\
&+\int_M\langle w-v,D_{\beta}f\cdot\tilde{J}(v)D_{\beta}(w-v)\rangle\,dM\\
:=&I_1+I_2+I_3+I_4.
\endaligned
\end{equation}
Noting that
\[J(w)P(w)-J(v)P(v)=\Big[\int_0^1D(J\cdot P)(v+s(w-v))\,ds\Big](w-v)\]
and
\[\tilde{J}(w)-\tilde{J}(v)=\Big[\int_0^1(D\tilde{J})(v+s(w-v))\,ds\Big](w-v),\]
we have
\begin{equation}\label{U:2}
\aligned
I_1\leqslant&\int_M|w-v|\cdot\eta\cdot\sup\limits_{\substack{|y|\leqslant\max\{||v||_{C^0},||w||_{C^0}\}\\y\in N}}|D(J\cdot P)(y)|\cdot|w-v|\cdot||w||_{C^2}\,dM\\
=&\eta\sup\limits_{\substack{|y|\leqslant\max\{||v||_{C^0},||w||_{C^0}\}\\y\in N}}|D(J\cdot P)(y)|\cdot||w||_{C^2}\int_M|w-v|^2\,dM.
\endaligned
\end{equation}
By integrating by parts and that $J$ is antisymmetric, we derive
\[I_2\leqslant c_2\int_M|w-v|\cdot|D(w-v)|\,dM\]
where $c_2$ depends on $C^1$-norm of $v$ and $C^1$-norm of $f$. Obviously, there hold
\begin{equation}\label{U:3}
I_3\leqslant||Df||_{\infty}||w||_{C^1}\sup\limits_{\substack{|y|\leqslant\max\{||v||_{C^0},||w||_{C^0}\}\\ y\in N}}|D\tilde{J}(y)|\int_M|w-v|^2\,dM
\end{equation}
and
\begin{equation}\label{U:4}
I_4\leqslant c_4\int_M|w-v|\cdot|D(w-v)|\,dM
\end{equation}
where $c_4$ depends upon the $C^0$-norm of $v$ and $C^1$-norm of $f$. Therefore, in view of (\ref{U:2}), (\ref{U:3}) and (\ref{U:4}) we have from (\ref{U:1})
\begin{equation}\label{U:5}
\aligned
&\frac{1}{2}\frac{d}{dt}\int_M|w-v|^2\,dM\\
\leqslant& c_5(\int_M|w-v|^2\,dM+\int_M|D(w-v)|^2\,dM)\\
\leqslant& \frac{c_5}{\delta}(\int_M|w-v|^2\,dM+\int_Mf|D(w-v)|^2\,dM),
\endaligned
\end{equation}
where $c_5$ depends on $C^1$-norm of $v$, $C^1$-norm of $f$ and $C^2$-norm of $w$. However,

\begin{equation}\label{U:6}
\aligned
&\frac{1}{2}\frac{d}{dt}\int_Mf|D(w-v)|^2\,dM\\
=&\frac{1}{2}\int_Mf_t|D(w-v)|^2\,dM+\int_Mf\langle D_{\beta}D_t(w-v),D_{\beta}(w-v)\rangle\,dM\\
=&\frac{1}{2}\int_Mf_t|D(w-v)|^2\,dM-\int_M\langle D_t(w-v),D_{\beta}(fD_{\beta}(w-v))\rangle\,dM\\
=&\frac{1}{2}\int_Mf_t|D(w-v)|^2\,dM\\
&-\int_M\langle f[J(w)P(w)-J(v)P(v)]\Delta w,D_{\beta}(fD_{\beta}(w-v))\rangle\,dM\\
&-\int_M\langle f\cdot J(v)P(v)\Delta(w-v),f\Delta(w-v)+D_{\beta}f\cdot D_{\beta}(w-v)\rangle\,dM\\
&-\int_M\langle D_{\theta}f(\tilde{J}(w)-\tilde{J}(v))D_{\theta}w,D_{\beta}(fD_{\beta}(w-v))\rangle\,dM\\
&-\int_M\langle D_{\theta}f\cdot\tilde{J}(v)D_{\theta}(w-v),f\Delta(w-v)+D_{\beta}f\cdot D_{\beta}(w-v)\rangle\,dM\\
:=&\frac{1}{2}\int_Mf_t|D(w-v)|^2\,dM+I_5+I_6+I_7+I_8.
\endaligned
\end{equation}

For $I_5$, we have
\begin{equation}\label{U:7}
\aligned
I_5&=\int_M\langle D_{\beta}\{f[J(w)P(w)-J(v)P(v)]\Delta w\},\, fD_{\beta}(w-v)\rangle\,dM\\
&\leqslant c'_5\int_M(|D(w-v)|^2+|D(w-v)|\cdot|w-v|)\,dM\\
&\leqslant c''_5\int_M(|D(w-v)|^2+|w-v|^2)\,dM,
\endaligned
\end{equation}
where $c'_5$ and $c''_5$ depends on the $C^1$-norm of $f$, $C^3$-norm of $w$ and $C^0$-norm of $v$. For $I_7$, we have
\begin{equation}\label{U:8}
\aligned
I_7&=\int_M\langle D_{\beta}[D_{\theta}f(\tilde{J}(w)-\tilde{J}(v))D_{\theta}w],\, f D_{\beta}(w-v)\rangle\,dM\\
&\leqslant c'_7\int_M(|D(w-v)|^2+|D(w-v)|\cdot|w-v|)\,dM\\
&\leqslant c''_7\int_M(|D(w-v)|^2+|w-v|^2)\,dM,
\endaligned
\end{equation}
where $c'_7$ and $c''_7$ depends on $C^2$-norm of $f$, $C^2$-norm of $w$ and $C^0$-norm of $v$. On the other side, there holds true
\begin{equation}\label{U:9}
\aligned
I_6+I_8=&-\int_M\langle f\cdot J(v)P(v)\Delta(w-v),\,f\Delta(w-v)\rangle\,dM\\
&-\int_M\langle f\cdot J(v)P(v)\Delta(w-v),\,D_{\beta}f\cdot D_{\beta}(w-v)\rangle\,dM\\
&-\int_M\langle D_{\theta}f\cdot\tilde{J}(v)D_{\theta}(w-v),\,f\Delta(w-v)\rangle\,dM\\
&-\int_M\langle D_{\theta}f\cdot\tilde{J}(v)D_{\theta}(w-v),\,D_{\beta}f\cdot D_{\beta}(w-v)\rangle\,dM.
\endaligned
\end{equation}
Because of that $J(v)$ is antisymmetric, we know that the first integral on the right-hand side of the above equality vanishes. Let
\[B(v):=Id-P(v),\]
where $Id$ is identity operator. Then, for the second term of the right hand side of (\ref{U:9}) we have
\begin{equation}\label{Un:1}
\aligned
&-\int_M\langle f\cdot J(v)P(v)\Delta(w-v),\, D_{\beta}f\cdot D_{\beta}(w-v)\rangle\,dM\\
=&-\int_M\langle f\cdot J(v)P(v)\Delta(w-v),\, D_{\beta}f\cdot P(v)D_{\beta}(w-v)\rangle\,dM,
\endaligned
\end{equation}
for the third term of the right hand side of (\ref{U:9}) we have
\begin{equation}\label{Un:2}
\aligned
&-\int_M\langle D_{\theta}f\cdot\tilde{J}(v)D_{\theta}(w-v),\,f\Delta(w-v)\rangle\,dM\\
=&-\int_M\langle D_{\theta}f\cdot J(v)P(v)D_{\theta}(w-v)+D_{\theta}f\cdot B(v)D_{\theta}w,\, f\Delta(w-v)\rangle\,dM\\
=&-\int_M\langle D_{\theta}f\cdot J(v)P(v)D_{\theta}(w-v),f\cdot P(v)\Delta(w-v)\rangle\,dM\\
&-\int_M\langle D_{\theta}f\cdot B(v)D_{\theta}w,\, f\cdot B(v)\Delta(w-v)\rangle\,dM\\
=&\int_M\langle D_{\theta}f\cdot P(v)D_{\theta}(w-v),\, f\cdot J(v)P(v)\Delta(w-v)\rangle\,dM\\
&-\int_M\langle D_{\theta}f\cdot B(v)D_{\theta}w,\, f\cdot B(v)\Delta(w-v)\rangle\,dM,
\endaligned
\end{equation}
and for the fourth term of the right hand side of (\ref{U:9}) we have
\begin{equation}\label{Un:3}
\aligned
&-\int_M\langle D_{\theta}f\cdot\tilde{J}(v)D_{\theta}(w-v),D_{\beta}f\cdot D_{\beta}(w-v)\rangle\,dM\\
=&-\int_M|D_{\theta}f\cdot B(v)D_{\theta}(w-v)|^2\,dM.
\endaligned
\end{equation}
Combining (\ref{Un:1}),(\ref{Un:2}),(\ref{Un:3}) and (\ref{U:9}) we obtain
\[
\aligned
&I_6+I_8\\
=&-\int_M\langle D_{\theta}f\cdot B(v)D_{\theta}w,f\cdot B(v)\Delta(w-v)\rangle\,dM-\int_M|D_{\theta}f\cdot B(v)D_{\theta}(w-v)|^2\,dM\\
=&-\int_M\langle D_{\theta}f\cdot B(v)D_{\theta}w,f\Delta(w-v)\rangle\,dM-\int_M|D_{\theta}f\cdot B(v)D_{\theta}(w-v)|^2\,dM\\
=&-\int_M\langle D_{\theta}f\cdot (B(v)-B(w))D_{\theta}w,f\Delta(w-v)\rangle\,dM-\int_M|D_{\theta}f\cdot B(v)D_{\theta}(w-v)|^2\,dM.
\endaligned
\]
By integrating by parts, the first term on the right-hand side of the above equality is not bigger than
\[c'_8\int_M(|v-w|\cdot|Dv-Dw|+|Dv-Dw|^2)\,dM.\]
By the H\"older inequality, we infer the above integral is less than or equal to
\[c''_8\int_M(|v-w|^2+|Dv-Dw|^2)\,dM\]
where $c'_8$ and $c''_8$ depends on $C^2$-norm of $f$, $C^2$-norm of $w$ and $C^1$-norm of $v$. Thus, we obtain
\begin{equation}\label{Un:4}
I_6+I_8\leqslant c''_8\int_M(|v-w|^2+|Dv-Dw|^2)\,dM.
\end{equation}

Combing (\ref{U:7}), (\ref{U:8}), (\ref{Un:4}) and (\ref{U:6}) we obtain
\begin{equation}\label{Un:5}
\aligned
\frac{1}{2}\frac{d}{dt}\int_Mf|Dw-Dv|^2\,dM
\leqslant&c_9\int_M(|w-v|^2+|Dw-Dv|^2)\,dM\\
\leqslant&c'_9\int_M(|w-v|^2+f|Dw-Dv|^2)\,dM,
\endaligned
\end{equation}
where $c_9$ and $c'_9$ depends upon the $C^0$-norm of $f_t$, $C^2$-norm of $f$, $C^3$-norm of $w$ and $C^1$-norm of $v$. Letting
\[\alpha(t):=\int_M(|w-v|^2+f|Dw-Dv|^2)\,dM,\]
we get from (\ref{U:5}) and (\ref{Un:5})
\[\frac{d\alpha}{dt}(t)\leqslant c(t)\alpha(t),\]
where $$c(t):=2\max\{c'_9,\frac{c_5}{\delta}\}.$$ Then, it follows
\[\alpha(t)\leqslant\alpha(0)\exp\Big\{\int_0^tc(s)ds\Big\}.\]
Since $\alpha\geqslant0$ and $\alpha(0)=0$, it follows that $\alpha(t)\equiv0$ which means that $w\equiv v$.\endproof
\\

\textbf{Proof of Theorem \ref{NSFtheorem1.2}.} Since $f$ is periodic with respect to space variables, there exist $\delta\in(0,1]$ and $\eta\in[1,\infty)$ such that $\delta<f(x,t)<\eta$. We define
\[
\Omega_k:=[-kR_1,\, kR_1]\times\cdots\times[-kR_m,\, kR_m]
\]
where $k\in\mathbb{N}^+$. Letting $f_k:=f|_{\Omega_k}$, we get that $f_k$ can be regarded as a function defined on $\widetilde{\mathbb{T}}^m_k$ (that is to say, $f_k$ belongs to $C^1([0,T_*],C^{r+1}(\widetilde{\mathbb{T}}^m_k))$) and $\delta<f_k<\eta$.

By Lemma 3.4 of \cite{DW1}, if $r>\frac{m}{2}$, then there exist $\{u_{i0}\}\subseteq W^{r,2}(\mathbb{R}^m,N)\cap C_0^{\infty}(\mathbb{R}^m,\mathbb{R}^L)$ such that $u_{i0}\longrightarrow u_0$ strongly in $W^{r,2}(\mathbb{R}^m,N)$. We assume that the support of $u_{i0}$ lies in $\Omega_{k_i}$ for some $k_i$ and $k_i\nearrow\infty$ as $i\longrightarrow\infty$. So $u_{i0}$ can be regarded as a smooth map from $\widetilde{\mathbb{T}}^m_{k_i}$ into $N$.

Let us consider the following Cauchy problem:
\[
\left\{
\begin{array}{llll}
\partial_t u=J(u)\tau_{f_{k_i}}(u)\\
u(\cdot ,0)=u_{i0}\in C^{\infty}(\widetilde{\mathbb{T}}^m_{k_i},N)
\end{array}
\right.
\]
By Theorem \ref{NSFtheorem1.1}, there is a $u_i\in L^{\infty}([0,T_i],W^{r,2}(\widetilde{\mathbb{T}}^m_{k_i},N))$, for some $T_i>0$, being the solution of above system. Taking $\varepsilon=0$ in Lemma \ref{NSFlemma3.1}, Lemma \ref{NSFlemma3.2} and Lemma \ref{NSFlemma3.3}, by Lemma \ref{NSFlemma3.4} we know that
\[
T_i\geqslant T(f_{k_i},m_0,||\nabla u_{i0}||^2_{H^{m_0,2}(f_{k_i0})},T_*,\widetilde{\mathbb{T}}^m_{k_i},N)>0
\]
and for all $t\in[0,T_i]$,
\[
||u_i(t)||_{W^{r,2}(\widetilde{\mathbb{T}}^m_{k_i},\mathbb{R}^L)}\leqslant C(\widetilde{\mathbb{T}}^m_{k_i},N,f_{k_i},r-1,T_*,||\nabla u_{i0}||_{H^{r-1,2}(f_{k_i0})}):=F_{k_i}(r-1)
\]
where $f_{k_i0}(x):=f_{k_i}(x,0)$. By Lemma \ref{NSFlemma2.3} one may see that the constant in (\ref{equ:int}) does not depend on the diameter of $\widetilde{\mathbb{T}}^m_{k_i}$. So the above quantities
$$T(f_{k_i},m_0,||\nabla u_{i0}||^2_{H^{m_0,2}(f_{k_i0})},T_*,\widetilde{\mathbb{T}}^m_{k_i},N)$$
and $F_{k_i}(r-1)$ are independent of the diameter of $\widetilde{\mathbb{T}}^m_{k_i}$. By the expression of
$$T(f_{k_i},m_0,||\nabla u_{i0}||^2_{H^{m_0,2}(f_{k_i0})},T_*,\widetilde{\mathbb{T}}^m_{k_i},N)$$
and $F_{k_i}(r-1)$, one can see easily that $T(f_{k_i},m_0,||\nabla u_{i0}||^2_{H^{m_0,2}(f_{k_i0})},T_*,\widetilde{\mathbb{T}}^m_{k_i},N)$ and $F_{k_i}(r-1)$ depend continuously on $\partial_tf_{k_i}$, $f_{k_i}$, $\nabla f_{k_i},\cdots$, $\nabla^{r+1}f_{k_i}$, $||\nabla u_{i0}||^2_{H^{m_0,2}(f_{k_i0})}$(or $||\nabla u_{i0}||_{H^{r-1,2}(f_{k_i0})}$). Since $f$ is periodic with respect to space variables, all of $|\partial_tf|, f, |\nabla f|,\cdots, |\nabla^{r+1}f|$ are controlled by some $\tilde{M}\in(0,\infty)$. By the definition of $f_k$ it is easy to know for all $i$, $|\partial_tf_{k_i}|$, $f_{k_i}$, $|\nabla f_{k_i}|,\cdots$, $|\nabla^{r+1}f_{k_i}|$ are bounded by $\tilde{M}$.

Since $u_{i0}\longrightarrow u_0$ in $W^{r,2}(\mathbb{R}^m,N)$ and $r>\frac{m}{2}$, applying Lemma \ref{lem:equ} we get
\[
||\nabla u_{i0}-\nabla u_0||_{H^{r-1,2}}\longrightarrow0\s as\s i\longrightarrow\infty.
\]
So
\[
||\nabla u_{i0}-\nabla u_0||_{H^{r-1,2}(f_0)}\longrightarrow0.
\]
Moreover,
\[
||\nabla u_{i0}||_{H^{r-1,2}(f_{k_i0})}\longrightarrow||\nabla u_0||_{H^{r-1,2}(f_0)}\s as\s i\longrightarrow\infty.
\]
Combining above discussion, there are $T_{m_0}>0$ and $F(r-1)<\infty$ such that
\[
T(f_{k_i},m_0,||\nabla u_{i0}||^2_{H^{m_0,2}(f_{k_i0})},T_*,\widetilde{\mathbb{T}}^m_{k_i},N)\geqslant T_{m_0}
\]
and for all $t\in[0,T_{m_0}]$,
\[
||u_i(t)||_{W^{r,2}(\widetilde{\mathbb{T}}^m_{k_i},\mathbb{R}^L)}\leqslant F(r-1).
\]

We regard each $u_i$ as a map from $\Omega_{k_i}\times[0,T_{m_0}]$ into $N$. By the proof of  Theorem 1.2 of \cite{DW1} we know that there exists a $u\in L^{\infty}([0,T_{m_0}],W^{r,2}(\mathbb{R}^m,N))$ such that for any compact domain $\mathcal{C}\subseteq\mathbb{R}^m$,
\[
u_i\rightharpoonup u\s weakly*\s in\s L^{\infty}([0,T_{m_0}],W^{r,2}(\mathcal{C},N))
\]
upon extracting a subsequence and re-indexing if necessary. Because $W^{r,2}(\mathcal{C},N)\hookrightarrow C^2(\mathcal{C},N)$ compactly, by Aubin-Lions Lemma, there exists a subsequence which is still denoted by $\{u_i\}$ such that $u_i\longrightarrow u$ strongly in $L^{\infty}([0,T_{m_0}],C^2(\mathcal{C},N))$. Therefore,
\[
J(u_i)\tau_{f_{k_i}}(u_i)=\tilde{J}(u_i)\tau_{f_{k_i}}(u_i)\longrightarrow\tilde{J}(u)\tau_f(u)=J(u)\tau_f(u)\s strongly\s in\s L^{\infty}([0,T_{m_0}]\times\mathcal{C},\mathbb{R}^L).
\]
It is easy to see that $u$ is a strong solution to the following Cauchy problem:
\[
\left\{
\begin{array}{llll}
\partial_t u=J(u)\tau_f(u)\s on\s \mathbb{R}^m\times(0,T_{m_0}]\\
u(\cdot ,0)=u_0\in W^{r,2}(\mathbb{R}^m,N).
\end{array}
\right.
\]
\\

The uniqueness will be stated in the next section. As for smoothness, it is easy to see that if $u_0\in\mathscr{H}$ and $f\in C^1([0,T_*],C^{\infty}(\mathbb{R}^m))$, then $u\in L^{\infty}([0,T_{m_0}],\mathscr{H})$. In the next, we consider \[\underbrace{\nabla_t\cdots\nabla_t}\limits_{k-times}u.\]
Using the above equation inductively to transform derivatives about time variable into derivatives about space variables, we can know that
\[\underbrace{\nabla_t\cdots\nabla_t}\limits_{k-times}u\in L^{\infty}([0,T_{m_0}],\mathscr{H}).\]
This completes the proof.\endproof\\

\section{Uniqueness}
In this section, we always assume $N$ is a complete manifold with curvature bounded by $K_0$ and injectivity radius bounded from below by $i_0>0$. First, we need the following lemma.

\begin{lem}
Under the assumptions of Theorem \ref{thm0.1} or Theorem \ref{thm0.2}, there exists a $T'>0$ such that $d(u_1(t,x),u_2(t,x))<\delta_0$ for any $(t,x)\in[0,T']\times M$, where $\delta_0:=\min\{\frac{i_0}{2},\frac{1}{4\sqrt{K_0}}\}.$
\end{lem}
\textbf{Proof.} Since
$$d(u_1,u_2)\leqslant d(u_1,u_0)+d(u_2,u_0),$$
we only need to prove that $u_{\lambda}$ stays sufficently close to $u_0$ in $[0,T']$ where $\lambda=1,2$.

\textbf{Case1:} If $f$ and $u_{\lambda}$ satisfies the assumptions of Theorem \ref{thm0.1}, then
$$\tau_f(u_{\lambda})=f\tau(u_{\lambda})+\nabla f\cdot du_{\lambda}\in L^{\infty}([0,T]\times M).$$
Let $\gamma_{\lambda}(s):=u_{\lambda}(s,x)$ where $s\in[0,t]$. This is a curve connected $u_0(x)$ and $u_{\lambda}(t,x)$, and
\[
\aligned
&d(u_{\lambda}(t,x),u_0(x))\\
\leqslant&\int_0^t\Big|\frac{d\gamma_{\lambda}}{dt}(s)\Big|\,ds=\int_0^t|\partial_tu_{\lambda}(s,x)|\,ds=\int_0^t|J(u_{\lambda}(s,x))\tau_{f(s,x)}(u_{\lambda}(s,x))|\,ds\\
\leqslant&t\cdot||\tau_f(u_{\lambda})||_{L^{\infty}([0,T]\times M)}\leqslant T'\cdot||\tau_f(u_{\lambda})||_{L^{\infty}([0,T]\times M)}.
\endaligned
\]
By choosing
$$T':=\frac{\delta_0}{4\max\{||\tau_f(u_{\lambda})||_{L^{\infty}([0,T]\times M)}|\lambda=1,2\}}$$
we obtain the required result.

\textbf{Case2:} If $f$ and $u_{\lambda}$ satisfies the assumptions of Theorem \ref{thm0.2}, we embed $N$ into Euclidean space $\mathbb{R}^L$. Then
\[
\aligned
&\frac{d}{dt}||u_{\lambda}(t)-u_0||^2_2\\
=&2\int_M\langle u_{\lambda}(t)-u_0,\partial_tu_{\lambda}(t)\rangle\,dM\leqslant2||u_{\lambda}(t)-u_0||_2\cdot||\partial_tu_{\lambda}(t)||_2\\
=&2||u_{\lambda}(t)-u_0||_2\cdot||J(u_{\lambda})\tau_f(u_{\lambda})(t)||_2=2||u_{\lambda}(t)-u_0||_2\cdot||\tau_f(u_{\lambda})(t)||_2.
\endaligned
\]
So we get
\[
\frac{d}{dt}||u_{\lambda}(t)-u_0||^2_2\leqslant C_1,
\]
where
$$C_1:=\sup\Big\{2(||u_{\lambda}(t)||_2+||u_0||_2)\cdot||\tau_f(u_{\lambda})(t)||_2\Big|t\in[0,T],\s\lambda=1,2\Big\}.$$
Since $u_{\lambda}(0,\cdot)=u_0$,
\[
||u_{\lambda}(t,\cdot)-u_0||_2\leqslant \sqrt{C_1\cdot t}.
\]
By Theorem 5 in \cite{Cantor}, we have
\[
||u_{\lambda}(t)-u_0||_{\infty}\leqslant C_2||u_{\lambda}(t)-u_0||_2^a\cdot||u_{\lambda}(t)-u_0||^{1-a}_{W^{[\frac{m}{2}]+1,2}}\leqslant C_3\cdot t^{a/2}
\]
where
\[a=1-\frac{m}{2([\frac{m}{2}]+1)}\in(0,1)\]
and
\[C_3:=C_2\cdot C_1^{a/2}\Big[\sup\Big\{||u_{\lambda}(t)||_{W^{[\frac{m}{2}]+1,2}}+||u_0||_{W^{[\frac{m}{2}]+1,2}}\Big|t\in[0,T],\s\lambda=1,2\Big\}\Big]^{1-a}.\]
Following the proof of Lemma 3.1 in \cite{SW}, we can see that, since $N$ has bounded geometry, there holds
\[
d(u_{\lambda}(t),u_0)\leqslant C_4||u_{\lambda}(t)-u_0||_{\infty}\leqslant C_4C_3\cdot t^a\leqslant C_4C_3\cdot(T')^a.
\]
We only need to choose $T'=(\frac{\delta_0}{4C_3C_4})^{\frac{1}{a}}$ to complete the proof.\endproof
\\

\textbf{Proof of Theorem \ref{thm0.1}}

For $\lambda=1,2$, we define
$$\nabla_{\lambda}u_{\lambda}:=\frac{\partial u^{\alpha}_{\lambda}}{\partial x^i}dx^i\otimes\frac{\partial}{\partial y^{\alpha}}$$
and
$$\nabla_{\lambda,t}u_{\lambda}:=\frac{\partial u_{\lambda}^{\alpha}}{\partial t}\frac{\partial}{\partial y^{\alpha}}.$$

Following the idea to estimate $Q_1$ in \cite{SW}, we can construct a global bundle morphism
$$\mathcal{P}: u_2^*TN\longrightarrow u^*_1TN$$
and extend it to a bundle morphism from $u_2^*TN\otimes T^*M$ to $u_1^*TN\otimes T^*M$. Let $\tilde{\nabla}:=\nabla\oplus\nabla$ be the covariant derivative on $N\times N$ induced by the Levi-Civita connection $\nabla$ on $N$. Using the NSF equation and integrating by parts, we have
\[
\aligned
&\frac{d}{dt}\int_Md^2(u_1,u_2)\,dM\\
=&\int_M\langle\tilde{\nabla}d^2,(\nabla_{1,t}u_1,\nabla_{2,t}u_2)\rangle\,dM\\
=&\int_M\langle\tilde{\nabla}d^2,(J(u_1)\tau_f(u_1),J(u_2)\tau_f(u_2))\rangle\,dM\\
=&-\int_Mf(\tilde{\nabla}^2d^2)\Big((\nabla_1u_1,\nabla_2u_2),(J(u_1)\nabla_1u_1,J(u_2)\nabla_2u_2)\Big)\,dM\\
\leqslant&\eta\int_M\Big|(\tilde{\nabla}^2d^2)\Big((\nabla_1u_1,\nabla_2u_2),(J(u_1)\nabla_1u_1,J(u_2)\nabla_2u_2)\Big)\Big|\,dM.
\endaligned
\]
In view of \textbf{Estimate of $Q_1$} in \cite{SW}, we have
\[\frac{1}{2}\Big|(\tilde{\nabla}^2d^2)\Big((\nabla_1u_1,\nabla_2u_2),(J(u_1)\nabla_1u_1,J(u_2)\nabla_2u_2)\Big)\Big|\leqslant|\mathcal {P}\nabla_2u_2-\nabla_1u_1|^2+C_5\cdot d^2(u_1,u_2),\]
where $C_5$ depends on $L^{\infty}$-norm of $\nabla_{\lambda}u_{\lambda}$. So
\begin{equation}\label{Formula1}
\aligned
&\frac{d}{dt}\int_Md^2(u_1,u_2)\,dM\\
\leqslant&2\eta\cdot C_5\int_Md^2(u_1,u_2)\,dM+2\eta\int_M|\mathcal{P}\nabla_2u_2-\nabla_1u_1|^2\,dM\\
\leqslant&2\eta\cdot C_5\int_Md^2(u_1,u_2)\,dM+\frac{2\eta}{\delta}\int_Mf|\mathcal{P}\nabla_2u_2-\nabla_1u_1|^2\,dM.
\endaligned
\end{equation}
Review that
\[\nabla_{\lambda,t}u_{\lambda}=J(u_{\lambda})\tau_f(u_{\lambda}).\]
Differentiating the equation, we get
\begin{equation}\label{Formula2}
\aligned
\nabla_{\lambda,t}\nabla_{\lambda,k}u_{\lambda}=J(u_{\lambda})&\Big[g^{ij}\nabla_j\nabla_kf\cdot\nabla_{\lambda,i}u_{\lambda}+g^{ij}\nabla_jf\cdot\nabla_{\lambda,i}\nabla_{\lambda,k}u_{\lambda}\\
&+\nabla_kf\cdot g^{ij}\nabla_{\lambda,i}\nabla_{\lambda,j}u_{\lambda}+f\cdot g^{ij}\nabla_{\lambda,j}\nabla_{\lambda,i}\nabla_{\lambda,k}u_{\lambda}\\
&-f\nabla_{\lambda,h}u_{\lambda}\cdot g^{hs}Ric^M_{sk}-f\cdot g^{ij}R^N(\nabla_{\lambda,j}u_{\lambda},\nabla_{\lambda,k}u_{\lambda})\nabla_{\lambda,i}u_{\lambda}\Big].
\endaligned
\end{equation}
Recall that in Preliminary, we have chosen a frame $\{f_{\alpha}\}^n_{\alpha=1}$ such that $J$ is reduced to a constant skew-symmetric matrix $J_0$. Using (\ref{Formula2}), one can obtain
\[
\aligned
&\frac{d}{dt}\Big(\frac{1}{2}\int_Mf|\mathcal{P}\nabla_2u_2-\nabla_1u_1|^2\,dM\Big)\\
=&\frac{1}{2}\int_Mf_t|\mathcal{P}\nabla_2u_2-\nabla_1u_1|^2\,dM+\int_Mg^{kl}\langle(\nabla_{1,t}\mathcal{P}-\mathcal{P}\nabla_{2,t})\nabla_{2,k}u_2,f(\mathcal{P}\nabla_{2,l}u_2-\nabla_{1,l}u_1)\rangle\,dM\\
&+I_1+I_2+I_3-\int_Mf^2\cdot g^{kl}\langle \mathcal{P}\nabla_{2,l}u_2-\nabla_{1,l}u_1,g^{hs}Ric^M_{sk}J_0(\mathcal{P}\nabla_{2,h}u_2-\nabla_{1,h}u_1)\rangle\,dM\\
&-\int_Mf^2\cdot g^{kl}g^{ij}\langle\mathcal{P}\nabla_{2,l}u_2-\nabla_{1,l}u_1,J_0[\mathcal{P}R^N(\nabla_{2,j}u_2,\nabla_{2,k}u_2)\nabla_{2,i}u_2\\
&-R^N(\nabla_{1,j}u_1,\nabla_{1,k}u_1)\nabla_{1,i}u_1]\rangle\,dM
\endaligned
\]
where
\[I_1:=\int_Mg^{kl}g^{ij}f\cdot\nabla_kf\langle\mathcal{P}\nabla_{2,l}u_2-\nabla_{1,l}u_1,J_0(\mathcal{P}\nabla_{2,i}\nabla_{2,j}u_2-\nabla_{1,i}\nabla_{1,j}u_1)\rangle\,dM,\]
\[I_2:=\int_Mg^{kl}g^{ij}f\cdot\nabla_jf\langle\mathcal{P}\nabla_{2,l}u_2-\nabla_{1,l}u_1,J_0(\mathcal{P}\nabla_{2,i}\nabla_{2,k}u_2-\nabla_{1,i}\nabla_{1,k}u_1)\rangle\,dM\]
and
\[I_3:=\int_Mg^{kl}g^{ij}f^2\langle\mathcal{P}\nabla_{2,l}u_2-\nabla_{1,l}u_1,J_0(\mathcal{P}\nabla_{2,j}\nabla_{2,i}\nabla_{2,k}u_2-\nabla_{1,j}\nabla_{1,i}\nabla_{1,k}u_1)\rangle\,dM\]
Note that
\[
g^{kl}g^{hs}Ric^M_{sk}\langle \mathcal{P}\nabla_{2,l}u_2-\nabla_{1,l}u_1,J_0(\mathcal{P}\nabla_{2,h}u_2-\nabla_{1,h}u_1)\rangle=0
\]
and
\[
\aligned
I_1=\int_M&g^{kl}g^{ij}f\cdot\nabla_kf\langle\mathcal{P}\nabla_{2,l}u_2-\nabla_{1,l}u_1,\\
&J_0(\mathcal{P}\nabla_{2,i}\nabla_{2,j}u_2-\nabla_{1,i}\mathcal{P}\nabla_{2,j}u_2+\nabla_{1,i}\mathcal{P}\nabla_{2,j}u_2-\nabla_{1,i}\nabla_{1,j}u_1)\rangle\,dM\\
=\int_M&g^{kl}g^{ij}f\cdot\nabla_kf\langle\mathcal{P}\nabla_{2,l}u_2-\nabla_{1,l}u_1,J_0(\mathcal{P}\nabla_{2,i}-\nabla_{1,i}\mathcal{P})\nabla_{2,j}u_2\rangle d\,M\\
-&\int_Mg^{kl}g^{ij}\langle\nabla_{1,i}[\nabla_kf\cdot f(\mathcal{P}\nabla_{2,l}u_2-\nabla_{1,l}u_1)],J_0(\mathcal{P}\nabla_{2,j}u_2-\nabla_{1,j}u_1)\rangle d\,M\\
=\int_M&g^{kl}g^{ij}f\cdot\nabla_kf\langle\mathcal{P}\nabla_{2,l}u_2-\nabla_{1,l}u_1,J_0(\mathcal{P}\nabla_{2,i}-\nabla_{1,i}\mathcal{P})\nabla_{2,j}u_2\rangle d\,M\\
-&\int_Mg^{kl}g^{ij}\langle(\nabla_i\nabla_kf\cdot f+\nabla_kf\cdot\nabla_if)(\mathcal{P}\nabla_{2,l}u_2-\nabla_{1,l}u_1),J_0(\mathcal{P}\nabla_{2,j}u_2-\nabla_{1,j}u_1)\rangle d\,M\\
-&\int_Mg^{kl}g^{ij}\nabla_kf\cdot f\langle\nabla_{1,i}\mathcal{P}\nabla_{2,l}u_2-\nabla_{1,i}\nabla_{1,l}u_1,J_0(\mathcal{P}\nabla_{2,j}u_2-\nabla_{1,j}u_1)\rangle d\,M\\
=\int_M&g^{kl}g^{ij}f\cdot\nabla_kf\langle\mathcal{P}\nabla_{2,l}u_2-\nabla_{1,l}u_1,J_0(\mathcal{P}\nabla_{2,i}-\nabla_{1,i}\mathcal{P})\nabla_{2,j}u_2\rangle d\,M\\
-&\int_Mg^{kl}g^{ij}\nabla_kf\cdot f\langle(\nabla_{1,i}\mathcal{P}-\mathcal{P}\nabla_{2,i})\nabla_{2,l}u_2,J_0(\mathcal{P}\nabla_{2,j}u_2-\nabla_{1,j}u_1)\rangle d\,M\\
-&\int_Mg^{kl}g^{ij}\nabla_kf\cdot f\langle\mathcal{P}\nabla_{2,i}\nabla_{2,l}u_2-\nabla_{1,i}\nabla_{1,l}u_1,J_0(\mathcal{P}\nabla_{2,j}u_2-\nabla_{1,j}u_1)\rangle d\,M
\endaligned
\]
Since $J_0$ is skew-symmetric, we obtain
\begin{eqnarray*}
I_2=-\int_Mg^{kl}g^{ij}f\cdot\nabla_jf\langle J_0(\mathcal{P}\nabla_{2,l}u_2-\nabla_{1,l}u_1),\mathcal{P}\nabla_{2,i}\nabla_{2,k}u_2-\nabla_{1,i}\nabla_{1,k}u_1\rangle\,dM.
\end{eqnarray*}
And integrating by parts, we have
\begin{equation}\label{Formula3}
\aligned
I_3=&\int_Mg^{kl}g^{ij}f^2\langle\mathcal{P}\nabla_{2,l}u_2-\nabla_{1,l}u_1,J_0(\mathcal{P}\nabla_{2,j}-\nabla_{1,j}\mathcal{P})\nabla_{2,i}\nabla_{2,k}u_2\rangle\,dM\\
&-\int_Mg^{kl}g^{ij}\langle2f\cdot\nabla_jf(\mathcal{P}\nabla_{2,l}u_2-\nabla_{1,l}u_1),J_0(\mathcal{P}\nabla_{2,i}\nabla_{2,k}u_2-\nabla_{1,i}\nabla_{1,k}u_1)\rangle\,dM\\
&-\int_Mg^{kl}g^{ij}f^2\langle \nabla_{1,j}(\mathcal{P}\nabla_{2,l}u_2-\nabla_{1,l}u_1),J_0(\mathcal{P}\nabla_{2,i}\nabla_{2,k}u_2-\nabla_{1,i}\nabla_{1,k}u_1)\rangle\,dM\\
=&\int_Mg^{kl}g^{ij}f^2\langle\mathcal{P}\nabla_{2,l}u_2-\nabla_{1,l}u_1,J_0(\mathcal{P}\nabla_{2,j}-\nabla_{1,j}\mathcal{P})\nabla_{2,i}\nabla_{2,k}u_2\rangle\,dM\\
&+\int_M2f\cdot\nabla_jf\cdot g^{kl}g^{ij}\langle J_0(\mathcal{P}\nabla_{2,l}u_2-\nabla_{1,l}u_1),\mathcal{P}\nabla_{2,i}\nabla_{2,k}u_2-\nabla_{1,i}\nabla_{1,k}u_1\rangle\,dM\\
&-\int_Mg^{kl}g^{ij}f^2\langle (\nabla_{1,j}\mathcal{P}-\mathcal{P}\nabla_{2,j})\nabla_{2,l}u_2,J_0(\mathcal{P}\nabla_{2,i}\nabla_{2,k}u_2-\nabla_{1,i}\nabla_{1,k}u_1)\rangle\,dM\\
&-\int_Mg^{kl}g^{ij}f^2\langle\mathcal{P}\nabla_{2,j}\nabla_{2,l}u_2-\nabla_{1,j}\nabla_{1,l}u_1,J_0(\mathcal{P}\nabla_{2,i}\nabla_{2,k}u_2-\nabla_{1,i}\nabla_{1,k}u_1)\rangle\,dM
\endaligned
\end{equation}
Note that the last term on the right-hand side of the last equality sign of (\ref{Formula3}) vanishes. Moreover,
\[
\aligned
&\int_Mg^{kl}g^{ij}f^2\langle (\nabla_{1,j}\mathcal{P}-\mathcal{P}\nabla_{2,j})\nabla_{2,l}u_2,J_0(\mathcal{P}\nabla_{2,i}\nabla_{2,k}u_2-\nabla_{1,i}\nabla_{1,k}u_1)\rangle\,dM\\
=&\int_Mg^{kl}g^{ij}f^2\langle (\nabla_{1,j}\mathcal{P}-\mathcal{P}\nabla_{2,j})\nabla_{2,l}u_2,\\
&J_0[(\mathcal{P}\nabla_{2,i}-\nabla_{1,i}\mathcal{P})\nabla_{2,k}u_2+\nabla_{1,i}(\mathcal{P}\nabla_{2,k}u_2-\nabla_{1,k}u_1)]\rangle\,dM\\
=&\int_Mg^{kl}g^{ij}f^2\langle (\nabla_{1,j}\mathcal{P}-\mathcal{P}\nabla_{2,j})\nabla_{2,l}u_2,\nabla_{1,i}[J_0(\mathcal{P}\nabla_{2,k}u_2-\nabla_{1,k}u_1)]\rangle\,dM\\
=&-\int_Mg^{kl}g^{ij}\langle2f\cdot\nabla_if (\nabla_{1,j}\mathcal{P}-\mathcal{P}\nabla_{2,j})\nabla_{2,l}u_2+f^2\nabla_{1,i}[(\nabla_{1,j}\mathcal{P}-\mathcal{P}\nabla_{2,j})\nabla_{2,l}u_2],\\
&J_0(\mathcal{P}\nabla_{2,k}u_2-\nabla_{1,k}u_1)\rangle\,dM
\endaligned
\]
Combining the above discussion, we get
\[
\aligned
&\frac{d}{dt}(\frac{1}{2}\int_Mf|\mathcal{P}\nabla_2u_2-\nabla_1u_1|^2\,dM)\\
=&\frac{1}{2}\int_Mf_t|\mathcal{P}\nabla_2u_2-\nabla_1u_1|^2\,dM+\int_Mf\cdot g^{kl}\langle(\nabla_{1,t}\mathcal{P}-\mathcal{P}\nabla_{2,t})\nabla_{2,k}u_2,\mathcal{P}\nabla_{2,l}u_2-\nabla_{1,l}u_1\rangle\,dM\\
&+\int_Mf\cdot\nabla_kf\cdot g^{ij} g^{kl}\langle J_0(\mathcal{P}\nabla_{2,i}-\nabla_{1,i}\mathcal{P})\nabla_{2,j}u_2,\mathcal{P}\nabla_{2,l}u_2-\nabla_{1,l}u_1\rangle\,dM\\
&+\int_Mf\cdot\nabla_kf\cdot g^{ij} g^{kl}\langle J_0(-\mathcal{P}\nabla_{2,i}+\nabla_{1,i}\mathcal{P})\nabla_{2,l}u_2,\mathcal{P}\nabla_{2,j}u_2-\nabla_{1,j}u_1\rangle\,dM\\
&+\int_Mf^2\cdot g^{ij}g^{kl}\langle J_0(\mathcal{P}\nabla_{2,j}-\nabla_{1,j}\mathcal{P})\nabla_{2,i}\nabla_{2,k}u_2,\mathcal{P}\nabla_{2,l}u_2-\nabla_{1,l}u_1\rangle\,dM\\
&+\int_M2f\cdot\nabla_if\cdot g^{ij}g^{kl}\langle (-\mathcal{P}\nabla_{2,j}+\nabla_{1,j}\mathcal{P})\nabla_{2,l}u_2,J_0(\mathcal{P}\nabla_{2,k}u_2-\nabla_{1,k}u_1)\rangle\,dM\\
&+\int_Mf^2\cdot g^{ij}g^{kl}\langle (-\mathcal{P}\nabla_{2,i}\nabla_{2,j}+\nabla_{1,i}\nabla_{1,j}\mathcal{P})\nabla_{2,l}u_2,J_0(\mathcal{P}\nabla_{2,k}u_2-\nabla_{1,k}u_1)\rangle\,dM\\
&+\int_Mf^2\cdot g^{ij}g^{kl}\langle (\mathcal{P}\nabla_{2,i}-\nabla_{1,i}\mathcal{P})\nabla_{2,j}\nabla_{2,l}u_2,J_0(\mathcal{P}\nabla_{2,k}u_2-\nabla_{1,k}u_1)\rangle\,dM\\
&-\int_Mf^2\cdot g^{ij}g^{kl}\langle\mathcal{P}\nabla_{2,l}u_2-\nabla_{1,l}u_1,\\
&J_0[\mathcal{P}R^N(\nabla_{2,j}u_2,\nabla_{2,k}u_2)\nabla_{2,i}u_2-R^N(\nabla_{1,j}u_1,\nabla_{1,k}u_1),\nabla_{1,i}u_1]\rangle\,dM.
\endaligned
\]

Let $\psi:=\mathcal{P}\nabla_2u_2-\nabla_1u_1$ and $B_{\mathbf{i}}:=\mathcal{P}\nabla_{2,\mathbf{i}}-\nabla_{1,\mathbf{i}}\mathcal{P}$, where $0\leqslant\mathbf{i}\leqslant m$ and $x^0:=t$. $\Delta_{\lambda}$ denotes the Laplacian operator on $u_{\lambda}^*TN$ or $u_{\lambda}^*TN\otimes T^*M$. Then we have
\[
\aligned
&\frac{d}{dt}\Big(\frac{1}{2}\int_Mf|\psi|^2\,dM\Big)\\
=&\frac{1}{2}\int_Mf_t|\psi|^2\,dM-\int_Mg^{kl}\langle B_0\nabla_{2,k}u_2,f\psi_l\rangle\,dM+\int_Mg^{kl}g^{ij}\nabla_kf\cdot f\langle\psi_l,J_0B_i\nabla_{2,j}u_2\rangle\,dM\\
&-\int_Mg^{kl}g^{ij}\nabla_kf\cdot f\langle\psi_j,J_0B_i\nabla_{2,l}u_2\rangle\,dM-2\int_Mg^{kl}g^{ij}\nabla_if\cdot f\langle J_0\psi_k,B_j\nabla_{2,l}u_2\rangle\,dM\\
&-\int_Mf^2\cdot g^{ij}g^{kl}\langle\psi_l,J_0[\mathcal{P}R^N(\nabla_{2,j}u_2,\nabla_{2,k}u_2)\nabla_{2,i}u_2-R^N(\nabla_{1,j}u_1,\nabla_{1,k}u_1),\nabla_{1,i}u_1]\rangle\,dM\\
&+\int_Mf^2\cdot g^{kl}\langle(\Delta_1\mathcal{P}-\mathcal{P}\Delta_2)\nabla_{2,l}u_2,J_0\psi_k\rangle\,dM.
\endaligned
\]
By (3.5) of \cite{SW} and Lemma 3.2 of \cite{SW},
\[
|B_{\mathbf{i}}|\leqslant C_5\cdot\Big(|\nabla_{1,\mathbf{i}}u_1|+|\nabla_{2,\mathbf{i}}u_2|\Big)\cdot d(u_1,u_2),
\]
where $C_5$ depends on the geometry of $N$ and the derivative of exponential map in some domain(the details of this domain can be seen \cite{SW}). By (3.6) of \cite{SW} and Lemma 3.2 of \cite{SW},
\[
\aligned
&|(\Delta_1\mathcal{P}-\mathcal{P}\Delta_2)\nabla_2u_2|\\
\leqslant&C_6\Big\{d(u_1,u_2)+|\psi|+d(u_1,u_2)\Big[|\tau(u_1)|+|\tau(u_2)|\\
&+g^{ij}(|\nabla_{1,i}u_1|+|\nabla_{2,i}u_2|)(|\nabla_{1,j}u_1|+|\nabla_{2,j}u_2|)\Big]\Big\}\\
\leqslant&C_6\Big\{d(u_1,u_2)+|\psi|+d(u_1,u_2)\Big[|\tau(u_1)|+|\tau(u_2)|+m(|\nabla_1u_1|+|\nabla_2u_2|)^2\Big]\Big\}
\endaligned
\]
where $C_6$ depends only on $|\nabla_{\lambda}u_{\lambda}|$, $N$ and the derivative of exponential map in some domain(the details of this domain can be seen in \cite{SW}). Note that
\[
\aligned
|B_0|\leqslant&C_5\cdot(|\nabla_{1,t}u_1|+|\nabla_{2,t}u_2|)\cdot d(u_1,u_2)=C_5\cdot(|\tau_f(u_1)|+|\tau_f(u_2)|)\cdot d(u_1,u_2)\\
\leqslant&C_5\cdot(f|\tau(u_1)|+f|\tau(u_2)|+|\nabla f|\cdot|\nabla_1u_1|+|\nabla f|\cdot|\nabla_2u_2|)\cdot d(u_1,u_2).
\endaligned
\]
Let $\tilde{B}=\sum\limits_{i=1}^mB_idx^i$, then
\[|\tilde{B}|\leqslant C_5\cdot(|\nabla_1u_1|+|\nabla_2u_2|)\cdot d(u_1,u_2).\]
By the assumption of bounded geometry of $N$, it is easy to see that
\[
|g^{ij}\mathcal{P}R^N(\nabla_{2,j}u_2,\nabla_2u_2)\nabla_{2,i}u_2-g^{ij}R^N(\nabla_{1,j}u_1,\nabla_1u_1)\nabla_{1,i}u_1|\leqslant C_7\cdot(d(u_1,u_2)+|\psi|),
\]
where $C_7$ depends on $R^N$, $\nabla^NR^N$ and the $L^{\infty}$-norm of $\nabla_{\lambda}u_{\lambda}$. Thus
\[
\aligned
\frac{d}{dt}\Big(\frac{1}{2}\int_Mf|\psi|^2\,dM\Big)\leqslant&\frac{1}{2}\int_M|f_t|\cdot|\psi|^2\,dM+\int_M|B_0|\cdot f\cdot|\nabla_2u_2|\cdot|\psi|\,dM\\
&+4\int_M|\nabla f|\cdot f\cdot|\psi|\cdot|\tilde{B}|\cdot|\nabla_2u_2|\,dM\\
&+\int_Mf^2|(\Delta_1\mathcal{P}-\mathcal{P}\Delta_2)\nabla_2u_2|\cdot|\psi|\,dM\\
&+C_7\int_Mf^2\cdot|\psi|\cdot[d(u_1,u_2)+|\psi|]\,dM\\
\leqslant&C_5\int_M|\psi|\cdot d(u_1,u_2)\cdot|\nabla_2u_2|\Big[f^2(|\tau(u_1)|+|\tau(u_2)|)\\
&+f\cdot|\nabla f|\cdot(|\nabla_1u_1|+|\nabla_2u_2|)\Big]\,dM\\
&+4C_5\int_M|\psi|\cdot d(u_1,u_2)\cdot f|\nabla f|\cdot|\nabla_2u_2|(|\nabla_1u_1|+|\nabla_2u_2|)\,dM\\
&+\frac{1}{2}C_8\int_Mf|\psi|^2\,dM+(C_6+C_7)\int_Mf^2|\psi|^2\,dM\\
&+\int_Mf^2\cdot|\psi|\cdot d(u_1,u_2)\Big\{C_6+C_7\\
&+C_6\Big[|\tau(u_1)|+|\tau(u_2)|+m(|\nabla_1u_1|+|\nabla_2u_2|)^2\Big]\Big\}\,dM
\endaligned
\]
where $C_8:=\frac{1}{\delta}\cdot\sup\limits_{(t,x)\in[0,T]\times M}\{|f_t(t,x)|\}$. Since $u_{\lambda}\in\dot{W}^{1,\infty}(M,N)$, we have
\[
\aligned
\frac{d}{dt}\Big(\frac{1}{2}\int_Mf|\psi|^2\,dM\Big)\leqslant&C_9\int_M\sqrt{f}\cdot|\psi|\cdot d(u_1,u_2)\Big(|\tau(u_1)|+|\tau(u_2)|+|\nabla f|\Big)\,dM\\
&+(\frac{1}{2}C_8+C_9)\int_Mf|\psi|^2\,dM+C_9\int_M|\psi|\cdot d(u_1,u_2)\sqrt{f}\cdot|\nabla f|\,dM\\
&+C_9\int_M\sqrt{f}\cdot|\psi|\cdot d(u_1,u_2)\Big\{1+|\tau(u_1)|+|\tau(u_2)|\Big\}\,dM\\
\leqslant&\frac{1}{2}C_8\int_Mf\cdot|\psi|^2\,dM+C_9\Big[\frac{1}{2}\int_Mf\cdot|\psi|^2\,dM\\
&+\frac{1}{2}||d(u_1,u_2)(\sqrt{m}|\nabla^2_1u_1|+\sqrt{m}|\nabla^2_2u_2|+|\nabla f|)||_2^2\Big]\\
&+C_9\Big[\frac{1}{2}\int_Mf\cdot|\psi|^2\,dM+\frac{1}{2}\Big|\Big|d(u_1,u_2)\cdot|\nabla f|\Big|\Big|^2_2\Big]\\
&+C_9\Big[\frac{1}{2}||d(u_1,u_2)(1+\sqrt{m}|\nabla^2_1u_1|+\sqrt{m}|\nabla^2_2u_2|)||^2_2\\
&+\frac{1}{2}\int_Mf\cdot|\psi|^2\,dM\Big]
\endaligned
\]
where $C_9$ depends only on $||\nabla_{\lambda}u_{\lambda}||_{\infty}$, $||f||_{\infty}$, $C_5$, $C_6$ and $C_7$.

Since $\nabla^2_{\lambda}u_{\lambda}\in L^{\infty}$ and $\nabla f\in L^{\infty}$,
\[\frac{d}{dt}\Big(\frac{1}{2}\int_Mf|\psi|^2\,dM\Big)\leqslant C_{10}\Big[\frac{1}{2}\int_Mf|\psi|^2\,dM+||d(u_1,u_2)||_2^2\Big]\]
where $C_{10}$ depends only on $C_8$, $C_9$, $||\nabla^2_{\lambda}u_{\lambda}||_{\infty}$ and $||\nabla f||_{\infty}$.
In conclusion, we obtain
\begin{equation}\label{Formula4'}
\frac{d}{dt}\Big(\frac{1}{2}\int_Mf|\psi|^2\,dM\Big)\leqslant C\Big[\frac{1}{2}\int_Mf|\psi|^2\,dM+||d(u_1,u_2)||_2^2\Big].
\end{equation}
Combining $(\ref{Formula1})$ and $(\ref{Formula4'})$, we get
\[\frac{d}{dt}\Big(\frac{1}{2}\int_Mf|\psi|^2\,dM+||d(u_1,u_2)||_2^2\Big)\leqslant C\cdot\Big(\frac{1}{2}\int_Mf|\psi|^2\,dM+||d(u_1,u_2)||_2^2\Big),\]
where $C$ depends on the norm of $u_{\lambda}$ in $\mathscr{S}_{\infty}$. Since $d(u_1(0),u_2(0))=0$ and $\frac{1}{2}\int_Mf(0)|\psi(0)|^2\,dM=0$,
\[\frac{1}{2}\int_Mf|\psi|^2\,dM+||d(u_1,u_2)||_2^2=0\s in\s [0,T'].\]
So $u_1=u_2$ a.e. in $[0,T']\times M$. Since uniqueness is a local property, this completes the proof.
\endproof
\\

\textbf{Proof of Theorem \ref{thm0.2}}\\
Using the same method by which we prove Theorem \ref{thm0.1}, one can obtain
\begin{equation}\label{Formula1'}
\aligned
&\frac{d}{dt}\int_Md^2(u_1,u_2)\,dM\\
\leqslant&2\eta\cdot C_5\int_Md^2(u_1,u_2)\,dM+2\eta\int_M|\mathcal{P}\nabla_2u_2-\nabla_1u_1|^2\,dM\\
\leqslant&2\eta\cdot C_5\int_Md^2(u_1,u_2)\,dM+\frac{2\eta}{\delta}\int_Mf|\mathcal{P}\nabla_2u_2-\nabla_1u_1|^2\,dM,
\endaligned
\end{equation}
and
\[
\aligned
\frac{d}{dt}\Big(\frac{1}{2}\int_Mf|\psi|^2\,dM\Big)\leqslant&\frac{1}{2}C_8\int_Mf\cdot|\psi|^2\,dM+C_9\Big[\frac{1}{2}\int_Mf\cdot|\psi|^2\,dM\\
&+\frac{1}{2}||d(u_1,u_2)(\sqrt{m}|\nabla^2_1u_1|+\sqrt{m}|\nabla^2_2u_2|+|\nabla f|)||_2^2\Big]\\
&+C_9\Big[\frac{1}{2}\int_Mf\cdot|\psi|^2\,dM+\frac{1}{2}\Big|\Big|d(u_1,u_2)\cdot|\nabla f|\Big|\Big|^2_2\Big]\\
&+C_9\Big[\frac{1}{2}||d(u_1,u_2)(1+\sqrt{m}|\nabla^2_1u_1|+\sqrt{m}|\nabla^2_2u_2|)||^2_2\\
&+\frac{1}{2}\int_Mf\cdot|\psi|^2\,dM\Big]
\endaligned
\]
where $C_8:=\frac{1}{\delta}\cdot\sup\limits_{(t,x)\in[0,T]\times M}\{|f_t(t,x)|\}$ and $C_9$ depends only on $||\nabla_{\lambda}u_{\lambda}||_{\infty}$, $||f||_{\infty}$, $C_5$, $C_6$ and $C_7$.

Since $\nabla^2_{\lambda}u_{\lambda}\in L^{\infty}([0,T],L^m(M,N))$ and $\nabla f\in L^{\infty}([0,T],L^m(M))$, we have
\[
\aligned
&||d(u_1,u_2)(\sqrt{m}|\nabla^2_1u_1|+\sqrt{m}|\nabla^2_2u_2|+|\nabla f|)||_2\\
\leqslant&||d(u_1,u_2)||_{\frac{2m}{m-2}}\Big(\sqrt{m}||\nabla_1^2u_1||_m+\sqrt{m}||\nabla_2^2u_2||_m+||\nabla f||_m\Big),
\endaligned
\]
and
\[\Big|\Big|d(u_1,u_2)\cdot|\nabla f|\Big|\Big|_2\leqslant||d(u_1,u_2)||_{\frac{2m}{m-2}}\cdot||\nabla f||_m,\]
and
\[
\aligned
&||d(u_1,u_2)(1+\sqrt{m}|\nabla^2_1u_1|+\sqrt{m}|\nabla^2_2u_2|)||_2\\
\leqslant&||d(u_1,u_2)||_2+||d(u_1,u_2)||_{\frac{2m}{m-2}}\cdot\sqrt{m}(||\nabla_1^2u_1||_m+||\nabla_2^2u_1||_m).
\endaligned
\]
By Sobolev embedding and the estimate in Lemma 2.2 of \cite{SW}, we have
\[
\aligned
||d(u_1,u_2)||_{\frac{2m}{m-2}}\leqslant&C_{11}||d(u_1,u_2)||_{W^{1,2}}\\
\leqslant&C_{12}\cdot(||d(u_1,u_2)||_2+||\psi||_2)\leqslant C_{13}\cdot\Big(||d(u_1,u_2)||_2+\sqrt{\frac{1}{2}\int_Mf|\psi|^2\,dM}\Big)
\endaligned
\]
In conclusion, we obtain
\begin{equation}\label{Formula4}
\frac{d}{dt}\Big(\frac{1}{2}\int_Mf|\psi|^2\,dM\Big)\leqslant C\Big[\frac{1}{2}\int_Mf|\psi|^2\,dM+||d(u_1,u_2)||_2^2\Big].
\end{equation}
Combining $(\ref{Formula1'})$ and $(\ref{Formula4})$, we get
\[\frac{d}{dt}\Big(\frac{1}{2}\int_Mf|\psi|^2\,dM+||d(u_1,u_2)||_2^2\Big)\leqslant C\cdot\Big(\frac{1}{2}\int_Mf|\psi|^2\,dM+||d(u_1,u_2)||_2^2\Big),\]
where $C$ depends on the norm of $u_{\lambda}$ in $\mathscr{S}_m$. Since $d(u_1(0),u_2(0))=0$ and $\frac{1}{2}\int_Mf(0)|\psi(0)|^2\,dM=0$,
\[\frac{1}{2}\int_Mf|\psi|^2\,dM+||d(u_1,u_2)||_2^2=0\s in\s [0,T'].\]
So $u_1=u_2$ a.e. in $[0,T']\times M$. Since uniqueness is a local property, this completes the proof.\endproof
\section{Appendix}
In this section, we always assume that $N$ is compact. Otherwise, we replace $N$ by $\bar{\Omega}$, where $\Omega:=\{y\in N|dist_N(y,u_0(M))<1\}$. By Proposition 8.3 in Chapter 15 of \cite{T} we can see that (\ref{T:1}) has a local $C^{\infty}$-smooth solution. Since $f>\delta>0$, it is easy to see that  (\ref{T:1}) is parabolic strongly in the sense of Petrowski. We need to show that the solution to Cauchy problem (\ref{T:1}) lies on $N$ if the initial value map $u_0(M)\subseteq N$, as well as the Cauchy problem  (\ref{T:1}) admits local solutions. For this goal, we first define a tubular neighborhood of $N$, say $N_d\subseteq\mathbb{R}^L$, $d>0$, as follows
\begin{eqnarray*}
N_d:=\{y\in \mathbb{R}^L|\mbox{dist}(y,N)<d\}
\end{eqnarray*}
such that for any $y\in N_d$, there exists a unique point $q=\pi(y)\in N$ satisfying $\mbox{dist}(y,N)=|y-q|$. By choosing $d$ small, the projection map $\pi:N_d\longrightarrow N$ is smooth.
\medskip

Define $\rho:N_d\longrightarrow \mathbb{R}^L$ by $\rho(y)=y-\pi(y)$ for $y\in N_d$. Obviously,
\begin{eqnarray*}
\rho+\pi=Id
\end{eqnarray*}
is an identity map from $N_d\subseteq \mathbb{R}^L$ into $N_d\subseteq \mathbb{R}^L$. Since $\pi^2=\pi$, we know $\rho\circ\pi=0$. Hence, it follows that there holds true
\begin{eqnarray*}
d\rho\circ d\pi=0
\end{eqnarray*}
which implies that for any $y,y'\in\mathbb{R}^L$ and any operator $\mathcal{B}$ from $d\pi(\mathbb{R}^L)$ into $d\pi(\mathbb{R}^L)$
\begin{eqnarray*}
\langle d\rho(y),\mathcal{B}\circ d\pi(y')\rangle=0.
\end{eqnarray*}
This is equivalent to
\begin{eqnarray*}
d\rho\circ\mathcal{B}\circ d\pi=0.
\end{eqnarray*}
It is easy to know that for $y\in N$, $d\pi(y)$ is just the orthogonal projection operator $P(y)$ which maps $\mathbb{R}^L$ into $T_yN$ and $d\rho(y)$ is just $Id-P(y)=B(y)$.

Meanwhile, for a map $u:M\longrightarrow N\hookrightarrow\mathbb{R}^L$, take $A(u):T_uN\times T_uN\longrightarrow(T_uN)^{\bot}$ as the second fundamental form, then by the properties of submanifold we have
\begin{eqnarray*}
\tau(u)=\nabla_i\nabla_iu=\Delta u+A(u)(\nabla u,\nabla u).
\end{eqnarray*}
We consider the following Cauchy problem which is equivalent to the corresponding Cauchy problem of (\ref{T:1}):
\begin{eqnarray}\label{A:1}
\left\{
\begin{array}{llll}
\aligned
u_t=&\varepsilon f\Delta u+\varepsilon f\lambda(|\rho(u)|^2)A(\pi\circ u)(\nabla(\pi\circ u),\nabla(\pi\circ u))\\
&+\varepsilon\cdot\lambda(|\rho(u)|^2)\cdot df\cdot\nabla(\pi\circ u)+\lambda(|\rho(u)|^2)V(\pi\circ u),
\endaligned
\\
u(\cdot ,0)=u_0.
\end{array}
\right.
\end{eqnarray}
Here $\lambda(s)$ is a smooth function on $[0,\,d]$ such that $0\leqslant\lambda(s)\leqslant1$,
\begin{eqnarray*}
\left\{
\begin{array}{llll}
\lambda(s^2)=1,\s\s if\,\,\,s\in[0,\frac{d}{2});\\
\lambda(s^2)=0, \s\s if\,\,\,s>\frac{3}{4}d.
\end{array}
\right.
\end{eqnarray*}
and $V(\pi\circ u):=J(\pi\circ u)\tau_f(\pi\circ u)$. In $\mathbb{R}^L$, we write
\begin{eqnarray*}
\aligned
\tau_f(\pi\circ u)&=f\tau(\pi\circ u)+df\cdot\nabla(\pi\circ u)\\
&=f[d\pi(\Delta u)+\nabla d\pi(Du,Du)+A(\pi\circ u)(\nabla(\pi\circ u),\nabla(\pi\circ u))]+df\cdot\nabla(\pi\circ u)
\endaligned
\end{eqnarray*}
where $Du$ is the covariant differential induced by the map $u:M\longrightarrow\mathbb{R}^L$. Then, the term of second order derivatives in (\ref{A:1}) is of the following form
\begin{eqnarray*}
\varepsilon f\Delta u+f\cdot\lambda(|\rho\circ u|^2)J(\pi\circ u)d\pi(\Delta u).
\end{eqnarray*}
Hence, the principle symbol of (\ref{A:1}) can be written by
\begin{eqnarray*}
\aligned
&f[\varepsilon|\tilde{\eta}|^2+\langle\tilde{\eta},\lambda(|\rho\circ u|^2)J(\pi\circ u)\circ d\pi(\tilde{\eta})\rangle]\\
=&f\cdot\varepsilon|\tilde{\eta}|^2+f\lambda(|\rho\circ u|^2)\langle(d\rho+d\pi)\tilde{\eta},J(\pi\circ u)\circ d\pi(\tilde{\eta})\rangle\\
=&f\cdot\varepsilon|\tilde{\eta}|^2>\delta\varepsilon|\tilde{\eta}|^2.
\endaligned
\end{eqnarray*}
Here, we have used $d\rho+d\pi=Id$ and the fact
\begin{eqnarray*}
\langle d\rho(\tilde{\eta}),J(\pi\circ u)\circ d\pi(\tilde{\eta})\rangle=0
\end{eqnarray*}
which follows $\langle d\rho(\tilde{\eta}),d\pi(\tilde{\eta}')\rangle\equiv0$ and $d\rho\circ J(\pi\circ u)\circ d\pi\equiv0$.

Indeed, we have
\begin{eqnarray*}
\aligned
&\langle d\rho(\tilde{\eta}),J(\pi\circ u)\circ d\pi(\tilde{\eta})\rangle\\
=&\langle d\rho(\tilde{\eta}),(d\rho+d\pi)\circ J(\pi\circ u)\circ d\pi(\tilde{\eta})\rangle\\
=&\langle d\rho(\tilde{\eta}),d\rho\circ J(\pi\circ u)\circ d\pi(\tilde{\eta})\rangle+\langle d\rho(\tilde{\eta}),d\pi\circ J(\pi\circ u)\circ d\pi(\tilde{\eta})\rangle\\
=&0.
\endaligned
\end{eqnarray*}
So, the following operator
\begin{eqnarray*}
\sum\limits_{j,k}A^{jk}(t,x,u)\partial_j\partial_ku=f[\varepsilon Id+\lambda(|\rho\circ u|^2)J(\pi\circ u)d\pi]\Delta
\end{eqnarray*}
satisfies the so-called strong parabolicity hypothesis:
\begin{eqnarray*}
\sum\limits_{j,k}A^{jk}(u)\xi_j\xi_k\geqslant\delta\varepsilon|\xi|^2Id.
\end{eqnarray*}
This means that (\ref{A:1}) is a strong parabolic equation in the sense of Petrowski. Therefore, according to Proposition 8.3 in chapter 15 of \cite{T} we know that, if $u_0\in H^s(M)$, $s > \frac{\dim(M)}{2}$(and $s\ge 1$ if $\dim(M)=1$), then the Cauchy problem (\ref{A:1}) admits a unique solution
\begin{eqnarray*}
u_{\varepsilon}\in C([0,\, T_{\varepsilon}), H^s(M,\mathbb{R}^L))\cap C^\infty((0,\, T_{\varepsilon})\times M).
\end{eqnarray*}
Therefore, when $t$ is small enough, for $s$ is large enough, and $0<\varepsilon<1$ we can see easily that
\begin{eqnarray*}
\sup\limits_{x\in M}|u_{\varepsilon}(x,t)-u_0(x)|<d.
\end{eqnarray*}

It remains that we need to verify the solution $u_{\varepsilon}$ will always be $N$-valued if the initial map $u_0$ is $N$-valued. We consider $|\rho\circ u_{\varepsilon}(x,t)|=\mbox{dist}(u_{\varepsilon}(x,t),N)$. Obviously,
\begin{eqnarray*}
\rho\circ u_{\varepsilon}\in(T_{\pi\circ u_{\varepsilon}}N)^{\bot}.
\end{eqnarray*}
Hence, to show that $u_{\varepsilon}$ is $N$-valued, it suffices to show that, for all $t\in[0,T_{\varepsilon})$,
$$\int_Mf^{-1}|\rho\circ u_{\varepsilon}(t)|^2\,dM=0.$$
Indeed, noting $d\rho+d\pi=id$ on $\mathbb{R}^L$ where $T_{\pi\circ u_{\varepsilon}}N+(T_{\pi\circ u_{\varepsilon}}N)^{\bot}$ is regarded as $\mathbb{R}^L$, we could deduce
\begin{eqnarray}\label{A:2}
\aligned
&\frac{1}{2}\frac{d}{dt}\int_Mf^{-1}|\rho\circ u_{\varepsilon}|^2\,dM\\
=&-\frac{1}{2}\int_Mf^{-2}f_t|\rho\circ u_{\varepsilon}|^2\,dM+\int_Mf^{-1}\langle\rho\circ u_{\varepsilon},d\rho_{u_{\varepsilon}}(u_{\varepsilon t})\rangle\,dM\\
=&-\frac{1}{2}\int_Mf^{-2}f_t|\rho\circ u_{\varepsilon}|^2\,dM+\int_Mf^{-1}\langle\rho\circ u_{\varepsilon},u_{\varepsilon t}\rangle\,dM\\
&-\int_Mf^{-1}\langle\rho\circ u_{\varepsilon},d\pi_{u_{\varepsilon}}(u_{\varepsilon t})\rangle\,dM.
\endaligned
\end{eqnarray}
The last term on the right hand side of (\ref{A:2}) vanishes since $\rho\circ u_{\varepsilon}\in(T_{\pi\circ u_{\varepsilon}}N)^{\bot}$ and $d\pi_{u_{\varepsilon}}(u_{\varepsilon t})\in T_{\pi\circ u_{\varepsilon}}N$.
\medskip

Substituting (\ref{A:1}) into (\ref{A:2}) yields
\begin{eqnarray}\label{A:3}
\aligned
&\frac{1}{2}\frac{d}{dt}\int_Mf^{-1}|\rho\circ u_{\varepsilon}|^2\,dM=-\frac{1}{2}\int_Mf^{-2}f_t|\rho\circ u_{\varepsilon}|^2\,dM\\
&+\int_M\langle\rho\circ u_{\varepsilon},\,\varepsilon\Delta u_{\varepsilon}+\varepsilon\lambda(|\rho\circ u_{\varepsilon}|^2)A(\pi\circ u_{\varepsilon})(\nabla(\pi\circ u_{\varepsilon}),\,\nabla(\pi\circ u_{\varepsilon}))\rangle\,dM\\
&+\int_Mf^{-1}\langle\rho\circ u_{\varepsilon},\,\varepsilon\cdot\lambda(|\rho\circ u_{\varepsilon}|^2)df\cdot\nabla(\pi\circ u_{\varepsilon})\rangle\,dM\\
&+\int_Mf^{-1}\langle\rho\circ u_{\varepsilon},\,\lambda(|\rho\circ u_{\varepsilon}|^2)V(\pi\circ u_{\varepsilon})\rangle\,dM.
\endaligned
\end{eqnarray}
The third term on the right hand side of (\ref{A:3}) vanishes since $\rho\circ u_{\varepsilon}\in(T_{\pi\circ u_{\varepsilon}}N)^{\bot}$ and $df\cdot\nabla(\pi\circ u_{\varepsilon})\in T_{\pi\circ u_{\varepsilon}}N$.
\medskip

Using the relation $\rho+\pi=id$, we get
\begin{eqnarray}\label{A:4}
\aligned
&\frac{1}{2}\frac{d}{dt}\int_Mf^{-1}|\rho\circ u_{\varepsilon}|^2\,dM=-\frac{1}{2}\int_Mf^{-2}f_t|\rho\circ u_{\varepsilon}|^2\,dM\\
&+\varepsilon\int_M\langle\rho\circ u_{\varepsilon},\Delta(\pi\circ u_{\varepsilon})+\lambda(|\rho\circ u_{\varepsilon}|^2)A(\pi\circ u_{\varepsilon})(\nabla(\pi\circ u_{\varepsilon}),\nabla(\pi\circ u_{\varepsilon}))\rangle\,dM\\
&+\varepsilon\int_M\langle\rho\circ u_{\varepsilon},\Delta(\rho\circ u_{\varepsilon})\rangle\,dM+\int_Mf^{-1}\langle\rho\circ u_{\varepsilon},\lambda(|\rho\circ u_{\varepsilon}|^2)V(\pi\circ u_{\varepsilon})\rangle\,dM.
\endaligned
\end{eqnarray}
Since $\lambda(|\rho\circ u_{\varepsilon}|^2)=1$, by the definitions we know that, when $t$ is small enough, there hold true
\begin{eqnarray*}
\Delta(\pi\circ u_{\varepsilon})+\lambda(|\rho\circ u_{\varepsilon}|^2)A(\pi\circ u_{\varepsilon})(\nabla(\pi\circ u_{\varepsilon}),\nabla(\pi\circ u_{\varepsilon}))=\tau(\pi\circ u_{\varepsilon})\in T_{\pi\circ u_{\varepsilon}}N
\end{eqnarray*}
and
\begin{eqnarray*}
V(\pi\circ u_{\varepsilon})\in T_{\pi\circ u_{\varepsilon}}N.
\end{eqnarray*}
Hence the second term and the fourth term of the right hand side of (\ref{A:4}) vanish. Therefore, by integrating by parts we obtain from (\ref{A:4})

\begin{eqnarray}\label{A:5}
\aligned
&\frac{1}{2}\frac{d}{dt}\int_Mf^{-1}|\rho\circ u_{\varepsilon}|^2\,dM\\
=&-\frac{1}{2}\int_Mf^{-2}f_t|\rho\circ u_{\varepsilon}|^2\,dM+\varepsilon\int_M\langle\rho\circ u_{\varepsilon},\Delta(\rho\circ u_{\varepsilon})\rangle\,dM\\
=&-\frac{1}{2}\int_Mf^{-2}f_t|\rho\circ u_{\varepsilon}|^2\,dM-\varepsilon\int_M|\nabla(\rho\circ u_{\varepsilon})|^2\,dM.
\endaligned
\end{eqnarray}
Recall that $$C_1(t)=\max\limits_{x\in M}\{|f_t(x,t)|/f(x,t)\}$$
and let
$$E^N_f(u_{\varepsilon}):=\frac{1}{2}\int_Mf^{-1}|\rho\circ u_{\varepsilon}|^2\,dM.$$
Since $\rho(u_{\varepsilon})(0)=\rho(u_0)=0$, we have
\begin{eqnarray*}
\left\{
\begin{array}{llll}
\frac{d}{dt}E^N_f(u_{\varepsilon})\leqslant C_1(t)E^N_f(u_{\varepsilon}),\\
E^N_f(u_{\varepsilon})(0)=0.
\end{array}
\right.
\end{eqnarray*}
It follows from the Gronwall inequality
\begin{eqnarray*}
E^N_f(u_{\varepsilon})(t)\leqslant0.
\end{eqnarray*}
By the definition of $E^N_f(u_{\varepsilon})(t)$, we know that it is not negative. So $E^N_f(u_{\varepsilon})(t)=0$ and it implies $(\rho\circ u_{\varepsilon})(t)=0$ for all $t\in[0,T_{\varepsilon})$. Thus $u_{\varepsilon}(t)$ is an $N$-valued map. This completes the proof.\endproof

\medskip

Zonglin Jia

{\small\it Academy of Mathematics and Systems Sciences, Chinese Academy of Sciences, Beijing 100190,  P.R. China.}

{\small\it Email: 756693084@qq.com}
\\

Youde Wang

{\small\it Academy of Mathematics and Systems Sciences, Chinese Academy of Sciences, Beijing 100190,  P.R. China.}

{\small\it Email: wyd@math.ac.cn}

\end{document}